\documentclass[aos,seceqn,nameyear,dvips]{arximspdf}
\usepackage{graphics}
\doi{10.1214/09-AOS767}
\volume{38}
\issue{3}
\pubyear{2010}
\firstpage{1793}
\lastpage{1844}

\makeatletter

\newtheorem{theorem}{Theorem}
\newtheorem{lemma}{Lemma}
\newtheorem{corollary}{Corollary}
\newproclaim{remark}{Remark}
\newproclaim{example}{Example}

\newproclaim{Condition}{Condition}

\mathchardef\varLambda="0103
\mathchardef\varXi="0104
\mathchardef\varOmega="010A
\mathchardef\varDelta="0101
\makeatother

\begin{document}
\begin{frontmatter}

\title{On convergence rates equivalency and sampling strategies in
functional deconvolution models}
\runtitle{Functional deconvolution}

\begin{aug}
\author[A]{\fnms{Marianna} \snm{Pensky}\corref{}\thanksref{t1}\ead[label=e1]{mpensky@pegasus.cc.ucf.edu}} and
\author[B]{\fnms{Theofanis} \snm{Sapatinas}\ead[label=e2]{T.Sapatinas@ucy.ac.cy}}
\runauthor{M. Pensky and T. Sapatinas}
\affiliation{University of Central Florida and University of Cyprus}
\address[A]{Department of Mathematics\\
University of Central Florida\\
Orlando, Florida 32816-1353\\
USA\\
\printead{e1}}
\address[B]{Department of Mathematics\\
\quad and Statistics\\
University of Cyprus\\
P.O. Box 20537\\
Nicosia CY 1678\\
Cyprus\\
\printead{e2}}
\end{aug}

\thankstext{t1}{Supported in part by NSF Grants DMS-05-05133 and DMS-06-52524.}

\received{\smonth{8} \syear{2009}}
\revised{\smonth{11} \syear{2009}}

%
\begin{abstract}
Using the asymptotical minimax framework, we examine convergence rates
equivalency between a continuous functional deconvolution model and
its real-life discrete counterpart over a wide range of Besov balls
and for the $L^2$-risk. For this purpose, all possible models are
divided into three groups. For the models in the first group, which we
call \textit{uniform}, the convergence rates in the discrete and the
continuous models coincide no matter what the sampling scheme is
chosen, and hence the replacement of the discrete model by its
continuous counterpart is legitimate. For the models in the second
group, to which we refer as \textit{regular}, one can point out the best
sampling strategy in the discrete model, but not every sampling scheme
leads to the same convergence rates; there are at least two sampling
schemes which deliver different convergence rates in the discrete model
(i.e., at least one of the discrete models leads to convergence rates
that are different from the convergence rates in the continuous model).
The third group consists of models for which, in general, it is
impossible to devise the best sampling strategy; we call these models
\textit{irregular}.

We formulate the conditions when each of these situations takes place.
In the regular case, we not only point out the number and the selection
of sampling points which deliver the fastest convergence rates in the
discrete model but also investigate when, in the case of an arbitrary
sampling scheme, the convergence rates in the continuous model coincide
or do not coincide with the convergence rates in the discrete model.
We also study what happens if one chooses a uniform, or a more general
pseudo-uniform, sampling scheme which can be viewed as an intuitive
replacement of the continuous model. Finally, as a representative of
the irregular case, we study functional deconvolution with a
boxcar-like blurring function since this model has a number of
important applications. All theoretical results presented in the paper
are illustrated by numerous examples; many of which are motivated
directly by a multitude of inverse problems in mathematical physics
where one needs to recover initial or boundary conditions on the basis
of observations from a noisy solution of a partial differential
equation. The theoretical performance of the suggested estimator in the
multichannel deconvolution model with a boxcar-like blurring function
is also supplemented by a limited simulation study and compared to an
estimator available in the current literature. The paper concludes that
in both regular and irregular cases one should be extremely careful
when replacing a discrete functional deconvolution model by its
continuous counterpart.
\end{abstract}

%

\end{frontmatter}

\section{Introduction}
\label{sec:intro}

We consider the estimation problem of the unknown response function
$f(\cdot)$ based on observations from the following noisy
convolutions:
%
%
\begin{equation}\label{convcont}
y(u, t) = f*g (u, t) + \frac{1}{\sqrt{n}}
z(u,t),\qquad u \in U, t \in T,
\end{equation}
where
$U=[a,b]$, $-\infty< a \leq b < \infty$ and $T=[0,1]$. Here, $z(u,
t)$ is assumed to be a two-dimensional Gaussian white noise, that is, a
generalized two-dimensional Gaussian field with covariance function
${\mathbb E}[z (u_1, t_1) z (u_2, t_2)] = \delta(u_1-u_2) \delta(t_1-t_2)$,
where $\delta(\cdot)$ denotes the Dirac $\delta$-function, and
\[
f*g(u, t) = \int_T f(x)g(u,t-x) \,dx
\]
with the
blurring (or kernel) function $g(\cdot,\cdot)$ also assumed to be
known.

The model (\ref{convcont}) has been recently introduced by Pensky and
Sapatinas (\citeyear{PenskyS09}) and can be viewed as a \textit{functional
deconvolution} model. If $a=b$, it reduces to the standard
deconvolution model which attracted the attention of a number of
researchers, for example, Donoho (\citeyear{D95}), Abramovich and
Silverman (\citeyear{AS98}),
Kalifa and Mallat (\citeyear{KM03}), Johnstone et al. (\citeyear
{Johnstoneetal04}), Donoho and Raimondo
(\citeyear{DR04}), Johnstone and Raimondo
(\citeyear{JR04}), Neelamani, Choi and Baraniuk (\citeyear{NCB04}),
Kerkyacharian, Picard and Raimondo
(\citeyear{KPR07}), Cavalier and Raimondo (\citeyear{CR07}) and
Chesneau (\citeyear{C08}).

The functional deconvolution model (\ref{convcont})
can be viewed as a generalization of a multitude of inverse problems
in mathematical physics where one needs to recover initial or
boundary conditions on the basis of observations of a noisy solution
of a partial differential equation. Lattes and Lions (\citeyear{LL67})
initiated
research in the problem of recovering the initial condition for
parabolic equations based on observations in a fixed-time strip,
while this problem and the problem of recovering the boundary
condition for elliptic equations based on observations in an
internal domain were studied in Golubev and Khasminskii (\citeyear
{GK99}); the
latter problem was also discussed in Golubev (\citeyear{G04}). These
and other
specific models in mathematical physics were discussed in detail in
Pensky and Sapatinas (\citeyear{PenskyS09}).

However, model (\ref{convcont}) is just an idealization
of a real-life situation. One can make observations only at particular points
$(u_l, t_i)$, $l=1,2,\ldots, M$, $i=1,2,\ldots, N$,
so that the actual problem can be formulated as follows: recover the
unknown response function $f(\cdot)$
from observations $y(u_l, t_i)$, where
%
%
\begin{equation}\label{convdis}
y(u_l, t_i) = \int_T f(x) g(u_l, t_i-x) \,dx + \varepsilon_{li},\qquad
u_l \in U, t_i = i/N,
\end{equation}
with $\varepsilon_{li}$ being standard Gaussian random variables, independent
for different $l$ and~$i$.
Model (\ref{convdis}) can be viewed as a \textit{discrete} version of the
\textit{continuous} functional deconvolution model (\ref{convcont}).

It is well documented in the literature that asymptotic equivalence
between discrete and continuous
models holds in some nonparametric models. In particular, Brown and Low
(\citeyear{BL96}) and Brown et al. (\citeyear{Brownetal02})
in the univariate case and Reiss (\citeyear{R08}) in the multivariate case
established, under some restrictions, asymptotic equivalence (in the Le Cam
sense) between nonparametric regression and Gaussian white noise
models. Although, to the best of our knowledge, such an asymptotic
equivalence between continuous and discrete models, in the functional
deconvolution setting, has not yet been explored, it has been
documented in the literature a convergence rate equivalency, in the
asymptotical minimax sense, between standard continuous and discrete
deconvolution models, that is, when $a=b$, $M=1$ and $N=n$ in (\ref
{convcont}) and (\ref{convdis}), over a wide range of Besov balls and for
the $L^r$-risks, $1 \leq r < \infty$ [e.g., Chesneau (\citeyear
{C08}), Pensky and
Sapatinas (\citeyear{PenskyS09}) and Petsa and Sapatinas (\citeyear{PetsaS09})].

For the above reason, and using the asymptotical minimax framework, one
may attempt to study the continuous functional deconvolution model
(\ref{convcont})
instead of its discrete counterpart (\ref{convdis}), assuming that the
convergence rates between these models coincide. However, in this case,
this equivalence has only a limited scope. Indeed, Pensky and Sapatinas
(\citeyear{PenskyS09}) only
touched upon the issue, showing that, under very restrictive
conditions, a convergence rate equivalence between the continuous
functional deconvolution model (\ref{convcont})
and its discrete counterpart (\ref{convdis}) models holds when $n = N M$,
over a wide range of Besov balls and for the $L^2$-risk. Nevertheless,
in majority
of practical situations, these conditions are violated and it remains
to be seen how legitimate the replacement of the real life model (\ref
{convdis})
by its idealization (\ref{convcont}) is, even in the case of inverse
problems in mathematical physics, presented in Pensky and Sapatinas
(\citeyear{PenskyS09}).
In fact, in many situations, the convergence rates in the two models
depend on the choice of $M$ and the selection of sampling points $u_1,
u_2, \ldots, u_M$
and may coincide with the convergence rates in the continuous model for
one selection
and be different for another. Also, from a practical point of view, the
objective
is not to find $M$ and $u_1,u_2, \ldots, u_M$ which make the two models
equivalent from the convergence rate viewpoint, but rather
to point out $M$ and $u_1,u_2, \ldots, u_M$ which deliver the fastest
possible convergence rates in the real life model (\ref{convdis}).
Note that the discrete model (\ref{convdis}) can also be viewed as a
multichannel deconvolution
model where the number of channels $M=M_n$ is fixed or, possibly, $M_n
\rightarrow
\infty$ as the sample size $n \rightarrow\infty$; the case when $M
\geq2$ (finite) was considered in, for example, Casey and Walnut
(\citeyear{CW94}) and
De Canditiis and Pensky (\citeyear{CP04}, \citeyear{CP06}). Hence if
the kernel $g(\cdot,\cdot
)$ is fixed,
the choice of $M$ and the selection of sampling points $u_1,u_2, \ldots
, u_M$ which provide the fastest convergence rates is of extreme
importance in signal processing.

Using the asymptotical minimax framework, our objective is to evaluate
how legitimate
it is to replace the real-life discrete model (\ref{convdis}) by its
continuous counterpart (\ref{convcont}). For this purpose, we shall
divide all possible models into three groups. For the models in the
first group, which we call \textit{uniform}, the convergence rates
in discrete and continuous functional deconvolution models coincide no
matter what the sampling scheme is chosen, and hence the replacement of
the discrete model by its continuous counterpart is legitimate.
For the models in the second group, to which we refer as \textit{regular},
one can point out the best
sampling strategy in the discrete model (i.e., the strategy which leads
to the fastest convergence rate),
but not every sampling scheme leads to the same convergence rates:
there are at least two sampling schemes which deliver
different convergence rates in the discrete model (i.e., at least one
of the discrete models leads to convergence rates that are different
from the convergence rates in the continuous functional deconvolution
model). The third group consists of models
for which, in general, it is impossible to devise the best sampling strategy;
we call these models \textit{irregular}.

We formulate the conditions when each of these situations takes place.
In the regular case, we not only point out the choice of $M$ and the
selection of sampling points $u_1,u_2, \ldots, u_M$
which deliver the fastest convergence rates in the discrete model but also
investigate when, in the case of an arbitrary sampling scheme,
the convergence rates in the continuous functional deconvolution model
coincide or do not coincide with the convergence rates of its discrete
counterpart.
We also study what happens if one chooses a uniform, or a more general
pseudo-uniform, sampling scheme
which can be viewed as an intuitive replacement of the continuous model.
Finally, as a representative of the irregular case, we study functional
deconvolution with a boxcar-like kernel since this model has a number
of important applications.
All theoretical results presented are illustrated by numerous examples;
many of which are motivated directly by a multitude of inverse problems
in mathematical physics where
one needs to recover initial or boundary conditions on the basis of observations
from a noisy solution of a partial differential equation.


%

As in Pensky and Sapatinas (\citeyear{PenskyS09}), we consider functional
deconvolution in a periodic setting; that is, we assume that $f(\cdot)$
and, for fixed
$u \in U$, $g(u, \cdot)$ are periodic functions with
period on the unit interval $T$. Note that the periodicity
assumption appears naturally in the above mentioned special models
which (\ref{convcont}) and (\ref{convdis}) generalize and allows one to
explore ideas considered in the above cited papers to the proposed
functional deconvolution framework. Moreover, not only for
theoretical reasons but also for practical convenience [see
Johnstone et al. (\citeyear{Johnstoneetal04}), Sections 2.3, 3.1 and 3.2],
we use band-limited wavelet bases and in particular the
periodized Meyer wavelet basis for which fast algorithms exist [see
Kolaczyk (\citeyear{K94}) and Donoho and Raimondo (\citeyear{DR04})].
In order to also allow
inhomogeneous functions $f(\cdot)$ into our study, we consider a wide
range of Besov balls, as it is common in the wavelet literature, and, for
simplicity, we work with the $L^2$-risk only. However, the results
of this paper can be extended to a more general class of $L^r$-risks,
$1 \leq r < \infty$, using the unconditionality and Temlyakov
properties of Meyer wavelets
[e.g., Johnstone et al. (\citeyear{Johnstoneetal04}) and Petsa and
Sapatinas (\citeyear{PetsaS09})].

The rest of the paper is organized as follows. In Section \ref
{estconstr}, we
describe the construction of wavelet estimators of $f(\cdot)$
derived by Pensky and Sapatinas (\citeyear{PenskyS09}) both for the continuous
functional deconvolution model
(\ref{convcont}) and its discrete counterpart (\ref{convdis}). In
Section \ref
{risk_lower_upper},
for the continuous model and for a discrete model with any particular
design of sampling points,
using the asymptotical minimax framework, we provide lower bounds for
the $L^2$-risk over a wide range
of Besov balls and show that those bounds
are attained by the wavelet estimators constructed in Section \ref{estconstr}.
Section \ref{discrcont} is devoted to the discussion of the interplay
between continuous and discrete functional deconvolution models. First,
Section \ref{suf_uniform} formulates the necessary and
sufficient conditions for the
convergence rates in continuous and discrete functional deconvolution
models to coincide and to be
independent of the choice of $M$ and the selection of points $u_1,u_2,
\ldots, u_M$.
Then Section \ref{some_examples} provides examples where these
conditions do or do not take place, and Section \ref{possible_cases}
sorts all possible situations into the uniform, regular and irregular cases.
Section \ref{risk_regular} studies the regular case. In particular,
Section \ref{best_discrete} designs the best
possible sampling strategy. Section \ref{examples} provides some
motivating examples.
Section \ref{equirates} investigates the relation between
the $L^2$-risks in the continuous and the discrete models under an
arbitrary sampling scheme and
formulates conditions when the convergence rates do or do not coincide.
Section \ref{uniform_strategies}
formulates sufficient conditions when the convergence rates in both
models coincide
under a pseudo-uniform sampling scheme in the discrete model. Section
\ref{ex_revisited} provides a variety of examples where the convergence
rates coincide or differ
depending on what sampling scheme is employed.
Section \ref{box-car} explores the interplay between
continuous and discrete functional deconvolution models in the case of
a boxcar-like blurring
function. Section \ref{sim} supplements the theory with a limited
simulation study
in the case of a boxcar-like blurring function
and compares performance of the suggested estimator
to the estimator proposed by De Canditiis and Pensky (\citeyear{CP06}).
Concluding remarks are given in Section \ref{discussion}. Finally,
the \hyperref[append]{Appendix} provides the proofs of the
theoretical results obtained in the previous sections.

In the rest of the paper, the continuous functional deconvolution model
(\ref{convcont}) is referred to as the ``continuous model'' and its
discrete counterpart (\ref{convdis}) is referred to as the ``discrete model.''

\section{Wavelet estimators}
\label{estconstr}

For both the continuous model and the discrete model, we use the
wavelet estimator derived in Pensky and Sapatinas (\citeyear
{PenskyS09}), described as follows.

Let $\varphi^*(\cdot)$ and $\psi^*(\cdot)$ be the Meyer scaling and
mother wavelet functions, respectively, in the real line [see, e.g.,
Meyer (\citeyear{M92})]
and obtain a periodized version of Meyer
wavelet basis as in Johnstone et al. (\citeyear{Johnstoneetal04});
that is, for $j\geq0$ and
$k=0,1,\ldots,2^j-1$,
\begin{eqnarray*}
\varphi_{jk}(x) &=& \sum_{i \in\mathbb Z} 2^{j/2} \varphi^*\bigl(2^j (x
+i) - k\bigr),\\
\psi_{jk}(x) &=& \sum_{i \in\mathbb Z} 2^{j/2} \psi^*\bigl(2^j (x +i) -
k\bigr),\qquad x \in T.
\end{eqnarray*}

Denote $\langle f,g \rangle= \int_{T}f(t)\overline{g(t)} \,dt$,
the inner product in the Hilbert space $L^2(T)$.
Let $e_m(t) = e^{i 2 \pi m t}$, $m \in\mathbb Z$, and, for any (primary
resolution level) $j_0 \geq0$ and any $j \geq j_0$, let
$\varphi_{mj_0k}= \langle e_m, \varphi_{j_0k} \rangle, \psi_{mjk}=
\langle
e_m, \psi_{jk}\rangle, f_m= \langle e_m, f \rangle$
be the Fourier coefficients of $\varphi_{jk}(\cdot)$, $\psi
_{jk}(\cdot)$
and $f(\cdot)$, respectively. For each $u \in U$, denote the
functional Fourier
coefficients by
\[
y_{m}(u)= \langle e_m, y(u,\cdot) \rangle,\qquad
g_{m}(u)= \langle e_m, g(u,\cdot) \rangle.
\]
In what follows we assume that function $g(u,t)$ is such that
$g_{m}(u)$ are
continuous functions of $u$ for every $m$. (This condition is not
restrictive and holds
in all examples considered below.)

If we have the continuous model (\ref{convcont}), then, by using
properties of the Fourier transform, for each $u \in U$, we have
$h_{m}(u)= g_{m}(u)f_m$ and
%
%
\begin{equation} \label{finaleq}
y_{m}(u)= g_{m}(u)f_m+ n^{-1/2} z_{m}(u),
\end{equation}
where $z_{m}(u)$ are generalized one-dimensional (complex-valued) Gaussian pro\-cesses
such that
${\mathbb E}[z_{m_1}(u_1) \overline{z_{m_2}(u_2)}] = \delta_{m_1,m_2} \delta(u_1-u_2)$,
where $\delta_{m_1,m_2}$ is Kronecker's delta. If we have the discrete
model (\ref{convdis}), then, by using
properties of the discrete Fourier transform, for each
$l=1,2,\ldots,M$, (\ref{finaleq}) takes the form
%
%
\begin{equation}\label{finaleqdis}
y_{m}(u_l) = g_{m}(u_l) f_m+ N^{-1/2} z_{ml},
\end{equation}
where $z_{ml}$ are standard (complex-valued) Gaussian random variables, independent for
different $m$ and $l$.


Estimate the Fourier coefficients $f_m$ of $f(\cdot)$ by
%
%
\begin{eqnarray}
\label{fmexprc}
\hat{f}_m & = & \biggl( \int_a^b\overline{g_{m}(u)}y_{m}(u)\,du \biggr) \bigg/ \biggl(
\int_a^b|g_{m}(u)
|^2 \,du \biggr)\nonumber\\[-8pt]\\[-8pt]
\eqntext{\mbox{in the continuous model},}\\
\label{fmexprd}
\hat{f}_m & = & \Biggl( \sum_{l=1}^M\overline{g_{m}(u_l)} y_{m}(u_l) \Biggr) \Bigg/
\Biggl( \sum_{l=1}^M|g_{m}(u_l)|^2 \Biggr)\nonumber\\[-8pt]\\[-8pt]
\eqntext{\mbox{in the discrete model.}}
\end{eqnarray}
Here, we adopt the convention that when $a=b$,
$\hat{f}_m $ takes the form $\hat{f}_m = (\overline
{g_{m}(a)}y_{m}(a)/|g_{m}(a)|^2$
and somewhat abuse notation using $f_m$ for both
functional Fourier coefficients and their discrete counterparts.

Note that, using the periodized Meyer wavelet basis described above
and for any $j_0 \geq0$, any (periodic) $f(\cdot) \in L^2(T)$ can
be expanded as
%
%
\begin{equation}\label{funf}
f(t) = \sum_{k=0}^{2^{j_0}-1} a_{j_0k} \varphi_{j_0k}
(t) + \sum_{j=j_0}^\infty\sum_{k=0}^{2^j -1} b_{jk} \psi_{jk}(t),\qquad
t \in T.
\end{equation}
Furthermore, by Plancherel's formula, the scaling
coefficients, $a_{j_0k}=\langle f, \varphi_{j_0k} \rangle$, and the wavelet
coefficients, $b_{jk}=\langle f,\psi_{jk}\rangle$, of $f(\cdot)$ can
be represented as
%
%
\begin{equation} \label{alkandblk}
a_{j_0k}= \sum_{m \in C_{j_0}} f_m
\overline{\varphi_{mj_0k}},\qquad b_{jk}= \sum_{m \in C_j} f_m
\overline{\psi_{mjk}},
\end{equation}
where $C_{j_0} = \{ m\dvtx
\varphi_{mj_0k}\neq0 \}$ and, for any $j \geq j_0$, $C_j = \{ m\dvtx
\psi_{mjk}\neq0 \}$. We estimate $a_{j_0k}$ and
$b_{jk}$ by substituting $f_m$ in (\ref{alkandblk}) with (\ref
{fmexprc}) or
(\ref{fmexprd}), that is,
%
%
\begin{equation}\label{coefest}
\hat{a}_{j_0k}= \sum_{m \in C_{j_0}} \hat{f}_m
\overline{\varphi_{mj_0k}},\qquad \hat{b}_{jk}= \sum_{m \in C_j}
\hat{f}_m
\overline{\psi_{mjk}}.
\end{equation}

We now construct a (block thresholding) wavelet estimator of
$f(\cdot)$, suggested by Pensky and Sapatinas (\citeyear{PenskyS09}).
For this
purpose, we divide the wavelet coefficients at each resolution level
into blocks of length $\ln n$. Let $A_j$ and $U_{jr}$ be the following
sets of indices:
$
A_j= \{ r \mid r=1,2,\ldots, 2^j/\ln n \}, $ $ U_{jr}= \{
k \mid k = 0,1, \ldots, 2^j-1; (r-1) \ln n \leq k \leq r \ln n -1
\}.
$
Denote
%
%
\begin{equation}\label{bjr}
B_{jr}= \sum_{k \in U_{jr}}b_{jk}^2,\qquad \hat{B}_{jr}= \sum_{k \in
U_{jr}}\hat{b}_{jk}^2.
\end{equation}
Finally, for any $j_0 \geq0$, $f(\cdot)$ is
constructed as
%
%
\begin{eqnarray}\label{fest}
\hat{f}_n(t) &=& \sum_{k=0}^{2^{j_0} -1} \hat{a}_{j_0k}
\varphi_{j_0k} (t)\nonumber\\[-8pt]\\[-8pt]
&&{} + \sum_{j=j_0}^{J-1} \sum_{r \in A_j}\sum_{k \in
U_{jr}}\hat{b}_{jk}{\mathbb I}(|\hat{B}_{jr}|
\geq\lambda_{j}) \psi_{jk}(t),\qquad t \in T,\nonumber
\end{eqnarray}
where ${\mathbb I}(A)$ is the
indicator function of the set $A$, and the resolution levels $j_0$
and $J$ and the thresholds $\lambda_{j}$ will be defined in
Section \ref{upbounds}.

In what follows, the symbol $C$ is used for a generic positive
constant, independent of $n$, while the symbol $K$ is used for a
generic positive
constant, independent of $m$, $n$, $M$ and $u_1,u_2,\ldots,u_M$, which
either of them may take different values at
different places.


\section{Minimax lower and upper bounds for the $L^2$-risk over Besov balls}
\label{risk_lower_upper}

Among the various characterizations of Besov spaces for periodic
functions defined on $L^p(T)$ in terms of wavelet bases, we recall
that for an $r$-regular multiresolution analysis with $0< s < r$ and
for a Besov ball, $B_{p,q}^s (A) =\{f(\cdot) \in L^p(T)\dvtx f \in
B_{p,q}^s, \|f\|_{B_{p,q}^s} \leq A \}$, of radius $A>0$ with $1
\leq p,q \leq\infty$, one has that, with $s' = s+1/2-1/p$,
%
%
\begin{eqnarray}
\label{bpqs}
B_{p,q}^s (A) &=& \Biggl\{ f(\cdot) \in L^p(T)\dvtx \Biggl(
{\sum_{k=0}^{2^{j_0}-1}}|a_{j_{0}k}|^p \Biggr)^{{1}/{p}}\nonumber\\[-8pt]\\[-8pt]
&&\hspace*{4.8pt}{} + \Biggl(
\sum_{j=j_0}^{\infty} 2^{js'q} \Biggl( {\sum_{k=0}^{2^j-1}|b_{jk}}|^p
\Biggr)^{q/p} \Biggr)^{{1/q}} \leq A \Biggr\}\nonumber
\end{eqnarray}
with respective sum(s) replaced by maximum if $p=\infty$ or
$q=\infty$ [see, e.g., Johnstone et al. (\citeyear{Johnstoneetal04}),
Section 2.4].
(Note that, for the Meyer wavelet basis,
considered in Section \ref{estconstr}, $r=\infty$.)

We construct below asymptotical minimax lower bounds for the
$L^2$-risk based on observations from either the continuous model
or the discrete model. For this purpose,
we define the corresponding minimax $L^2$-risks over the set
$\Omega$ as
%
%
\begin{eqnarray}
\label{risk:cont}
R_n^c (\Omega) & = & \inf_{\tilde{f}_n^c} \sup_{f \in\Omega}
{\mathbb E}
\| \tilde{f}_n^c - f \|^2, \\
\label{risk:point_dis}
R_n^d (\Omega, \underline{u}, M) & = & \inf_{\tilde{f}_n^d}
\sup_{f \in\Omega} {\mathbb E}\| \tilde{f}_n^d - f \|^2, \\
\label{risk:dis}
R_n^d (\Omega) & = & \inf_{\underline{u},M} R_n^d (\Omega, \underline{u},
M),
\end{eqnarray}
where the infimum in (\ref{risk:cont}) is taken over all possible
estimators (i.e., measurable functions) $\tilde{f}_n^c(\cdot)$ of
$f(\cdot)$ from the continuous
model, the infimum in (\ref{risk:dis}) is taken over all
possible estimators $\tilde{f}_n^d(\cdot)$ of $f(\cdot)$ from the
discrete model, based on a sample at $M$ points
$\underline{u}= (u_1,u_2, \ldots, u_M)$ and the infimum in (\ref{risk:dis})
is evaluated over all possible estimators $\tilde{f}_n^d(\cdot)$
of $f(\cdot)$ and the choices of $M$ and $\underline{u}$.
Note that, since the asymptotical minimax convergence rates for the $L^2$-risk
in the discrete model depends on $M$ and $\underline{u}$ if these
quantities are fixed, we are interested in the selection of $M$ and
$\underline{u}$, minimizing the asymptotical minimax convergence rates
for the $L^2$-risk.

Denote
$
s^* = s+1/2-1/p', p'=\min(p,2)
$
and, for $\kappa=1,2$, define
%
%
\begin{equation}\label{taum}\quad
\tau_\kappa^c(m) = {\int_{a}^{b}}|g_m(u)|^{2 \kappa} \,du \quad\mbox
{and}\quad
\tau_\kappa^d(m,\underline{u},M) = \frac{1}{M}
\sum_{l=1}^{M}|g_m(u_l)|^{2 \kappa},
\end{equation}
where $\tau_1^c (m) = |g_m(a)|^{2}$ when $a=b$.

Pensky and Sapatinas (\citeyear{PenskyS09}) constructed asymptotical
minimax lower and
upper bounds for the $L^2$-risk for
the continuous model. The corresponding bounds for the discrete model
were obtained
under the very restrictive conditions that the upper and the lower
bounds on $\tau_1^d(m,\underline{u},M)$
do not depend on $n$, $M$ and $\underline{u}$. Below we shall need asymptotic
lower and upper bounds for the $L^2$-risk in the case of
much more general expressions for $\tau_1^c(m)$ and $\tau
_1^d(m,\underline{u},M)$, than in Pensky and Sapatinas (\citeyear{PenskyS09}).

\subsection{Minimax lower bounds: Particular choice of sampling points}
\label{lowbounds_uniform}

Let, with some abuse of notation, $\tau_1 (m) = \tau_1^c(m)$, $R_n^*
(B_{p,q}^s (A)) = R_n^c (B_{p,q}^s (A)) $, in the continuous model, and
$\tau_1 (m) = \tau_1^d(m,\underline{u},M)$, $R_n^* (B_{p,q}^s (A))
=\break
R_n^d (B_{p,q}^s (A), \underline{u}, M) $, in the discrete model.

Assume that for some constants $\nu\in\mathbb R$, $\lambda\in
\mathbb R$, $\alpha
\geq0$ and $\beta
> 0$, independent of $m$ and $n$, and for some sequence $\varepsilon_n>0$,
independent of $m$,
%
%
\begin{equation} \label{cond1}
\tau_1(m) \leq K \varepsilon_n |m|^{-2\nu}({\ln}|m|)^{-\lambda} \exp
(-\alpha
|m|^\beta),\qquad
\nu>0 \mbox{ if } \alpha=0.
\end{equation}
Denote ${n^*}= n \varepsilon_n$ and assume that the sequence
$\varepsilon_n$ is such that
%
%
\begin{equation} \label{ns_prop1}
{n^*}= n \varepsilon_n \rightarrow\infty\qquad\mbox{as } n \rightarrow
\infty.
\end{equation}
Then the following statement is true.
\begin{theorem} \label{th:lower} Let $\{\phi_{j_0,k}(\cdot),\psi
_{j,k}(\cdot)\}$ be the periodic Meyer wavelet basis discussed in Section
\ref{estconstr}.
Let $s >\max(0,1/p-1/2)$, $1 \leq p \leq\infty$, $1 \leq q \leq
\infty
$ and
$A>0$. Let assumptions (\ref{cond1}) and (\ref{ns_prop1}) hold. Then,
as $n
\rightarrow\infty$,
%
%
\begin{equation}\label{low1}
R_n^* (B_{p,q}^s (A)) \geq\cases{
C ({n^*})^{-{2s}/({2s+2\nu+1})} (\ln{n^*})^{{2s \lambda
}/({2s+2\nu+1})},\cr
\hspace*{92.8pt} \mbox{if $\alpha=0, \nu(2-p) < p{s^*}$,}\cr
C \biggl( \dfrac{\ln{n^*}}{{n^*}} \biggr)^{{2{s^*}}/({2s^*+2\nu})}
(\ln{n^*})^{{2{s^* \lambda}}/({2s^*+2\nu})}, \cr
\hspace*{92.8pt}\mbox{if $\alpha=0, \nu(2-p) \geq p{s^*}$,}\cr
C (\ln{n^*})^{-{2{s^*}}/{\beta}}, \qquad\mbox{if $\alpha>0$.}}
\end{equation}
\end{theorem}

\subsection{Minimax upper bounds: Particular choice of sampling points}
\label{upbounds}

Assume that for the constants $\nu\in\mathbb R$, $\lambda\in\mathbb
R$, $\alpha
\geq0$ and
$\beta>0$ and the sequence $\varepsilon_n >0$ in (\ref{cond1})
%
%
\begin{equation}\label{cond2}
\tau_1 (m) \geq K \varepsilon_n |m|^{-2\nu} ({\ln}|m|)^{-\lambda}
\exp(-\alpha|m|^\beta),\qquad \nu>0 \mbox{ if }
\alpha=0.
\end{equation}

Let $\hat{f}_n(\cdot)$ be the wavelet estimator defined by (\ref{fest}).
Let, as before, ${n^*}= n \varepsilon_n$ satisfy condition (\ref
{ns_prop1}), and
assume that in the case of $\alpha=0$ in (\ref{cond2})
the sequence $\varepsilon_n$ is such that
%
%
\begin{equation} \label{ns_prop2}
-h_1 \ln n \leq\ln(1/\varepsilon_n) \leq(1-h_2) \ln n
\end{equation}
for some constants $h_1 >0$ and $h_2 \in(0,1)$. Observe that condition
(\ref{ns_prop2}) implies
(\ref{ns_prop1}) and that $\ln{n^*}\asymp\ln n$ as $n \rightarrow
\infty
$. Here, and in what follows, $u(n) \asymp v(n)$ means that there exist
constants $C_1>0$ and $C_2>0$,
independent of $n$, such that $0<C_1 v(n) \leq u(n) \leq C_2
v(n)<\infty
$ for $n$ large enough.

Choose ${j_0}$ and $J$ such that
%
%
\begin{eqnarray}
\label{jpower}
2^{j_0} &=& \ln({n^*}),\qquad 2^J = ({n^*})^{{1}/({2\nu+1})}\qquad \mbox{if }
\alpha
=0, \\
\label{jexp}
2^{j_0} &=& \frac{3}{8 \pi} \biggl( \frac{\ln({n^*})}{2\alpha}\biggr)
^{{1}/{\beta}},\qquad 2^J=2^{{j_0}}\qquad \mbox{if } \alpha>0.
\end{eqnarray}
[Since ${j_0}>J-1$ when $\alpha>0$, the estimator
(\ref{fest}) only consists of the first (linear) part, and hence
$\lambda_j$ does not need to be selected in this case.] Set, for some
constant $\mu>0$, large enough,
%
%
\begin{equation}\label{lamj}
\lambda_j = \mu^2 ({n^*})^{-1} \ln({n^*}) 2^{2 \nu j} j^\lambda\qquad
\mbox{if } \alpha=0.
\end{equation}
Note that the choices of ${j_0}$, $J$ and $\lambda_j$ are independent
of the
parameters, $s$,
$p$, $q$ and $A$ of the Besov ball $B_{p,q}^s (A)$; hence the estimator
(\ref{fest}) is
adaptive with respect to these parameters.

Set $(x)_{+} = \max(0,x)$, and define
%
%
\begin{equation} \label{rovalue}
\varrho= \cases{
\dfrac{(2\nu+1)(2-p)_+}{p(2s+2\nu+1)}, &\quad if $\nu(2-p) <
p{s^*}$, \cr
\dfrac{(q-p)_+}{q}, &\quad if $\nu(2-p) = p{s^*}$,\cr
0, &\quad if $\nu(2-p) > p{s^*}$.} 
\end{equation}

For any $j\geq j_0$, let $|C_j|$ be the
cardinality of the set $C_j$; note that, for Meyer wavelets, $|C_j|
= 4 \pi2^j $ [see, e.g., Johnstone et al. (\citeyear
{Johnstoneetal04})]. Let also
%
%
\begin{equation} \label{Delt}
\Delta_\kappa(j) =
|C_j|^{-1} \sum_{m \in C_j} \tau_{\kappa} (m) [\tau_1 (m)]^{-2
\kappa},\qquad \kappa=1,2.
\end{equation}
Direct calculations yield that under conditions (\ref{cond2}) and
(\ref{ns_prop2}),
for some constants $c_1 >0$ and $c_2 >0$, independent of $n$, one has
%
%
\begin{equation} \label{delta1}
\Delta_1(j) \leq\cases{
c_1 (\varepsilon_n)^{-1} 2^{2 \nu j} j^\lambda, &\quad
if $\alpha=0$, \cr
c_2 (\varepsilon_n)^{-1} 2^{2 \nu j} j^\lambda\exp\biggl\{ \alpha\biggl( \dfrac
{8\pi}{3}
\biggr) ^\beta
2^{j \beta} \biggr\}, &\quad if $\alpha> 0$.} 
\end{equation}

The proof of the minimax upper bounds for the $L^2$-risk is based on
the following two lemmas.
\begin{lemma} \label{l:coef}
Let assumption (\ref{cond2}) hold, and let the estimators
$\hat{a}_{j_0k}$ and $\hat{b}_{jk}$ of the scaling and
wavelet coefficients $a_{j_0k}$
and $b_{jk}$, respectively, be given by the formula (\ref{coefest}) with
$\hat{f}_m $ defined by (\ref{fmexprc}) in the continuous model
and by
(\ref{fmexprd}) in the discrete model. Then, for
all $j \geq j_0$,
%
%
\begin{equation} \label{ha}\quad
{\mathbb E}|\hat{a}_{j_0k}- a_{j_0k}|^2 \leq C n^{-1} \Delta_1
({j_0}),\qquad
{\mathbb E}|\hat{b}_{jk}- b_{jk}|^{2} \leq C n^{-1} \Delta_1 (j).
\end{equation}
%
If $\alpha=0$ and assumption (\ref{ns_prop2}) holds, then, for any $j
\geq
j_0$, one has
%
%
\begin{equation}\label{hbb}
{\mathbb E}|\hat{b}_{jk}- b_{jk}|^{4} \leq C n (\ln n)^{3 \lambda
} ({n^*})^{-{3}/({2
\nu+1})}.
\end{equation}
\end{lemma}
\begin{lemma} \label{l:deviation}
Let the estimators $\hat{b}_{jk}$ of the wavelet coefficients
$b_{jk}$ be
given by the formula (\ref{coefest}) with $\hat{f}_m $ defined by
(\ref{fmexprc}) in the continuous model and by (\ref{fmexprd}) in the
discrete model. Let assumptions (\ref{cond2}) (if $\alpha=0$) and
(\ref{ns_prop2}) hold. If
%
%
\begin{equation} \label{mu_cond}
\mu\geq2 \sqrt{c_1} \bigl(\sqrt{6}+1\bigr)/ \sqrt{h_2},
\end{equation}
where $c_1$ and $h_2$ are defined in (\ref{delta1}) and (\ref{ns_prop2}),
respectively, then,
for all $j \geq j_0$,
%
%
\begin{equation} \label{probbound}
{\mathbb P}\biggl( {\sum_{k \in U_{jr}}}|\hat{b}_{jk}- b_{jk}|^2 \geq(4
{n^*})^{-1} \mu^2 2^{2\nu j}
j^\lambda \ln({n^*})
\biggr) \leq n^{-3}.
\end{equation}
\end{lemma}

Then the following statement is true.
\begin{theorem} \label{th:upper}
Let $\hat{f}_n(\cdot)$ be the
wavelet estimator defined by (\ref{fest}), with ${j_0}$ and $J$ given by
(\ref{jpower}) (if $\alpha= 0$) or (\ref{jexp}) (if $\alpha>0$) and
$\mu$
satisfying (\ref{mu_cond}).
Let $s > 1/p'$, $1 \leq p \leq\infty$, $1 \leq q \leq\infty$ and $A>0$.
Then, under the assumptions (\ref{cond2}) and (\ref{ns_prop1}) if
$\alpha
>0$, or (\ref{cond2}) and (\ref{ns_prop2})
if $\alpha=0$, as $n \rightarrow\infty$,
%
%
\begin{eqnarray}\label{up}
&&\sup_{f \in B_{p,q}^s (A)} {\mathbb E}\|\hat{f}_n-f\|^2 \nonumber\\[-8pt]\\[-8pt]
&&\qquad\leq
\cases{C ({n^*})^{-{2s}/({2s+2\nu+1})} ( \ln n)^{\varrho+
{2s\lambda
}/({2s+2\nu+1})}, \cr
\hspace*{99.37pt}\mbox{if $\alpha=0, \nu(2-p) < p{s^*}$},
\cr
C \biggl( \dfrac{\ln n}{{n^*}} \biggr)^{{2{s^*}}/({2s^*+2\nu})}
(\ln n )^{\varrho+{2{s^* \lambda}}/({2s^*+2\nu})},\cr
\hspace*{99.37pt}\mbox{if $\alpha=0, \nu(2-p) \geq p{s^*}$},
\cr
C (\ln({n^*}))^{-{2{s^*}}/{\beta}}, \qquad\mbox{if
$\alpha>0$.}}\nonumber
\end{eqnarray}
\end{theorem}
\begin{remark}
Note that in the continuous model, one can write a lower bound
for $\tau_1(m)$ in (\ref{cond1}) and
an upper bound for $\tau_1(m)$ in (\ref{cond2}) with $\varepsilon_n =1$,
so that ${n^*}=n$ in (\ref{low1}) and (\ref{up}). However, in the discrete
model this may not be possible.
Theorems \ref{th:lower} and \ref{th:upper} allow to account for the
dependence of $\tau_1(m)$ on $n$
in the case of the discrete model as well as for an extra logarithmic
factor in the expression of $\tau_1(m)$
which often appears in the case of the continuous model.
\end{remark}
\begin{remark}
Note that Theorems \ref{th:lower} and \ref{th:upper} can be
applied even if the values
of $\nu, \lambda$, $\alpha$ and $\beta$ in assumptions (\ref
{cond1}) and (\ref{cond2}) are different, that may also depend on $M$
and $\underline{u}$.
Then Theorem \ref{th:lower} provides asymptotical minimax lower bounds
for the $L^2$-risk while Theorem \ref{th:upper}
provides the corresponding upper bounds. If,
in the continuous model or in the discrete model with some particular
choice of $M$ and sampling points $\underline{u}$,
the values of $\nu, \lambda$, $\alpha$ and $\beta$ and the
functions $\varepsilon
_n$ in conditions
(\ref{cond1}) and (\ref{cond2}) coincide, then Theorems \ref
{th:lower} and
\ref{th:upper}
imply that the estimator $\hat{f}_n(\cdot)$ defined by (\ref{fest}) is
asymptotically optimal (in the minimax sense), or near-optimal
within a logarithmic factor, over a wide range of Besov balls
$B_{p,q}^s (A)$.
Therefore, in the rest of the paper, when we talk about \textit
{convergence rates}
we refer to the asymptotical minimax lower bounds for the $L^2$-risk
which are attainable, up to at most a logarithmic factor,
according to Theorems \ref{th:lower} and \ref{th:upper}.
\end{remark}

\section{The interplay between continuous and
discrete models: Uniform, regular and irregular cases}
\label{discrcont}

The convergence rates in
the discrete model depend on two aspects: the total
number of observations $n = NM$ and the behavior of $\tau
_1^d(m,\underline{u},M)$. In the continuous model, the values of $\tau
_1^c (m)$ are fixed; they depend on $m$ only, and
hence conditions (\ref{cond1}) and (\ref{cond2}) can be easily verified.
However, this is no longer true in the discrete model; in this case,
the values of $\tau_1^d(m,\underline{u},M)$ may depend on the choice of
$M$ and the
selection of points $\underline{u}$. If we require
the values of $\tau_1^d(m,\underline{u},M)$ to be independent of the
choice of
$M$ and the selection of points $\underline{u}$, then the convergence
rates in the discrete and
the continuous models coincide and are independent of the selection
of points $\underline{u}$. Moreover, in this case, the wavelet estimator
(\ref{fest}) is asymptotically optimal (in the minimax sense) no matter
what the choice of $M$ is. It is quite possible, however, that in
the discrete model, conditions (\ref{cond1}) and (\ref{cond2}) both hold
but with different values of $\nu$, $\lambda$, $\alpha$ and $\beta$ for
different choices of $M$ and $\underline{u}$.
In this case, the asymptotical minimax upper bounds for the $L^2$-risk
in the
discrete model may not coincide with the convergence rates in the
continuous model, at least for some sampling schemes.

\subsection{Necessary and sufficient conditions for convergence rates
equivalency between continuous and discrete models}
\label{suf_uniform}

Assume that there exist points $u_*, u^* \in[a,b]$, independent of
$m$, such that, for any $u \in[a,b]$,
%
%
\begin{equation}\label{uniform_cond1}
|g_m(u)| \leq K |g_m(u^*)| \quad\mbox{and}\quad |g_m(u)| \geq K |g_m(u_*)|.
\end{equation}
In this case,
\[
(b-a) K^2 |g_{m}(u_*)|^2 \leq\tau_1^c (m) \leq(b-a) K^2 |g_{m}(u^*)|^2
\]
and
$
K^2 |g_{m}(u_*)|^2 \leq\tau_1^d (m, \underline{u},M) \leq K^2
|g_{m}(u^*)|^2
$
for any $M$ and $\underline{u}$. Note that, based on the assumption on the
blurring function $g(\cdot,\cdot)$ made in Section \ref{estconstr},
points $u_*$ and $u^*$ satisfying condition (\ref{uniform_cond1})
always exist; however, they are not necessarily independent of $m$.

Here, and in what follows, $u_m \asymp v_m$ means that there exist
constants $C_1>0$ and $C_2>0$,
independent of $m$, such that $0<C_1 v_m \leq u_m \leq C_2 v_m<\infty$
for $|m|$ large enough.

The following statement, which substantially extends Proposition 1 of
Pensky and
Sapatinas (\citeyear{PenskyS09}), presents the necessary and
sufficient conditions
for the convergence rates in the discrete model
to be independent of the choice of $M$ and the selection of points
$\underline{u}$ and hence to coincide with the convergence rates in the
continuous model.
\begin{theorem} \label{th:condisrates}
Let there exist constants $\nu_1 \in\mathbb R$, $\nu_2 \in\mathbb
R$, $\alpha_1
\geq0$, $\alpha_2 \geq0$, $\beta_1 >0$ and $\beta_2 >0$,
independent of
$m$ and $n$, such that
%
%
\begin{eqnarray}
\label{disconHigh}
|g_{m}(u^*)|^2 &\asymp& |m|^{-2\nu_1} \exp(-\alpha_1
|m|^{\beta_1}),\qquad
\nu_1>0 \mbox{ if } \alpha_1 =0,
\\
\label{disconLow}
|g_{m}(u_*)|^2 &\asymp& |m|^{-2\nu_2} \exp(-\alpha_2
|m|^{\beta_2}),\qquad
\nu_2>0 \mbox{ if } \alpha_2 =0.
\end{eqnarray}
Then, the convergence rates obtained in Theorems \ref{th:lower} and
\ref
{th:upper}
in the discrete model are independent of the choice
of $M$ and the selection of points $\underline{u}$, and
hence coincide with the convergence rates
obtained in Theorems \ref{th:lower} and \ref{th:upper}
in the continuous model, if and only if
%
%
\begin{equation}
\label{th3cond}
\alpha_1 \alpha_2 > 0 \quad\mbox{and}\quad
\beta_1=\beta_2 \quad\mbox{or}\quad \alpha_1=\alpha_2 = 0 \quad\mbox{and}\quad \nu_1 =
\nu_2.
\end{equation}
\end{theorem}
\begin{remark}
Theorem \ref{th:condisrates} provides necessary and sufficient
conditions for the convergence rates in the continuous and the discrete
models to coincide, and to be
independent of the choice of $M$ and the selection of points
$\underline
{u}$. These conditions also guarantee asymptotical
optimality (in the minimax sense) of the wavelet estimator
(\ref{fest}) and can be viewed as some kind of uniformity conditions.
Under assumptions (\ref{disconHigh})--(\ref{th3cond}), asymptotically (up
to a constant factor) it makes absolutely no difference whether one
samples the discrete model $n$ times at one point, say, $u_1$
or, say, $\sqrt{n}$ times at $M = \sqrt{n}$ points $u_l$,
$l=1,2,\ldots
,M$. In other
words, each sample value $y(u_l, t_i)$, $l=1,2,\ldots, M$,
$i=1,2,\ldots, N$, asymptotically (up to a constant factor) gives
the same amount of information, and, therefore, the convergence rates
are not sensitive
to the choice of $M$ and the selection of points $\underline{u}$.
On the other hand, if the conditions of Theorem \ref{th:condisrates}
are violated, then the convergence rates in the discrete model depend
on the choice of $M$
and $\underline{u}$, and some recommendations on their
selection should be given. Furthermore, optimality (in the minimax
sense) issues become
much more complex when $\tau_1^d (m; \underline{u}, M)$ is not
uniformly bounded from
above or below.
\end{remark}

\subsection{Some illustrative examples}
\label{some_examples}

Theorem \ref{th:condisrates} provides necessary and sufficient
conditions for the continuous and the discrete models to be equivalent,
from the viewpoint of convergence rates,
no matter what the choice of $M$ and the selection of points
$\underline{u}$
are. The difficulty, however, is that many models do not satisfy those
conditions.
Below, we consider some illustrative examples that have recently been
studied in Pensky and Sapatinas (\citeyear{PenskyS09}).
\begin{example}[(Estimation of the initial condition in the heat
conductivity equation)]
\label{ex:1}
Let $h(t,x)$ be a solution of the heat conductivity equation,
\[
\frac{\partial h (t,x)}{\partial t} = \frac{\partial^2 h
(t,x)}{\partial x^2},\qquad x \in[0,1], t \in[a,b], a>0,
b < \infty,
\]
with initial condition $h(0,x) = f(x)$ and periodic boundary
conditions
$h(t,0) = h(t,1)$ and $\partial h(t,x)/\partial x |_{x=0}
= \partial h(t,x)/\partial x |_{x=1}$.
%

We assume that a noisy solution $y(t,x) = h(t,x) + n^{-1/2}z(t,x)$
is observed, where $z(t,x)$ is a generalized two-dimensional
Gaussian field with covariance function ${\mathbb E}[z (t_1, x_1) z (t_2,
x_2)] = \delta(t_1-t_2) \delta(x_1-x_2)$, and the goal is to recover
the initial condition $f(\cdot)$ on the basis of observations
$y(t,x)$. This problem was initially considered by Lattes and Lions
(\citeyear{LL67}) and further studied by
Golubev and Khasminskii (\citeyear{GK99}).

Then the functional
Fourier coefficients $g_m(\cdot)$ are of the form
\[
g_{m}(u)= \exp(- 4 \pi^2 m^2 u),
\]
so that $u_*=b$, $u^* =a$, $|g_m (u_*)| = \exp(- 4 \pi^2 b m^2)$
and $|g_m (u^*)| = \exp(- 4 \pi^2a  \times\break m^2)$ [see Example 1 in Pensky and
Sapatinas (\citeyear{PenskyS09})]. Hence Theorem~\ref{th:condisrates}
holds with $\nu_1 = \nu_2=0$, $\alpha_1 = 4 \pi^2 b$, $\alpha_2 =
4 \pi
^2 a$
and $\beta_1 = \beta_2 = 2$. Therefore, the convergence rates in the
continuous and the discrete models coincide
and are independent of the choice of $M$ and the selection of points
$\underline{u}$.
\end{example}
\begin{example}[(Estimation of the boundary condition for the
Dirichlet problem of the Laplacian on the unit circle)]
\label{ex:2}
Let $h(x,w)$ be a solution of the Dirichlet problem of the Laplacian
on a region $D$ on the plane
%
%
\begin{equation} \label{laplace}
\frac{\partial^2 h (x,w)}{\partial
x^2} + \frac{\partial^2 h (x,w)}{\partial w^2} =0,\qquad (x, w) \in D
\subseteq\mathbb R^2,
\end{equation}
with a boundary $\partial D$
and boundary condition $h(x,w) |_{\partial D} = F(x,w)$.
Consider the situation when $D$ is the unit
circle. Then it is advantageous to rewrite the function
$h(\cdot,\cdot)$ in polar coordinates as $h(x,w) = h(u,t)$ where $u
\in[0,1]$ is the polar radius and $t \in[0, 2\pi]$ is the polar
angle. Then the boundary condition
can be presented as $h(1, t) = f(t)$, and $h(u,\cdot)$ and $f(\cdot)$ are
periodic functions of $t$ with period $2\pi$.

Suppose that only a noisy version $y(u,t) = h(u,t) + n^{-1/2}z(u,t)$
is observed, where $z(u,t)$ is as in Example \ref{ex:1}, and that
observations are available only on the interior of the unit circle
with $u \in[0, r_0]$, $r_0 <1$, that is, $a=0, b=r_0<1$. The goal is
to recover the boundary condition $f(\cdot)$ on the basis of
observations $y(u,t)$. This problem was initially investigated in
Golubev and
Khasminskii (\citeyear{GK99}) and Golubev (\citeyear{G04}).

Then the functional Fourier coefficients $g_m(\cdot)$ are of the form
%
%
\begin{equation} \label{gmu_ex2}
|g_{m}(u)| = K u^{|m|} = K \exp\bigl( - |m| \ln(1/u)\bigr),\qquad u \in[0, r_0],
\end{equation}
so that $u_* = 0$, $u^* = r_0$,
$|g_m (u_*)| = 0$ and \mbox{$|g_m (u^*)| = K \exp(- |m| \ln(1/r_0))$} [see
Pensky and Sapatinas (\citeyear{PenskyS09}), Example 2].
Hence, the conditions of Theorem \ref{th:condisrates} do not hold,
and we cannot be certain that the convergence rates in the continuous
and the discrete models coincide
for any sampling scheme. Actually, it is easy to see that if sampling
is carried out
entirely at the single point $u_*=0$, then $\tau_1^d (m, u_*, 1)=0$,
and we cannot recover the boundary condition $f(\cdot)$.
\end{example}
\begin{example}[(Estimation of the speed of a wave on a finite
interval)]
\label{ex:3}
Let $h(t,x)$ be a solution of the wave equation
\[
\frac{\partial^2 h (t,x)}{\partial t^2} = \frac{\partial^2 h
(t,x)}{\partial x^2}
\]
with initial-boundary conditions
$h(0,x)=0$, $\partial h(t,x)/\partial t |_{t=0} =
f(x)$ and $h(t$, $0) = h(t,1)=0$.

Here $f(\cdot)$ is a function defined on the unit interval
$[0,1]$, and the goal is to recover the speed of a wave $f(\cdot)$ on
the basis of
observing a noisy solution $y(t,x) = h(t,x) + n^{-1/2}z(t,x)$ where
$z(t,x)$ is as in Example \ref{ex:1} with $t \in[a,b]$, $a
>0$, $b<1$.

Then the functional Fourier coefficients $g_m(\cdot)$
are of the form
%
%
\begin{eqnarray}
\label{waverefl}
g_0(u)=1 \quad\mbox{and}\quad g_{m}(u)= (2\pi m)^{-1} \sin(2 \pi m
u),\nonumber\\[-8pt]\\[-8pt]
\eqntext{m \in\mathbb Z\setminus\{0\}, u \in[a,b],}
\end{eqnarray}
[see Pensky and Sapatinas (\citeyear{PenskyS09}), Example 4].
It is easy to see that in order to satisfy the condition (\ref{uniform_cond1})
the points $u_*$ and $u^*$ should depend on $m$, and hence the
convergence rates depend on the selection
of $M$ and $\underline{u}$. Hence the convergence rates in the
continuous and
the discrete models may coincide for one selection of $M$ and
$\underline{u}$
and be different for another. Actually, it is easy to see that if $M=1$
and $u$ is an integer, then
$\tau_1^d (m, u, 1)=0$, and we cannot recover the speed of a wave
$f(\cdot)$.
\end{example}

\subsection{Possible cases}
\label{possible_cases}

Theorem \ref{th:condisrates} in Section \ref{suf_uniform} provides
necessary and sufficient conditions for the convergence rates in the
discrete model to be independent of the choice
of $M$ and the selection of points $\underline{u}$ and
hence to coincide with the convergence rates
in the continuous model. We can divide these
conditions into the following two groups.
\renewcommand{\theCondition}{I}
\begin{Condition}\label{ConditionI}
There exist constants $\nu_1 \in\mathbb R$,
$\alpha_1 \geq
0$ and $\beta_1 >0$ and a point $u^* \in[a,b]$, independent of $m$ and
$n$, such that
%
%
\begin{equation}\label{disconFanis-I}
|g_m (u)|^2 \leq K |g_{m}(u^*)|^2 \asymp|m|^{-2\nu_1} \exp(-\alpha_1
|m|^{\beta_1}),\qquad \nu_1>0 \mbox{ if } \alpha_1 =0.\hspace*{-28pt}
\end{equation}
\end{Condition}
\renewcommand{\theCondition}{I*}
\begin{Condition}\label{ConditionIzv}
There exist constants $\nu_2 \in\mathbb R$,
$\alpha_2 \geq0$
and $\beta_2>0$, and a point $u_* \in[a,b]$, independent of $m$ and
$n$, such that
%
%
\begin{equation}\label{discon-I}
|g_m (u)|^2 \geq K |g_{m}(u_*)|^2 \asymp|m|^{-2\nu_2} \exp(-\alpha_2
|m|^{\beta_2}),\qquad \nu_2>0 \mbox{ if } \alpha_2 =0.\hspace*{-28pt}
\end{equation}
\end{Condition}
\renewcommand{\theCondition}{II}
\begin{Condition}\label{ConditionII}
Either $\alpha_1 \alpha_2 > 0$ and
$\beta_1=\beta_2$ or $\alpha_1= \alpha_2 = 0$ and $\nu_1 = \nu_2$.
\end{Condition}

Consider now the following three cases.

\begin{enumerate}

\item\textit{The uniform case}: Conditions \ref{ConditionI}, \ref{ConditionIzv} and \ref{ConditionII} hold.

\item\textit{The regular case}: Condition \ref{ConditionI} holds but Condition \ref{ConditionII}
does not
hold. Condition~\ref{ConditionIzv} holds or, possibly, $|g_m(u_*)| =0$.

\item\textit{The irregular case}: Condition \ref{ConditionI} does not hold.
\end{enumerate}

It is easy to see that Examples \ref{ex:1}, \ref{ex:2} and \ref
{ex:3} of Section \ref{some_examples} correspond to the uniform case,
the regular case and the irregular case, respectively.

Theorem \ref{th:condisrates} shows that in the uniform case,
the convergence rates obtained in Theorems \ref{th:lower} and \ref
{th:upper} in
the discrete model are independent of the choice of
$M$ and the selection of points $\underline{u}$ and hence
coincide with the convergence rates obtained in Theorems \ref{th:lower}
and \ref{th:upper} in
the continuous model. In the uniform case one can replace the discrete model
by the continuous model, no matter what $M$ and $\underline{u}$ are.

In the regular case, one cannot guarantee that the convergence rates
between continuous and discrete models coincide.
However, as we shall show below, one can still locate a point $u^*$
which delivers the best possible convergence rates.
If sampling is done entirely at this point, then the discrete model can
sometimes deliver better convergence rates than the continuous model.
Nevertheless, if another sampling strategy is chosen, then the
convergence rates in the discrete model may be
worse than in the continuous model. Note that we do not require
Condition \ref{ConditionIzv} to hold.
This is due to the fact that Condition \ref{ConditionIzv} refers to the ``worst case
scenario'' when we sample at the points which
leads to the highest possible variance and, consequently, to the lowest
convergence rates.
One can also view $|g_m(u_*)| =0$ as an extreme case of Condition \ref{ConditionIzv}
when $\nu_2 = \infty$ or $\alpha_2 = \alpha_1$
and $\beta_2 = \infty$. It is easy to see that if, in the discrete
model, all sampling is carried out at $u_*$, then the
convergence rates will be worse than in the case of sampling entirely
at $u^*$ or than in the continuous model.
Hence, in the regular case, sampling strategy \textit{does} matter.

In the irregular case, it is impossible to pinpoint the best sampling
strategy which suits any
problem; this is due to the fact that Condition \ref{ConditionI} can be violated in a
variety of ways. For this reason, we study a particular example
of the irregular case, namely, functional deconvolution with a boxcar-like
blurring function; this important model occurs in the
problem of estimation of the speed of a wave on a finite interval (see
Example \ref{ex:3} in Section \ref{some_examples})
and, a discretized version of it, in many areas of signal and image
processing which include, for instance, LIDAR
(Light Detection and Ranging), remote sensing and reconstruction of
blurred images (see Section \ref{box-car}).

\section{The regular case} \label{risk_regular}

\subsection{The best discrete rates}
\label{best_discrete}

It is easy to see that, in the regular case, $\tau_1^d (m, u^*,1) \geq
K \tau_1^c (m)$.
Hence it follows from Theorems \ref{th:lower} and \ref{th:upper} that,
if the discrete model is sampled entirely at $u^*$ (i.e., $M=1$ and
$u_1 = u^*$),
then the asymptotical minimax lower and upper bounds
for the $L^2$-risk in the discrete model can be only lower than
the respective lower and upper bounds in the continuous model.

Denote by $\hat{f}_n^c(\cdot)$ the wavelet estimator of $f(\cdot)$
defined by (\ref{fest}) based on observations from the continuous model,
and let $\hat{f}_n^d(\cdot) =\hat{f}_n^d(\underline{u},M,\cdot)$ be
the corresponding wavelet estimator of $f(\cdot)$ based on
observations from the discrete model evaluated at
the point $\underline{u}$.
Denote $\hat{f}_n^{d*}(\cdot)= \hat{f}_n^d(u^*,1,\cdot)$.

Then the following statement is true.
\begin{theorem} \label{th:regular_rates} Let $\{\phi_{j_0,k}(\cdot
),\psi
_{j,k}(\cdot)\}$ be the periodic Meyer wavelet basis discussed in Section
\ref{estconstr} and assume that $s >\max(0,1/p-1/2)$ (for the lower
bounds) or
$s>1/p'$ (for the upper bounds), $1 \leq p \leq\infty$, $1 \leq q
\leq\infty$ and $A>0$. Then
%
%
\begin{equation} \label{th41}
R_n^c (B_{p,q}^s (A)) \geq C R_n^d (B_{p,q}^s (A), u^*, 1) \asymp R_n^d
(B_{p,q}^s (A)).
\end{equation}
Also, for any choice of $M$ and $\underline{u}$, we have
%
%
\begin{eqnarray}
\label{th42}
\sup_{f \in B_{p,q}^s (A)} {\mathbb E}\|\hat{f}_n^{d*} -f\|^2 &\leq& C
\sup_{f \in
B_{p,q}^s (A)} {\mathbb E}\|\hat{f}_n^{c} -f\|^2,
\\
\label{th43}
\sup_{f \in B_{p,q}^s (A)} {\mathbb E}\|\hat{f}_n^{d*} -f\|^2 &\leq& C
\sup_{f \in
B_{p,q}^s (A)} {\mathbb E}\|\hat{f}_n^{d} -f\|^2.
\end{eqnarray}
\end{theorem}

Theorem \ref{th:regular_rates} confirms that sampling entirely at the
single point $u^*$ leads to the highest possible convergence rates
in the discrete model. However, it does not provide an answer to the
question whether the inequalities in (\ref{th41}) and (\ref{th42})
are strict
or the convergence rates are the same in the continuous and the
discrete models with
sampling entirely at the single point $u^*$. To get a better insight
into the matter, let us consider a few more examples.

\subsection{More examples}
\label{examples}

\textit{Example} \ref{ex:2} (\textit{continued}).
In the case of estimation of the boundary condition for the
Dirichlet problem of the Laplacian on the unit circle, the functional
Fourier coefficients $g_m (\cdot)$ are
of the form (\ref{gmu_ex2}) with $r_0 <1$. Hence,
$u^* = r_0$ and
$\tau_1^d (m, u^*,1) \asymp|g_m (u^*)|^2 \asymp\exp( - |m| \ln(1/r_0)
)$.
On the other hand,
$
\tau_1^c (m) \asymp\int_0^{r_0} u^{2|m|} \,du = r_0^{2|m|+1}/(2|m|+1)
\asymp|m|^{-1} \exp( -|m| \ln(1/r_0)).
$
Hence, by Theorems \ref{th:lower} and \ref{th:upper}, the convergence
rates in the continuous model coincide with the convergence rates
in the discrete model if sampling is carried out entirely at the single
point $u^*$.
\begin{example}
\label{ex:4} Let the functional Fourier coefficients
$g_m(\cdot)$ satisfy
\[
|g_m(u)|^2 \asymp|m|^{-2u},\qquad 0 < a \leq u \leq b < \infty.
\]
Then, in the continuous model,
\[
\tau_1^c(m) = \int_{a}^{b} |g_m(u)|^2 \,du \asymp\int_{a}^{b}
\exp({-2u \ln}|m|)\,du \asymp|m|^{-2a} ({\ln}|m|)^{-1},
\]
implying that conditions (\ref{cond1}) and (\ref{cond2}) hold with
$\nu=a$, $\alpha=0$ and $\lambda=1$.
In the case of the discrete model, $u^* = a$ and
$
\tau_1^d (m, u^*,1) \asymp|g_m(u^*)|^2 \asymp|m|^{-2a}
$
and conditions (\ref{cond1}) and (\ref{cond2})
hold with $\nu=u^*$, $\alpha=0$ and $\lambda=0$. Hence, by Theorems
\ref
{th:lower} and \ref{th:upper},
the convergence rates in the continuous model are worse than the
convergence rates in the discrete (they differ by a logarithmic factor)
model when sampling
is carried out entirely at the single point $u^*$.
\end{example}
\begin{example}
\label{ex:5}
Let the functional Fourier coefficients
$g_m(\cdot)$ satisfy
%
%
\begin{equation}
\label{ex5gmu}
|g_m(u)|^2 \asymp\exp( - \alpha|m|^u ),\qquad 0 <a \leq u \leq b <
\infty,
\end{equation}
for some constant $\alpha>0$, independent of $m$.
Then $u^* =a$ and $\tau_1^d (m, u^*,1) = |g_m (u^*)|^2 \asymp\exp( -
\alpha|m|^a
)$.
On the other hand,
$\tau_1^c (m) \asymp\int_a^b \exp( - \alpha|m|^u ) \,du \asymp
({\ln}|m|)^{-1} \int_{|m|^a}^{|m|^b} z^{-1} \exp(-\alpha z) \,dz$,
so that
\[
\tau_1^c (m) \geq K |m|^{-b} ({\ln}|m|)^{-1} \int
_{|m|^a}^{|m|^b} \exp(-\alpha z) \,dz
\asymp|m|^{-b} ({\ln}|m|)^{-1} \exp(-\alpha|m|^a)
\]
and
\[
\tau_1^c (m) \leq K |m|^{-a} ({\ln}|m|)^{-1} \int_{|m|^a}^{\infty}
\exp(-\alpha z) \,dz
\asymp|m|^{-a} ({\ln}|m|)^{-1} \exp(-\alpha|m|^a).
\]

Hence, by Theorems \ref{th:lower} and \ref{th:upper}, the convergence
rates in the continuous and the discrete models
coincide if sampling is carried out entirely at the single point $u^*$.
\end{example}
\begin{example}
\label{ex:6}
Let the functional Fourier coefficients
$g_m(\cdot)$ satisfy
%
%
\begin{equation}
\label{ex6gmu}
|g_m(u)|^2 \asymp|m|^{- 2 \nu} \exp(- u|m|^\beta),\qquad 0 \leq u \leq
b <
\infty,
\end{equation}
for some constants $\nu>0$ and $\beta>0$, independent of $m$.
Then, $u^* =0$ and
%
%
\begin{equation}
\label{ex6bestpoint}
\tau_1^d (m, u^*,1) \asymp|g_m (u^*)|^2 \asymp|m|^{- 2 \nu}.
\end{equation}
On the other hand, it is easy to check that
%
%
\begin{equation}
\label{ex6ratescont}
\tau_1^c (m) \asymp|m|^{- 2 \nu} \int_0^b \exp(- u|m|^\beta) \,du
\asymp
|m|^{-(2\nu+ \beta)}.
\end{equation}
Hence, by Theorems \ref{th:lower} and \ref{th:upper}, the convergence rates
in the continuous model are worse than in the discrete model
when sampling is carried out entirely at the single point $u^* =0$.
\end{example}

\subsection{Conditions for convergence rates equivalency and
nonequivalency between continuous and discrete models}
\label{equirates}

We shall say that the convergence rates in the continuous
and the discrete models ``almost coincide'' if the convergence rates
coincide up to, at most, a logarithmic factor when the convergence
rates are polynomial [$\alpha(u) \equiv0$] or
up to, at most, a constant when the convergence rates are logarithmic
[$\alpha(u) > 0$].
We choose this distinction between the cases of polynomial and
logarithmic convergence rates
since in the polynomial case the upper bounds for the risks of the
adaptive estimator
may differ from the corresponding lower bounds for the risk by a
logarithmic factor.

Hence a question naturally arises: under which conditions on the choice
of $M$ and the selection of sampling points
$\underline{u}$ do the convergence rates in the discrete and the
continuous models almost coincide, and under which conditions this does
not happen? In order to answer this question, first we have to derive
upper and lower bounds for the $L^2$-risk in the continuous model.

In what follows we assume that
the functional Fourier coefficients $g_m(\cdot)$ satisfy the assumption
%
%
\begin{equation} \label{equi:reg}
|g_m(u)|^2 \asymp|m|^{-2\nu(u)} \exp\bigl(-\alpha(u) |m|^{\beta(u)}
\bigr),\qquad u \in U,
\end{equation}
for some continuous functions $\nu(\cdot)$, $\alpha(\cdot)$ and
$\beta
(\cdot)$ defined on $u \in U$,
such that either $\alpha(u) =0$ and $\nu(u) >0$
or $\alpha(u)>0$ and $ \beta(u) >0$, for all $u \in U$. Denote
%
%
\begin{equation} \label{vartheta}
\vartheta= \cases{
\dfrac{2s}{k(2s+2\nu(u^*)+1)}, &\quad if
$\nu(u^*)(2-p) < ps^*$,\vspace*{2pt}\cr
\dfrac{2s^*}{k(2s^*+2\nu(u^*))}, &\quad if
$\nu(u^*)(2-p) \geq ps^*$.}
\end{equation}
Then the following statement is valid.
%
%
\begin{lemma}
\label{th:answer} Let $\{\phi_{j_0,k}(\cdot),\psi_{j,k}(\cdot)\}$ be
the periodic Meyer wavelet basis discussed in Section
\ref{estconstr}, and assume that $s >\max(0,1/p-1/2)$ (for the lower
bounds) or
$s>1/p'$ (for the upper bounds), $1 \leq p \leq\infty$, $1 \leq q
\leq\infty$ and $A>0$.
Let also the functional Fourier coefficients
$g_m(\cdot)$ satisfy assumption (\ref{equi:reg}).
Denote
\[
u^* = \cases{
\displaystyle\mathop{\arg\min}_{u \in U} \nu(u), &\quad if $\alpha(u) \equiv0$, \cr
\displaystyle\mathop{\arg\min}_{u \in U} \beta(u), &\quad if $\alpha(u) > 0, \beta(u) \neq
\mathrm{const}$.}
\]
Assume further that, in the neighborhood of point $u=u^*$, the
function $\beta(\cdot)$ is continuously differentiable [if $\alpha(u)
> 0$, $u \in U$] or the function $\nu(\cdot)$ is $k$-times continuously
differentiable
[if $\alpha(u) = 0$, $u \in U$], where $k \geq1$ is such that
%
%
\begin{equation} \label{k_val}
\nu^{(s)}(u^*)=0,\qquad s=1,\ldots,k-1,\qquad \nu^{(k)}(u^*) \neq0
\end{equation}
with $\nu^{(s)}(\cdot)$ denoting the $s$th derivative of the function
$\nu(\cdot)$. Then
the asymptotical minimax lower and upper bounds for the $L^2$-risk
in the continuous model are as follows:
%
%
\begin{eqnarray}\label{low:reg_equi}\qquad
&&R_n^c (B_{p,q}^s (A)) \nonumber\\[-8pt]\\[-8pt]
&&\qquad\geq
\cases{Cn^{-{2s}/({2s+2\nu(u^*)+1})} (\ln n)^{\vartheta}, \cr
\hspace*{109pt}\mbox{if
$\alpha(u) = 0, \nu(u^*)(2-p) < p{s^*}$},
\cr
C \biggl( \dfrac{\ln n}{n} \biggr)^{{2{s^*}}/({2s^*+2\nu(u^*)})} (\ln
n)^{\vartheta}, \cr
\hspace*{109pt}\mbox{if $\alpha(u) = 0, \nu(u^*)(2-p) \geq
p{s^*}$},
\cr
C (\ln n)^{-{2{s^*}}/({\beta(u^*)})}, \qquad\mbox{if $\alpha(u)
>0$},}\nonumber
\end{eqnarray}
and
%
%
\begin{eqnarray}
\label{up:reg_equi}
&&\sup_{f \in B_{p,q}^s (A)} {\mathbb E}\|\hat{f}_n^c
-f\|^2 \nonumber\\[-8pt]\\[-8pt]
&&\qquad\leq \cases{
Cn^{-{2s}/({2s+2\nu(u^*)+1})} (\ln n)^{\rho+ \vartheta}, \cr
\hspace*{109pt}\mbox{if $\alpha(u) = 0, \nu(u^*)(2-p) < p{s^*}$},
\cr
C \biggl( \dfrac{\ln n}{n} \biggr)^{{2{s^*}}/({2s^*+2\nu(u^*)})} (\ln
n)^{\rho+ \vartheta}, \cr
\hspace*{109pt}\mbox{if $\alpha(u) = 0, \nu(u^*)(2-p)
\geq p{s^*}$},
\cr
C (\ln n)^{-{2{s^*}}/({\beta(u^*)})}, \qquad\mbox{if $\alpha(u)
>0$}.}\nonumber
\end{eqnarray}
Here $\rho$ is given by (\ref{rovalue}) with $\nu= \nu(u^*)$, and
$\vartheta$ is given by (\ref{vartheta}).
If $\nu(\cdot)$ is a constant function, then $k=\infty$ in (\ref{k_val})
and $\vartheta=0$.
\end{lemma}
%
%
\begin{remark} \label{rem:nonconst}
In Lemma \ref{th:answer}, we do not consider the case when
$\beta(u)$ is constant since this situation belongs to the uniform case
and the convergence rates in the continuous and the discrete models
coincide for any sampling scheme
due to Theorem \ref{th:condisrates}. Note also that the value of $u^*$
in Lemma \ref{th:answer} is always independent of $m$ and easy
to find.
\end{remark}

The utility of Lemma \ref{th:answer} is that it allows one to formulate
conditions
such that the convergence rates in the continuous model almost coincide
with the convergence rates in the discrete model for any particular
choice of a sampling scheme.
\begin{theorem}
\label{th:sufficient}
Let assumptions (\ref{equi:reg}) and (\ref{k_val}) hold.

\begin{longlist}
\item If $\alpha(u) \equiv0$, then the convergence rates in the
continuous and the discrete models coincide up to
at most a logarithmic factor if $M=M_n$ and $\underline{u}$ are such that
%
%
\begin{equation} \label{suf_pol}
\tau_1^d(m,\underline{u},M_n) \geq K \varepsilon_n |m|^{-2 \nu
(u^*)} ({\ln}|m|)^{-\lambda_1}
\end{equation}
for some constant $\lambda_1 \in\mathbb R$, independent of $m$ and
$n$, and for
some sequence $\varepsilon_n>0$, independent of $m$, satisfying
%
%
\begin{equation} \label{eps_log}
\lim_{n \rightarrow\infty} \varepsilon_n (\ln n)^{\lambda_2} >0
\end{equation}
for some constant $\lambda_2 \geq0$. If, moreover, $\varepsilon_n$,
$M=M_n$ and
$\underline{u}$ are such that opposite inequalities hold, that is,
%
%
\begin{eqnarray} \label{suf_opp}
\tau_1^d(m,\underline{u},M_n) &\leq& C \varepsilon_n |m|^{-2 \nu (u^*)}
({\ln}|m|)^{-\lambda_1} \quad\mbox{and}\nonumber\\[-8pt]\\[-8pt]
\lim_{n \rightarrow\infty} \varepsilon_n (\ln n)^{\lambda_2} &<&
\infty,\nonumber
\end{eqnarray}
for the same constants $\lambda_1$ and $\lambda_2$ as in formulae
(\ref{suf_pol}) and (\ref{eps_log}),
and if $k$ in (\ref{k_val}) is such that $k(\lambda_1 + \lambda_2) =1$,
then the convergence rates in the continuous and discrete models
coincide up to constant.

\item If $\alpha(u) > 0$, then the convergence rates in the continuous
and discrete models coincide up to
constant if $M=M_n$ and $\underline{u}$ are such that
%
%
\begin{equation} \label{suf_exp}
\tau_1^d(m,\underline{u},M_n) \geq K \varepsilon_n |m|^{-2 \nu}
\exp\bigl( -\alpha
|m|^{\beta(u^*)} \bigr)
( {\ln}|m| )^{-\lambda_1}
\end{equation}
for some constants $\nu\in\mathbb R$, $\lambda_1 \in\mathbb R$ and
$\alpha>0$,
independent of $m$ and $n$, and for some sequence $\varepsilon_n>0$,
independent of $m$, satisfying condition
(\ref{ns_prop2}).
\end{longlist}
\end{theorem}

Theorem \ref{th:sufficient} provides sufficient conditions for
a sampling scheme in the discrete model to lead to the convergence
rates which are optimal or near-optimal.
It follows from conditions (\ref{equi:reg}) and (\ref{k_val}) and Theorems
\ref{th:lower} and \ref{th:upper}
that, if the discrete model is sampled entirely at $u^*$, then
the convergence rates in the continuous and the discrete models almost coincide.
Namely, as $n \rightarrow\infty$,
%
%
\begin{eqnarray}\label{low:reg_equi1}
&&R_n^d (B_{p,q}^s (A)) \nonumber\\[-8pt]\\[-8pt]
&&\qquad\geq\cases{
Cn^{-{2s}/({2s+2\nu(u^*)+1})}, &\quad if $\alpha(u)= 0,
\nu(u^*)(2-p) < p{s^*}$,
\cr
C \biggl( \dfrac{\ln n}{n} \biggr)^{{2{s^*}}/({2s^*+2\nu(u^*)})}, &\quad
if $\alpha(u) = 0, \nu(u^*)(2-p) \geq p{s^*}$,
\cr
C (\ln n)^{-{2{s^*}}/({\beta(u^*)})}, &\quad if $\alpha(u)
>0$,}\nonumber\hspace*{-35pt}
\end{eqnarray}
%
and
%
%
\begin{eqnarray}\label{up:reg_equi1}\qquad
&&\sup_{f \in B_{p,q}^s (A)} {\mathbb E}\|\hat{f}_n^{d*} -f\|^2\nonumber\\[-8pt]\\[-8pt]
&&\qquad\leq
\cases{
Cn^{-{2s}/({2s+2\nu(u^*)+1})} ( \ln n )^{\varrho}, \cr
\hspace*{108.7pt}\mbox{if $\alpha(u) = 0, \nu(u^*)(2-p) < p{s^*}$},
\cr
C \biggl(\dfrac{\ln n}{n} \biggr)^{{2{s^*}}/({2s^*+2\nu(u^*)})}
( \ln n )^{\varrho}, \cr
\hspace*{108.7pt}\mbox{if $\alpha(u) = 0, \nu(u^*)(2-p)
\geq p{s^*}$},
\cr
C (\ln n)^{-{2{s^*}}/({\beta(u^*)})}, \qquad\mbox{if $\alpha(u)
>0$}.}
\nonumber
\end{eqnarray}

From the above, it also follows that
\[
\frac{R_n^c (B_{p,q}^s (A))}{R_n^{d} (B_{p,q}^s
(A))} \asymp\cases{
1, &\quad if $\alpha(u) >0$ and $\beta(u)>0, u \in U$,\cr
(\ln n)^{\vartheta}, &\quad if $\alpha(u)=0, u \in U$,}
\]
and hence the convergence rates in the discrete model cannot be better
than the convergence rates in the
continuous model if $\alpha(u) >0$ and cannot be better by more than a
logarithmic factor if $\alpha(u) \equiv0$.

We shall say that the convergence rates in the discrete model with
sampling at $M$ points $\underline{u}$ are ``inferior'' to the
convergence rates
in the
continuous model if the convergence rates differ by more than a
logarithmic factor for $\alpha(u) \equiv0$ or by more than a constant factor
if $\alpha(u) >0$.
The following statement shows when this happens.
\begin{theorem}
\label{th:necessary}
Let assumptions (\ref{ns_prop1}), (\ref{equi:reg}) and (\ref{k_val})
hold and
let
\[
\lim_{n \rightarrow\infty}\ln(\varepsilon_n)/\ln n =
\varepsilon_0 < \infty
\]
for some sequence
$\varepsilon_n>0$, independent of $m$.

\begin{longlist}
\item
Let $\alpha(u) \equiv0$, and let assumption (\ref{ns_prop2})
hold. If
$M=M_n$ and $\underline{u}$ are such that
%
%
\begin{equation} \label{necc_pol}
\tau_1^d(m,\underline{u},M_n) \leq K \varepsilon_n |m|^{-2 \nu}
({\ln}|m|)^{-\lambda}
\end{equation}
for some constants $\lambda\in\mathbb R$ and $\nu>0$, independent of
$m$ and
$n$, then the convergence rates in the discrete model are
inferior to the convergence rates in the continuous model if
%
%
\begin{eqnarray} \label{nec10}
\nu &>& \nu(u^*) \quad\mbox{and}\nonumber\\[-8pt]\\[-8pt]
\varepsilon_0 &<& \cases{
2\bigl(\nu- \nu(u^*)\bigr)/\bigl(2s+2\nu(u^*)+1\bigr), &\quad if $\nu(2-p) < p{s^*}$,\cr
2\bigl(\nu- \nu(u^*)\bigr)/\bigl(2{s^*}+2\nu(u^*)\bigr), &\quad if $\nu(2-p) \geq
p{s^*}$,}\nonumber
\end{eqnarray}
or
%
%
\begin{equation} \label{nec11}
\nu= \nu(u^*) \quad\mbox{and}\quad
\lim_{n \rightarrow\infty}\varepsilon_n (\ln n)^a =0 \qquad\mbox{for
any } a>0.
\end{equation}

\item Let $\alpha(u) > 0$ and $M=M_n$ and $\underline{u}$ be such that
%
%
\begin{equation} \label{necc_exp}
\tau_1^d(m,\underline{u},M_n) \leq K \varepsilon_n |m|^{-2 \nu}
\exp( -\alpha
|m|^{\beta} )
( {\ln}|m| )^{-\lambda}
\end{equation}
for some constants $\nu\in\mathbb R$, $\lambda\in\mathbb R$ and
$\alpha>0$,
independent of $m$ and $n$.
Then the convergence rates in the discrete model are
inferior to the convergence rates in the continuous model if
%
%
\begin{equation} \label{nec2}
\beta> \beta(u^*) \quad\mbox{and}\quad \varepsilon_0 \geq-1
\quad\mbox{or}\quad
\beta= \beta(u^*) \quad\mbox{and}\quad
\varepsilon_0 = -1.
\end{equation}
\end{longlist}
\end{theorem}

Theorems \ref{th:sufficient} and \ref{th:necessary} formulate
conditions in terms
of $\tau_1^d(m,\underline{u},M_n)$. The following corollaries contain
more specific results for various sampling schemes.
\begin{corollary} \label{cor1}
Let $M=M_n$ be finite. Then the necessary and sufficient condition for
the convergence rates in the continuous and the discrete models to
almost coincide
is that for at least one $l$, $l=1,2, \ldots, M$, one has $\nu(u_l) =
\nu(u^*)$ if $\alpha(u) \equiv0$
or $\beta(u_l) = \beta(u^*)$ if $\alpha(u) >0$.
\end{corollary}
\begin{corollary} \label{cor2}
If $\alpha(u) \equiv0$ and $M=M_n \leq C (\ln n)^{\lambda^*}$ for some
constant $\lambda^* \in[0, \infty)$, then the convergence rates in the
continuous and the discrete models almost coincide if one has $\nu(u_l)
= \nu(u^*)$ for at least one $l$, $l=1,2,\ldots, M$.
\end{corollary}
\begin{corollary} \label{cor3}
If $\alpha(u) >0$ and $M=M_n \leq C n^\tau$ for some constant $\tau
\in
[0,1)$, then the convergence rates in the continuous and the discrete
models almost coincide if one has $\beta(u_l) = \beta(u^*)$ for at
least one $l$, $l=1,2, \ldots, M$.
\end{corollary}

\subsection{Pseudo-uniform sampling strategies}
\label{uniform_strategies}

Theorems \ref{th:sufficient} and \ref{th:necessary} and Corollaries~\ref
{cor1}, \ref{cor2} and \ref{cor3}
in Section \ref{equirates} establish, in the case of an arbitrary
sampling scheme, when the
convergence rates in the continuous model almost coincide
with the convergence rates in the discrete model, or when the
convergence rates in the discrete model are inferior.

However, when the discrete model is replaced by the continuous model,
the underlying implicit assumption is that
sampling is carried out at $M=M_n$ equidistant points with $M_n
\rightarrow\infty$. In particular, the interval $[a,b]$ is partitioned
into $M$
equal subintervals of the length $\Delta= (b-a)/M$ and $u_l = \Delta
(l + d)$, $l=0,1, \ldots, M-1$,
where $d \in[0,1]$ is the parameter which allows one to accommodate
various sampling techniques (e.g.,
$d =0$, $d=1$ or $d=1/2$, respectively, when sampling is carried out at
the left, the right and the middle of each sub-interval).

Below, we study an extension of this sampling scheme. We avoid treating
$u_1, u_2, \ldots, u_M$ as a random sample since this
is not the case in both mathematical physics and signal processing
applications. Instead, in order to accommodate various
sampling strategies, we consider a continuously differentiable function
$S(x)$, $x \in[0, 1]$, such that
$0 \leq s_1 \leq S' (x) \leq s_2 < \infty$ and $S(0)=a$, $S(1) = b$.
Let $d\in[0,1]$, and let
%
%
\begin{equation} \label{sampling}
u_l = S\biggl( \frac{l-1+d}{M} \biggr),\qquad l=1,2, \ldots, M.
\end{equation}
Denote the inverse of $S(u)$ by $q(u) = S^{-1}(u)$, $u \in[a,b]$, and
observe that $q(u)$ is continuously differentiable
in $[a,b]$ with $0 \leq1/s_2 \leq q'(u) \leq1/s_1< \infty$. Many
functions $S(\cdot)$ satisfy these conditions, for example,
$S(x) = a + (b-a) x^h$, where $0<h<\infty$ (the case $h=1$ corresponds
to the uniform sampling).
\begin{theorem}
\label{th:uniform_sample}
Let assumptions (\ref{equi:reg}) and (\ref{k_val}) hold, and let $u_l$,
$l=1,2,\ldots,$ $M$, be defined by (\ref{sampling}) where the function
$S(x)$, $x \in[0,1]$, is continuously differentiable such that
$0 \leq s_1 \leq S' (x) \leq s_2 < \infty$ and $S(0)=a$, $S(1) = b$.
Then the convergence rates in the discrete and the continuous models
almost coincide if, for $M=M_n$,
%
%
\begin{eqnarray} \label{uniform_cond}
\alpha(u) &\equiv& 0 \quad\mbox{and}\quad \lim_{n \rightarrow\infty}M_n^{-1}
\ln n = \tau_1<\infty
\quad\mbox{or}\nonumber\\[-8pt]\\[-8pt]
\alpha(u) &>& 0 \quad\mbox{and}\quad \lim_{n \rightarrow\infty
}M_n^{-1} \ln\ln n = \tau
_2<\infty.\nonumber
\end{eqnarray}
If, moreover, $|g_m(u)|^2 = K |m|^{-2\nu(u)}$ for some continuously
differentiable function $\nu(u)$, $u\in U$,
and also
\[
\lim_{n \rightarrow\infty}M_n^{-1} (\ln n)^{1 + 1/k} =0,
\]
where $k$ is defined in (\ref{k_val}), then the convergence rates in the
discrete and the continuous models coincide up to a constant.
\end{theorem}
\begin{remark} \label{unisample_rates}
Note that if $\alpha(u) >0$ in (\ref{equi:reg}) and $d$ in (\ref{sampling})
is such that $\beta(u_l) = \beta(u^*)$
for some $l$, $l=1,2, \ldots, M$, then a combination of Theorems \ref
{th:sufficient} and \ref{th:uniform_sample}
yields that the convergence rates in the discrete and the continuous
models coincide
for any value of $M= M_n$. Note also that, although conditions (\ref
{uniform_cond}) in Theorem \ref{th:uniform_sample}
are sufficient for the convergence rates in the discrete and the
continuous models to almost coincide, examples in the
next section demonstrate that these conditions are also necessary or
close to being necessary;
if the conditions in (\ref{uniform_cond}), or some slightly weaker
conditions, are violated, then the convergence rates in the discrete
model are inferior to the convergence rates in the continuous model.
\end{remark}


\subsection{Examples revisited}
\label{ex_revisited}

\textit{Example} \ref{ex:2} (\textit{continued}).
Recall that $|g_{m}(u)|^2 \asymp\exp( - 2 \ln(1/u) |m|)$, $u \in
[0,r_0]$, so that $\beta=1$ and
$\alpha(u) = 2 \ln(1/u)$. Hence $u^* = r_0$, and if the discrete model
is sampled entirely at the single point $u^*$, then the convergence rates
in the continuous and the discrete models are given by formulae (\ref
{low:reg_equi}) and (\ref{up:reg_equi}) or
(\ref{low:reg_equi1}) and (\ref{up:reg_equi1}), respectively, and
they coincide.

However, the convergence rates in the discrete and the continuous
models coincide under much weaker conditions.
In fact, if $M_n = O(n^\tau)$ for some constant $\tau\in[0,1)$ and
$u_l = r_1>0$ for at least one $l$, $l=1,2,\ldots,M$, then
$
\tau_1^d(m,\underline{u},M ) \geq K n^{-\tau} \exp( - 2 \ln(1/r_1) |m|),
$
and, by Theorem \ref{th:sufficient}, the convergence rates in the
discrete and the continuous models coincide.
On the other hand, if $u_1 = \cdots= u_{M-1} =0$, $u_M=r_1 >0$ and
$M=M_n \asymp n/\ln n$, then
$
\tau_1^d(m,\underline{u},M ) \asymp n^{-1} \ln n \exp( - 2 \ln
(1/r_1) |m|),
$
and, by Theorem \ref{th:necessary}, the convergence rates
in the discrete model are inferior to the convergence rates in the
continuous model.

Now, consider the case of the pseudo-uniform sampling $u_l =
S((l-1+d)/M)$, $l=1,2, \ldots,
M$,
with $d \in[0, 1]$ and a function $S(x)$, $x\in[0,1]$, satisfying the
assumptions of Section
\ref{uniform_strategies}. We will show that the convergence rates in
the discrete and the continuous models coincide
no matter what the value of $M$ is. To verify this, note that
$
\tau_1^d(m,\underline{u},M ) = M^{-1} \sum_{l=1}^{M } u_l^{2|m|}
\leq r_0^{2|m|}.
$
On the other hand, it is easy to see that since $S((l-1+d)/M) =
S((l-1+d)/M) - S(0) \geq s_1 (l-1+d)/M$, one has
$
\tau_1^d(m,\underline{u},M ) \geq M^{-1} \sum_{l=M/2+1}^{M }
u_l^{2|m|} \geq
M^{-1} \sum_{l=M/2+1}^{M } ( M^{-1} s_1 (l-1+d) )^{2|m|}.
$
Here, $s_1 <1$, due to $S(0)=0$,\break $S(1) = r_0 <1$ and $0<s_1 \leq
S'(x)$, and, therefore,
$
\tau_1^d(m,\underline{u},M ) \geq\break M^{-1} \sum_{l=M/2+1}^{M } (0.5
s_1)^{2|m|} =
0.5 \exp( - 2 |m| \log(2/s_1) ).
$
Since $\ln(2/s_1) >0 $, the convergence rates in the discrete and the
continuous models coincide due to Theorems
\ref{th:lower} and \ref{th:upper}.

We conclude this example with a rather obvious observation. Reducing
the sampling interval from $[0,r_0]$ to
$[r_1, r_0]$, with $r_1 >0$, yields $u_* = r_1$ and Theorem~\ref
{th:condisrates} immediately becomes valid.
For this reason, although $|g_{m}(u)|$ does not satisfy condition (\ref
{equi:reg}) [since $|g_m(0)|=0$], the convergence rates
in the continuous and the discrete models coincide for the majority of
``reasonable'' sampling schemes.
Since, with the restriction $0 < r_1 \leq u$, the problem of the
estimation of the boundary condition
for the Dirichlet problem of the Laplacian on the unit circle simply
reduces to the uniform case,
we can consider the problem as an example of an ``almost uniform'' case
and conclude that
replacing the discrete model by the continuous model is a legitimate choice.

\textit{Example} \ref{ex:4} (\textit{continued}).
Recall that $|g_m(u)|^2 \asymp|m|^{-2u}$, $u \in[a,b]$, so that
$\alpha(u) =0$, $\nu(u) = 2u$,
$k=1$ and $u^* =a$. If $M=M_n = O( (\ln n)^{\lambda^*} )$ for some
constant $\lambda^* \geq0$ and
$u_l =a$ for at least one $l$, $l=1,2, \ldots, M$, then, by
Corollaries~\ref{cor1} and \ref{cor2}, the convergence rates in the discrete and the
continuous models almost coincide. On the other hand, if $u_1=a$ but
$u_l \geq a+d$, $d>0$,
for $l=2, 3,\ldots, M$, and $M=M_n$ is such that
$\lim_{n \rightarrow\infty}M_n (\ln n)^{-\lambda^*} = \infty$ for
any constant $\lambda^* >0$,
then the convergence rates in the discrete model are
inferior to those in the continuous model.

To verify this, note that under the assumptions above
$\tau_1^d (m, \underline{u}, M) \leq K ( M_n^{-1}\times |m|^{-2a} + |m|^{-2
(a+d)} )
\leq K
\max( M_n^{-1} |m|^{-2a}, |m|^{-2 (a+d)} )$.
Now, apply Theorem \ref{th:necessary}, first with $\varepsilon_n = M_n^{-1}$
and $\nu= \nu(u^*)$ and then with
$\varepsilon_n = 1$ and $\nu= \nu(u^*) +d$.

Now, consider the case of the pseudo-uniform sampling $u_l =
S((l-1+d)/M)$, $l=1,2, \ldots,
M$,
with $d \in[0, 1]$ and a function $S(x)$, $x\in[0,1]$, satisfying the
assumptions of Section
\ref{uniform_strategies}. By Theorem \ref{th:uniform_sample}, the
convergence rates in the discrete and the continuous models
coincide up to, at most, a logarithmic factor if $M=M_n$ is such that $
\lim_{n \rightarrow\infty}M_n^{-1} \ln n <\infty$.
If, moreover, $|g_m(u)|^2 = K |m|^{-2u}$ and $M=M_n$ is such that
$\lim_{n \rightarrow\infty}M_n^{-1} (\ln n)^{2} =0$, then the
convergence rates coincide up
to, at most, a constant.
In other words, in each case, the convergence rates in the discrete and the
continuous models almost coincide.

Let us show that the opposite is also true: if $d > 0$ and $M=M_n$ is
such that
%
%
\begin{equation} \label{necc_ex4}
\lim_{n \rightarrow\infty}M_n^{-1} (\ln\ln n)^{-1} \ln n = \infty,
\end{equation}
then the convergence rates in the discrete model are inferior to those
in the continuous model.
For this purpose, note that
$u_1 -a = S(d/M) - S(0) \geq s_1 d/M = d_1/M$, so that
$
\tau_1^d (m, \underline{u}, M) \leq K |m|^{-2(a + d_1/M)} =\break K |m|^{-2a}
\exp({-2 d_1 \ln}|m|/M).
$
Since ${\ln}|m| \asymp\ln n$ in this case, we have\break
$
\tau_1^d (m, \underline{u}, M) \leq K \varepsilon_n |m|^{-2a}
$
with $\varepsilon_n = \exp({- 2 d_1 \ln}|m|/M)$.
Now the fact that the convergence rates in the discrete model are
inferior to those in the continuous model
follows from Theorem \ref{th:necessary} and the observation
that condition (\ref{necc_ex4}) implies condition (\ref{nec11}). The latter
shows that the sufficient conditions
of Theorem \ref{th:uniform_sample}
are very close to being also necessary conditions in this case.

\textit{Example} \ref{ex:5} (\textit{continued}).
Recall that $|g_m(u)|^2 \asymp \exp( - \alpha|m|^u )$, $0 <a \leq
u \leq b < \infty$,
so that $\alpha(u) = \alpha>0$ and $u^* =a$. Note that, by Corollary
\ref{cor3}, if $M=M_n$ is such that
$M_n \leq C n^\tau$ for some constant $\tau\in[0,1)$ and $u_l =a$ for
at least one $l$, $l=1,2,\ldots,M$, then the convergence rates in the discrete
and the continuous models almost coincide. However, if $u_1=a$ and $u_l
\geq a+d$
for $l=2,3,\ldots, M$, and $M=M_n$ is such that $M_n \asymp n/ \ln n$,
then the convergence rates in the discrete model are inferior to those
in the continuous model.
To show this, note that
$
\tau_1^d (m, \underline{u}, M) \leq K [ \exp( -\alpha|m|^{a+d} ) + n^{-1}
\ln n \exp( -\alpha|m|^{a} ) ]
\asymp n^{-1} \ln n \exp( -\alpha|m|^{a} )
$
since, in this case, $|m| \asymp(\ln n)^{1/a}$ and thus
$\exp( -\alpha|m|^{a+d} ) = o ( n^{-1} \ln n \exp( -\alpha|m|^{a} ) )$
as \mbox{$n \rightarrow\infty$}.
Hence, application of Theorem \ref{th:necessary} with $\varepsilon_n
= \ln n/n$
yields that the convergence rates in the discrete model are inferior to
the convergence rates in the continuous model.

Now, consider the case of pseudo-uniform sampling. By Theorem \ref
{th:uniform_sample},
the convergence rates in the discrete and the continuous models
coincide if
$M=M_n$ is such that $\lim_{n \rightarrow\infty}M_n^{-1} \ln\ln n <
\infty$. Moreover, by
Remark \ref{unisample_rates},
the convergence rates in the discrete and the continuous models coincide
whenever $d=0$ in formula (\ref{sampling}), no matter what the value of
$M$ is.

Let us show that, if $d >0$, then the second condition in
(\ref{uniform_cond}) is necessary in order for the convergence rates in
the discrete and the
continuous models to coincide up to at most a constant. For this
purpose, we
assume that $M=M_n$ is such that $\lim_{n \rightarrow\infty}M_n^{-1}
\ln\ln n = \infty$ and
prove that the convergence rates in the discrete model are inferior to the
rates in the continuous model. For this purpose, observe that $u_l \geq
a + d/M$ for every $l$, $l=1,2,\ldots,M$,
so that
$
\tau_1^d (m, \underline{u}, M) \leq K \exp( - \alpha|m|^a e^{
{{d \ln}|m|}/{M}} ).
$
Now, recalling that, in this case, ${\ln}|m| \asymp\ln\ln n$ and $\ln
n^* \asymp\ln n$, and repeating the proof of Theorem \ref{th:lower}
with $\varepsilon_n =1$, we obtain that, for every~$n$, in both the
sparse and
the dense cases
as $n \rightarrow\infty$,
\[
R_n(B_{p,q}^s(A),\underline{u},M_n) \geq C (\ln n)^{-{2{s^*}}/({a
+ d/M_n})}.
\]
Hence, the convergence rates in the discrete case are inferior to those
in the continuous model whenever
\[
\lim_{n \rightarrow\infty}( \ln n )^{-{2{s^*}}/({a + d/M_n}) +
{2{s^*}}/{a}} =
\lim_{n \rightarrow\infty}\exp\biggl( \frac{2 {s^*}d}{a(aM_n + d)} \ln
\ln n \biggr) = \infty,
\]
which is true if $\lim_{n \rightarrow\infty}M_n^{-1} \ln\ln n =
\infty$ and $d>0$.

\textit{Example} \ref{ex:6} (\textit{continuation}). Recall that $|g_m(u)|^2
\asymp
|m|^{- 2 \nu} \exp(- u |m|^\beta)$, $u \in[0,b]$,
and that conditions of Lemma \ref{th:answer} do not hold since $\alpha
(u) =u\geq0$ and $\alpha(0) =0$.
We show that, in this example, the convergence rates in the discrete
and the continuous models do not coincide.
Recall that $u^*=0$ and, due to formulae (\ref{ex6bestpoint}) and
(\ref{ex6ratescont}), Theorem \ref{th:lower}
implies that, as $n \rightarrow\infty$,
%
%
\begin{equation} \label{convr1}\quad
R_n^c (B_{p,q}^s (A)) \geq
\cases{C n^{-{2s}/({2s+2\nu+ \beta+1})},
&\quad if $\nu(2-p) < p{s^*}$,\cr
C \biggl( \dfrac{\ln n}{n} \biggr)^{{2{s^*}}/({2s^*+2\nu+ \beta})},
&\quad
if $\nu(2-p) \geq p{s^*}$,}
\end{equation}
and
%
%
\begin{eqnarray} \label{convr2}
R_n^d (B_{p,q}^s (A)) &\asymp& R_n^d (B_{p,q}^s (A), u^*,
1)\nonumber\\[-8pt]\\[-8pt]
&\geq&
\cases{C n^{-{2s}/({2s+2\nu+1})},
&\quad if $\nu(2-p) < p{s^*}$,\cr
C \biggl( \dfrac{\ln n}{n} \biggr)^{{2{s^*}}/({2s^*+2\nu})},
&\quad if $\nu(2-p) \geq p{s^*}$;}\nonumber
\end{eqnarray}
that is, the convergence rates, in both discrete and continuous models,
are polynomial. However, if one samples the model at $u_l \geq d$,
$l=1,2, \ldots, M$, then
$\tau_1^d (m, \underline{u}, M) \leq C |m|^{- 2 \nu} \exp(- d
|m|^\beta)$
and the convergence rates in the discrete model
are logarithmic; that is, as $n \rightarrow\infty$,
%
%
\begin{equation} \label{convr3}
R_n^d (B_{p,q}^s (A), \underline{u}, M) \geq C (\ln n)^{-
{2{s^*}}/{\beta}}.
\end{equation}

Now, consider the pseudo-uniform sampling strategy $u_l =
S((l-1+d)/M)$, $l=1,2,\ldots,M$, with a continuous differentiable
function $S(x)$, $x \in[0,1]$,
such that $S(0)=0$, $S(1)=b$ and $0 \leq s_1 \leq S' (x) \leq s_2 <
\infty$. Since
$s_1 ((l-1+d)/M) \leq S((l-1+d)/M) \leq s_2 ((l-1+d)/M)$, $l=1,2,\ldots
,M$, one obtains, by direct calculations, that
%
%
\begin{equation} \label{tau_un_ex6}
\frac{K |m|^{- 2 \nu} e^{-s_2 d |m|^\beta/M}}{M ( 1 - e^{-s_2
|m|^\beta/M} ) }
\leq\tau_1^d (m, \underline{u}, M) \leq
\frac{ K |m|^{- 2 \nu} e^{-s_1 d |m|^\beta/M} } {M ( 1 - e^{-s_1
|m|^\beta/M} )}.
\end{equation}
Therefore, for $M=M_n$, the convergence rates in the discrete model
depend on the value of $d$ and
the asymptotic behavior of $|m|^\beta/M_n$. Let us now show that by
choosing different values of $d$ and $M_n$,
one can obtain each of the three convergence rates (\ref
{convr1})--(\ref{convr3}).

If $M_n$ is large (e.g., $M_n \geq C n^{1/(2\nu+ \beta+1)}$), so that
$|m|^\beta/M_n \rightarrow0$ as $n \rightarrow\infty$, then
$1 - e^{-s_i |m|^\beta/M_n} \asymp|m|^\beta/M_n$, $i=1,2$.
Therefore,
$\tau_1^d (m, \underline{u}, M_n) \asymp|m|^{-(2\nu+ \beta)}$
and hence the convergence rates in the discrete and the continuous
models coincide and are given by (\ref{convr1}).

If $M_n$ is small [e.g., $M_n =O(\ln n)$], so that $|m|^\beta/M_n
\rightarrow\infty$ as $n \rightarrow\infty$, then
(\ref{tau_un_ex6}) takes the form
\[
K M_n^{-1} |m|^{- 2 \nu} e^{-s_2 d |m|^\beta/M_n}
\leq\tau_1^d (m, \underline{u}, M_n) \leq
K M_n^{-1} |m|^{- 2 \nu} e^{-s_1 d |m|^\beta/M_n}.
\]
If $M = M_n$ is finite and $d>0$, then, by Theorems \ref{th:lower} and
\ref{th:upper},
the convergence rates in the discrete model are logarithmic, they are
given by the right-hand side of formula
(\ref{convr3}) and are inferior to the convergence rates in the
continuous model.

Finally,
if $M = M_n$ is finite and $d=0$, then the convergence rates in the
discrete model are provided by
the right-hand side of formula (\ref{convr2}) and are superior to
those in
the continuous model.
For moderate values of $M_n$, one can obtain convergence rates in
between (\ref{convr2})
and (\ref{convr3}).

\section{Irregular case: A boxcar-like blurring function} \label{box-car}

Suppose that the blurring function $g(\cdot,\cdot)$ in the
continuous model is of a boxcar-like,
for example,
%
%
\begin{equation}
\label{eq:box-car1} g(u,t) = 0.5 \gamma(u) {\mathbb I}(|t| < u),\qquad
u \in U, t \in T,
\end{equation}
where $\gamma(\cdot)$ is some positive function. In this case,
functional Fourier coefficients
$g_m(\cdot)$ satisfy
%
%
\begin{eqnarray}
\label{boxcartype}
g_0(u)=1 \quad\mbox{and}\quad g_m(u)= (2 \pi m)^{-1} \gamma(u) \sin
(2\pi mu),\nonumber\\[-8pt]\\[-8pt]
\eqntext{m \in\mathbb Z\setminus\{0\}, u \in[a,b].}
\end{eqnarray}
It is easy to see that estimation of the initial speed of a wave on a finite
interval (see Example \ref{ex:3} in Section \ref{some_examples})
leads to $g_m
(\cdot)$ of the form (\ref{boxcartype})
with $\gamma(u) =1$ [see (\ref{waverefl})].

Assume that
%
%
\begin{equation}
\label{bc-gammacond} \gamma_1 \leq\gamma(u) \leq\gamma_2 ,\qquad
u \in[a,b],
\end{equation}
for some $0 < \gamma_1 \leq\gamma_2 < \infty$. [Obviously, this is
true if $\gamma(\cdot)$ is a continuous function.] Under (\ref
{bc-gammacond}), it is easily seen that
%
%
\begin{equation}\label{tauone}
\tau_1^c(m) \asymp m^{-2},
\end{equation}
implying that conditions (\ref{cond1})
and (\ref{cond2}) hold with $\nu=1$ and $\alpha=0$. Consequently, in
this case, using the results of Theorems \ref{th:lower} and \ref
{th:upper}, we can
obtain the corresponding asymptotical minimax lower and upper
bounds for the $L^2$-risk.

Consider now the discrete model. Recall from
Section \ref{sec:intro} that this model can be viewed as a
discretization of the continuous model or as a multichannel
deconvolution problem with $M$ channels where
$n=NM$ denotes the total number of observations and, possibly,
$M=M_n \rightarrow\infty$ as $n\rightarrow\infty$. Note that
multichannel deconvolution with boxcar kernels [i.e., $\gamma(u) =1/u$,
for some \textit{fixed} $u >0$]
is the common problem in many areas of
signal and image processing which include, for instance, LIDAR
remote sensing and reconstruction of
blurred images. LIDAR is a lazer device which emits pulses;
reflections of which are gathered by a telescope aligned with the
lazer [see, e.g., Park, Dho and Kong (\citeyear{PDK97}) and Harsdorf
and Reuter (\citeyear{HR00})]. The return signal is used to determine distance
and the position of the reflecting material. However, if the system
response function of the LIDAR is longer than the time resolution
interval, then the measured LIDAR signal is blurred and the
effective accuracy of the LIDAR decreases. This loss of precision can
be corrected by deconvolution. In
practice, measured LIDAR signals are corrupted by additional noise
which renders direct deconvolution impossible. Moreover, if $M\geq2$
(finite) LIDAR devices are used to recover a signal, then we talk
about a \textit{multichannel} deconvolution problem, leading to the
discrete model described by~(\ref{convdis}).

For any choice of $M$ and selection of points $\underline{u}$,
under (\ref{bc-gammacond}), we easily see that
%
%
\begin{equation}
\label{eq:box-car-taulower}\qquad
\tau_1^d(m; \underline{u},M) =
\frac{1}{M} \sum_{l=1}^M \frac{\gamma^2(u_l)\sin^2(2 \pi m u_l)}{4
\pi^2 m^2} \asymp\frac{1}{m^2 M} \sum_{l=1}^M \sin^2(2 \pi m u_l).
\end{equation}
It follows from (\ref{eq:box-car-taulower}) that for any choice of $M$
and any selection of points $\underline{u}$, we have
%
%
\begin{equation}
\label{eq:fan-bx-1}
\tau_1^d(m; \underline{u},M) \leq K m^{-2}.
\end{equation}
Hence, in this case,
by Theorem \ref{th:lower}, the asymptotical minimax lower bounds for the
$L^2$-risk in this discrete model cannot be lower than the
asymptotical minimax lower bounds for the $L^2$-risk obtained in the
continuous model.

However, it is impossible to find a point $u^* \in[a,b]$, independent
of $m$, such that,
for any $u \in[a,b]$, one has
$\sin^2 (2\pi m u) \leq K \sin^2 (2\pi m u^*)$; in other words, in this
case, Condition \ref{ConditionI}
does not hold and we deal with the \textit{irregular} case here.
It turns out that in the case of a boxcar-type kernel, sampling at any
one point
is not at all the best strategy. Indeed, Johnstone and Raimondo
(\citeyear{JR04})
showed that
in the case of standard deconvolution [$M=1$, $\gamma(u)=1/u$, $u=u^* =
a = b$],
the degree of ill-posedness is $\nu=3/2$. The latter means that the
asymptotical minimax lower bounds for the $L^2$-risk is
given by Theorem \ref{th:lower} with $\alpha=0$ and $\nu= 3/2$.
Johnstone and Raimondo (\citeyear{JR04}) also demonstrated that if
$u^* = a$ is
selected to be a ``Badly
Approximable'' (BA) irrational number, then these lower bounds can be
attained over a wide range of ellipsiods using a nonlinear blockwise
estimator in the sequence space domain.

The convergence rates obtained above can be improved by sampling at
several different
points. De Canditiis and Pensky (\citeyear{CP06}) studied the
multichannel deconvolution
problem with the boxcar blurring function and derived that if $M$
is \textit{finite}, $M \geq2$, one of the $u_1,u_2,\ldots,u_M$ is a BA
irrational number, and $\underline{u}$ is a BA irrational
tuple, then in formula (\ref{delta1})
%
%
\begin{equation} \label{infMbox}
\Delta_1 (j) \leq C(M) j 2^{j(2+ 1/M)}
\end{equation}
[for the definitions of the BA irrational number and the BA irrational
tuple, see, e.g., Schmidt
(\citeyear{S80}), page 42 and also Section \ref{discussion}]. This
implies that
in this case, the degree of ill-posedness is at most $\nu\leq1 + 1/(2M)$,
meaning that if $M>1$, then $\nu$ is less than $3/2$ (that corresponds to
the case of sampling at a single BA irrational number). Furthermore, De
Canditiis and Pensky (\citeyear{CP06})
showed that the asymptotical upper bounds for the error [for the
$L^{r}$-risk, $1< r < \infty$ and for a \textit{fixed} response function
$f(\cdot)$] depend on $M$: the larger the
value of $M$ is the higher the asymptotical convergence rates will
be. Hence, in the multichannel boxcar deconvolution problem, it seems
to be
advantageous to take $M=M_n \rightarrow\infty$ as $n \rightarrow
\infty$ and to choose $\underline{u}$ to be a BA irrational
tuple.
However, the theoretical results obtained De Canditiis and Pensky
(\citeyear{CP06})
cannot be blindly applied to accommodate the case when $M=M_n
\rightarrow\infty$ as $n \rightarrow\infty$; this generalization
requires, possibly, nontrivial results in number theory (see the
discussion in Section \ref{discussion}).

On the other hand, if conditions (\ref{eq:box-car1}) and (\ref
{bc-gammacond}) hold
and $M=M_n \rightarrow\infty$ fast enough as $n\rightarrow
\infty$, then it is not needed to employ BA irrational tuples, as
we reveal below. If $M=M_n \rightarrow\infty$
fast enough as $n\rightarrow\infty$, then deconvolution with a
boxcar-like blurring function in the discrete model can provide
estimators with the same convergence rates
as in the continuous model.
The following statement shows that, if $M=M_n \rightarrow\infty$
fast enough as $n \rightarrow\infty$, then an appropriate
selection of points $\underline{u}$ can secure asymptotic relation
similar to (\ref{tauone})
thus ensuring equal convergence rates in both
the discrete and the continuous models.
\begin{lemma}
\label{prop:case1} Consider $g(\cdot,\cdot)$ to be of the form
(\ref{eq:box-car1}) with $\gamma(\cdot)$ satisfying (\ref{bc-gammacond}),
and let $0 < a < b < \infty$. Let $m \in A_j$, where $|A_j|=c2^j$,
for some $c>0$, with $(\ln n)^{\delta} \leq2^j \leq n^{1/3}$, $j
\geq j_0$, for some $\delta> 0$ and $j_0 \geq0$. Take
$u_l=a+(b-a)l/M$, $l=1,2,\ldots,M$. If $M \geq M_{0n} = (32 \pi/3)
(b-a)n^{1/3}$, then, for $n$ and $|m|$ large enough,
%
%
\begin{equation}
\label{eq:fan-bx-2} \tau_1^d(m; \underline{u},M) \geq K m^{-2}.
\end{equation}
\end{lemma}

Note that Lemma \ref{prop:case1} can be applied if $M = M_n \geq c_0
n^{1/3}$ for some constant $c_0 >0$, independent of $n$. Let $\Delta
= \min(3 c_0/(32 \pi), b-a)$. Set $M=M_n$, $u_l = a+ l \Delta/M$ and
observe that $u_l \in[a,b]$ for $l=1,2, \ldots, M$. Then the following
statements are valid.
\begin{theorem} \label{th:box-car-minimax} Let $\{\phi_{j_0,k}(\cdot
),\psi_{j,k}(\cdot)\}$ be the periodic Meyer wavelet basis discussed
in Section
\ref{estconstr}.
Consider $g(\cdot,\cdot)$ to be of the form (\ref{eq:box-car1}) with
$\gamma(\cdot)$ satisfying (\ref{bc-gammacond}), and let $0 < a < b <
\infty$. Let $R_n^{o} (B_{p,q}^s (A))$ to be either $R_n^{c}
(B_{p,q}^s (A))$ or $R_n^{d} (B_{p,q}^s (A))$.

(Lower bounds). Let $s >\max(0,1/p-1/2)$, $1 \leq p \leq\infty$, $1
\leq q \leq\infty$ and
$A>0$. Then, as $n \rightarrow\infty$,
%
%
\begin{equation} \label{eq:box-car-lowerF}
R_n^{o} (B_{p,q}^s (A)) \geq
\cases{
C n^{-{2s}/({2s+3})}, &\quad if $s
>3(1/p - 1/2)$,\cr
C \biggl( \dfrac{\ln n}{n} \biggr)^{{s'}/({s'+1})}, &\quad if
$s \leq3(1/p - 1/2)$.}
\end{equation}

(Upper bounds). Let $s
>1/p'$, $1 \leq p \leq\infty$, $1 \leq q \leq\infty$ and
$A>0$. Set $\nu=1$ and assume that $M = M_n \geq c_0 n^{1/3}$ for
some constant $c_0 >0$, independent of $n$. Let
$\hat{f}_n^{c}(\cdot)$ be the wavelet estimator defined by
(\ref{fest}), with ${j_0}$ and $J$ given by (\ref{jpower}), and let
$\hat{f}_n^{d}(\underline{u},M,\cdot)$ be the wavelet estimator
defined by (\ref{fest}), evaluated at the points $u_l = a+ l \Delta/M$,
$l=1,2,\ldots,M$, where $\Delta= \min(3 c_0/(32 \pi), b-a)$ and
${j_0}$ and $J$ are given by (\ref{jpower}). Let also
$\hat{f}_n^{o}(\cdot)$ be either $\hat{f}_n^{c}(\cdot)$ or
$\hat{f}_n^{d}(\underline{u},M,\cdot)$. Let $s > 1/p'$, $1 \leq p
\leq\infty$, $1 \leq q \leq\infty$ and $A>0$. Then, as $n
\rightarrow\infty$,
%
%
\begin{eqnarray}
\label{eq:box-car-lowerF1}
&&\sup_{f \in B_{p,q}^s (A)} {\mathbb E}
\|\hat{f}_n^{o}-f\|^2 \nonumber\\[-8pt]\\[-8pt]
&&\qquad\leq
\cases{
C n^{-{2s}/({2s+3})} ( \ln n)^{\varrho}, &\quad if $s
>3(1/p - 1/2)$,
\cr
C \biggl( \dfrac{\ln n}{n} \biggr)^{{s'}/({s'+1})} ( \ln n
)^{\varrho}, &\quad if $s \leq3(1/p - 1/2)$,}\nonumber
\end{eqnarray}
where $\varrho= 3(2/p-1)_+/(2s+3)$ if $s
>3(1/p - 1/2)$, $\varrho= (1 -p/q)_+$ if $s = 3(1/p - 1/2)$ and
$\varrho=0$ if $s < 3(1/p - 1/2)$.
\end{theorem}

\section{A limited simulation study}
\label{sim}

Here we present a limited simulation study in the multichannel
deconvolution model with a boxcar-like blurring function. We assess the
performance of the
suggested block thresholding wavelet estimator (BT) given by
(\ref{fest}), with equispaced selected points $u_l=l/M$,
$l=1,2,\ldots,M$, and compare it to the term-by-term thresholding
wavelet estimator (TT) proposed by De Canditiis and Pensky (\citeyear{CP06})
where the points, $u_l$, $l=1,2,\ldots,M$, were selected such that
one of the $u_l$'s is a BA irrational number, and $u_1,u_2,
\ldots, u_M$ is a BA irrational tuple [see De Canditiis and Pensky
(\citeyear{CP06}), Section 4].

Specifically, we assume that we observe
%
%
\begin{equation}
\label{eq:mod-simFanis} y(u_l, t_i) = \int_T f(x) g(u_l, t_i-x) \,dx
+ \sigma_{l} \varepsilon_{li},\qquad u_l \in U=[0,1], t_i =
i/N,\hspace*{-32pt}
\end{equation}
where
$
g(u_l,t)= (2u_l)^{-1} {\mathbb I}(|t|<u_l), u_l \in U=[0,1],
$
and $\varepsilon_{li}$ are standard Gaussian random variables,
independent for different $l$ and $i$. For simplicity, we assume
that $\sigma_{l}^2=\sigma^2$ for all $l=1,2,\ldots,M$.

The suggested algorithm consists of the following steps:
\begin{enumerate}
\item For each $M=4, 8, 16$, generate $M$ different equispaced
sequences, $y_{li}$
[$=y(u_l, i/N)$], $l=1,2,\ldots,M$, $i=1,2,\ldots, N$, following
model (\ref{eq:mod-simFanis}).

\item Generate functions
$g(u_{l},\cdot)$, $y(u_{l},\cdot)$ $\phi_{j_{0}k}(\cdot)$ and
$\psi_{jk}(\cdot)$, $j=j_0,j_0+1,\ldots,J-1$,
$k=0,1,\ldots,2^{j}-1$, at the same equispaced points, $t_i=i/N$,
$i=1,2,\ldots,N$.

\item Apply the discrete Fourier transform (FFT)
to $g_{l}$, $y_{l}$, $\phi_{j_{0}k}$ and $\psi_{jk}$,
$j=j_0,j_0+1,\ldots,J-1$, $k=0,1,\ldots,2^{j}-1$.

\item Estimate
$a_{j_{0}k}$ and $b_{jk}$ by, respectively,
$\hat{a}_{j_{0}k}$ and $\hat{b}_{jk}$, given by
(\ref{coefest}).

\item
Compute $\hat{B}_{jr}=\sum_{k \in
U_{jr}} \hat{b}_{jk}^{2}$.

\item
Compute the threshold
$
\lambda_{j}=\hat{\sigma}^2 d^* n^{-1} \ln n \Delta_{1}(j), j
\geq j_0$,
where $n=NM$, \mbox{$d^*=1$},
\begin{eqnarray*}
\hat{\sigma} &=& \sqrt{\frac{1}{M(N-2)}\sum_{l=1}^{M}
\sum_{i=2}^{N-1} \biggl(\frac{y_{l,i-1}}{\sqrt{6}}-\frac{2y_{li}}{\sqrt
{6}}+\frac{y_{l,i+1}}{\sqrt{6}} \biggr)^{2}},\\
\Delta_{1}(j)
&=&\frac{1}{|C_{j}|}\sum_{m \in C_{j}}\tau^{-1}_{1}(m)
\end{eqnarray*}
[see Pensky and Sapatinas (\citeyear{PenskyS09}), Remark 6, and M\"
{u}ller and Stadm\"
{u}ller (\citeyear{MS87})].

\item
Threshold the wavelet coefficients belonging to blocks with
$|\hat{B}_{jr}|<\lambda_{j}$.

\item Apply the inverse wavelet
transform to obtain $\hat{f}_n(\cdot)$ given by (\ref{fest}).
\end{enumerate}

%
\begin{figure}

\includegraphics{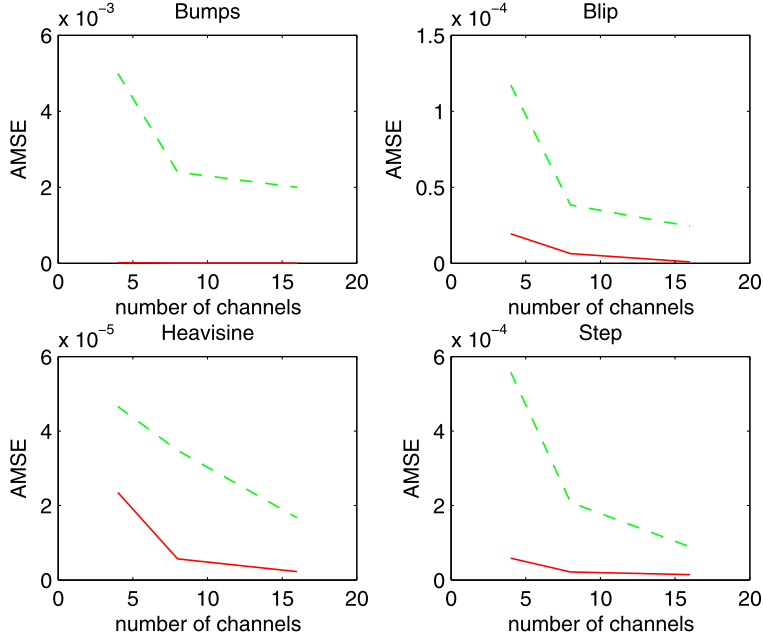}

\caption{AMSE for the \textup{Bumps}, \textup{Blip}, \textup
{Heavisine} and
\textup{Step} functions sampled at a fixed number of $N=128$ points,
based on RSNR${}={}$1, as the number of channels $M$ (and hence the
sample size $n$) increases. Solid line: BT wavelet estimator; Dash
line: TT wavelet estimator.}\label{plotdec1}
\end{figure}

We used the test functions ``Bumps,''
``Blip,'' ``Heavisine'' and ``Step,'' and set $j_0=3$.
For a fixed value of the (root) signal-to-noise ratio (RSNR${}={}$1),
we generated $S=100$ samples of size $n=NM$ from model
(\ref{eq:mod-simFanis}) in order to calculate the average
mean-squared error (AMSE) given by
\[
S^{-1} \sum_{m=1}^{S}\sum_{i=1}^{N}
\bigl(\hat{f}_{n}^m(t_{i})-f(t_{i})\bigr)^{2} \bigg/ \sum_{i=1}^{N}f^{2}(t_{i}),
\qquad
t_i=i/N.
\]
In Figure \ref{plotdec1}, for
a fixed number of data points $N=2^7$, we evaluate the AMSE as the
number of channels $M$, and hence the sample size $n$, increases
for the four signals mentioned above. Obviously, both BT and TT
wavelet estimators improve their performances as $n$ increases,
and the BT wavelet estimator appears to have smaller AMSE than the TT
wavelet estimator in all cases.

Although not reported here, we also evaluated the precision of
the suggested BT wavelet estimator for a wide variety of other
test functions [see the list of test functions in Appendix I of
Antoniadis, Bigot and Sapatinas (\citeyear{ABS01})] and RSNRs with
very good
performances. This numerical study confirms that under the
multichannel deconvolution model with a boxcar-like blurring function,
block thresholding
wavelet estimators with equispaced selection of points $u_l$,
$l=1,2,\ldots,M$, produce quite accurate estimates of $f(\cdot)$.

\section{Concluding remarks}
\label{discussion}

We considered the question of whether and when, in the functional
deconvolution setting, it is legitimate
to replace the real-life discrete deconvolution problem
by its continuous idealization.
In other words, using the asymptotical minimax framework, we studied
whether the continuous model and the discrete model
are equivalent for some or any sampling schemes from the viewpoint of
convergence rates, over a wider range of Besov balls and for the
$L^2$-risk. It is worth mentioning that when we talked about
convergence rates
we referred to the lower bounds which are attainable up to, at most, a
logarithmic factor
according to Theorems \ref{th:lower} and \ref{th:upper}.
In the cases when convergence rates in the discrete model depend on the
choice of a
sampling scheme, we also explored the optimal sampling strategies.
The conclusions of our investigation can be summarized as follows.

If Conditions \ref{ConditionI}, \ref{ConditionIzv} and \ref{ConditionII} are satisfied, then the convergence rates
in the discrete model are independent of the number $M$, and
the choice of sampling points $\underline{u}$ and coincide with the
convergence rates
in the continuous model. In this case, which we call \textit{uniform}, it
is legitimate
to replace discrete model (with any selection of sampling points) by
continuous model.

If Condition \ref{ConditionII} does not hold, then there exist at least
two different sampling schemes in discrete model which deliver two
different sets of convergence rates, and at least one of these
sampling schemes leads to the convergence rates different from the
continuous model.
However, if Condition \ref{ConditionI} holds, one can point out the sampling scheme
which delivers the fastest convergence rates, namely, sampling
entirely at ``the best possible'' point $u^*$. We refer to this case as
\textit{regular}
and explore when, under an arbitrary sampling scheme, convergence rates
in the discrete model
coincide or do not coincide with the convergence rates in the
continuous model. The case of sampling at
$u^*$ is studied as a particular case.

In addition, we consider convergence rates in the discrete model under
uniform or pseudo-uniform sampling strategies.
Indeed, when a discrete model is replaced by its continuous counterpart,
it is implicitly assumed that sampling is carried out at $M$
equidistant points in the interval $[a,b]$.
We formulate conditions when this replacement is legitimate and bring
examples when the uniform, or
a more general pseudo-uniform, sampling may lead to convergence rates
which differ from the convergence rates in the continuous model
and are lower than when sampling is carried out entirely at the ``best
possible'' point $u^*$.
Hence, even in the regular case, one should be extremely careful when
replacing a discrete model
by its continuous counterpart.

Finally, we study the case when Condition \ref{ConditionI} is violated. We referred to
this case as \textit{irregular}.
In this case, the convergence rates in the discrete model depend on a
sampling strategy,
and, in addition, one cannot design a sampling scheme which delivers
the highest convergence rates. Since
Condition \ref{ConditionI} can be violated in a variety of ways, in the irregular case
a general study is very complex.
For this reason, we study a particular example of the irregular case, namely,
functional deconvolution with a boxcar-like blurring function. This important
model occurs, in the problem of estimation of the speed of a
wave on a finite interval (Example \ref{ex:3}) as well as, a discrete
version of
it, in signal and image processing
(see Section \ref{box-car}). In the case of a boxcar-like kernel,
sampling at any one point is, by far,
not the best possible choice and delivers lower convergence rates than
the continuous model.
The best choice for this model is uniform sampling with a large value
of $M=M_n$.
Indeed, if $M = M_n \geq c_0 n^{1/3}$ for some constant $c_0 >0$,
independent of $n$, and the selection
points $u_1,u_2,\ldots,u_M$, are selected to be equispaced,
then, according to Theorem \ref{th:box-car-minimax}, the convergence
rates in the discrete
model with a boxcar-like blurring function coincide with the
convergence rates in the continuous model
and cannot be improved.

The assumption that $M=M_n$ grows at least at a rate of $n^{1/3}$ is
very natural in the
inverse mathematical physics problems: in fact, if one samples uniformly
in the rectangle $[0,1]\times[a,b]$, then $M_n \asymp\sqrt{n}$.
However, this assumption is hardly natural in a signal processing
setting where
$M$ corresponds to a number of physical devices, so even if $M=M_n
\rightarrow\infty$
as $n \rightarrow\infty$, it grows at a very slow rate.
For this reason, the question remains: if $M=M_n \rightarrow\infty$
at a rate slower than $O(n^{1/3})$ [e.g.,
$M=M_n = c_3 n^\upsilon$, where $0< \upsilon< 1/3$, or $M=M_n = c_4
(\ln n)^{\gamma}$, where $\gamma>0$, for some constants $c_3 >0$ and
$c_4 >0$, independent of $n$], can one select points $u_l \in
[a,b]$, $l=1,2, \ldots, M$, such that the convergence rates in the
discrete model coincide with the corresponding
convergence rates obtained in the continuous
model? And, if for some such $M=M_n$ the convergence rates in the
discrete and the continuous models are not the
same, what are the best convergence rates
that can be attained and the best selection of points $u_1, u_2, \ldots
, u_M$?

The solution of this question, possibly, rests on very nontrivial
results in number theory. Recall that De Canditiis and Pensky
(\citeyear{CP06})
showed that, if $M$ is finite, $M \geq2$, one of the $u_l$'s is a
BA irrational number, and $u_1,u_2, \ldots, u_M$ is a BA irrational
tuple, then (\ref{infMbox}) is valid. The constant $C(M)$ in
(\ref{infMbox}) depends on the value of $M$ and the choice of the BA
irrational tuple. Let us now elaborate more on this. Note that the
numbers $a_1,a_2,\ldots,a_M$ is a BA irrational tuple [see, e.g.,
Schmidt (\citeyear{S80}), page 42], if, for any integers
$p_1,p_2,\ldots,p_M$
and $q$, there exists constant $B_M$ such that
\[
\max(|a_1q-p_1|, |a_2q-p_2|,\ldots,|a_Mq-p_1| ) \geq B_M
q^{-1/M},
\]
where $B_M$ is a positive constant that depends only on $M$. Schmidt
[(\citeyear{S80}), page~43] showed that, for a finite value of $M$, a BA
irrational tuple always exists, and proposed an algorithm for
constructing it. It is easy to note that $B_M \rightarrow0$ as $M
\rightarrow\infty$. The value of $B_M$ affects the value of $C(M)$
in (\ref{infMbox}) and, therefore, the convergence rates in the
discrete model.

Unfortunately, we are not aware
of any results in number theory on how $B_M$ depends on $M$, and we
suspect that relevant results may not have yet been derived.
However, a partial answer to the above question, showing that $B_M \geq
C_0 \exp(-3M \ln M)$, for some $C_0 >0$, independent of $M$, $q$ and
$p_1,p_2,\ldots,p_M$, and the construction of minimax upper bounds for
the $L^2$-risk over a wide range of Besov balls, covering the case
$M=M_n=o((\ln n)^{u})$, where
$u \geq1/2$, have been recently obtained in Pensky and Sapatinas
(\citeyear{PS09}).


\begin{appendix}\label{append}
\section*{Appendix: Proofs}

Recall that the symbol $C$ is used for a generic positive constant,
independent of~$n$, while the symbol $K$ is used for a generic positive
constant, independent of $m$, $n$, $M$ and $u_1,u_2,\ldots,u_M$, which
either of them may take different values at
different places.
\begin{pf*}{Proof of Theorem \protect\ref{th:lower}}
The proof of the lower bounds falls into two parts. First, we consider
the lower bounds obtained when the worst functions $f$ (i.e., the
hardest functions to estimate) are represented by only one term in a
wavelet expansion (sparse case), and then when the worst functions $f$
are uniformly spread over the unit interval $T$ (dense case).

In the continuous model, one can always choose $\varepsilon_n=1$,
so the only difference with Pensky and Sapatinas (\citeyear
{PenskyS09}) is an extra
logarithmic factor.
Since the differences for the discrete model are much more significant,
we only consider below the proof for the discrete model.

\textit{Sparse case.} Let
the functions $f_{jk}$ be of the form $f_{jk}= \gamma_j\psi_{jk}$ and let
$f_0 \equiv0$. Note that by (\ref{bpqs}), in order $f_{jk}\in
B_{p,q}^s (A) $,
we need $\gamma_j\leq A 2^{-js'}$. Set $\gamma_j= c 2^{-j s'}$, where $c$
is a positive constant such that $c<A$, and apply the following
classical lemma on lower bounds:
\begin{lemma}[{[H\"{a}rdle et al. (\citeyear{Hardleetal98}), Lemma
10.1]}] \label{korost}
Let $V$ be a functional space, and let $d(\cdot, \cdot)$
be a distance on $V$. For $f, g \in V$, denote by $\varLambda_n
(f,g)$ the likelihood ratio $\varLambda_n (f,g) = d{\mathbb P}_{X_n^{(f)}}/
d{\mathbb P}_{X_n^{(g)}}$, where $d{\mathbb P}_{X_n^{(h)}}$ is the probability
distribution of the process $X_n$ when $h$ is true. Let $V$ contain
the functions $f_0, f_1, \ldots, f_\aleph$ such that
\textup{(a)} $d(f_k, f_{k'}) \geq\delta>0$ for $k=0,1,\ldots,\aleph$, $k
\neq k'$;
\textup{(b)} $\aleph\geq\exp(\lambda_n)$ for some $\lambda_n >0$;
\textup{(c)}~$\ln\varLambda_n (f_0, f_k) = u_{nk} - v_{nk}$, where $v_{nk}$
are constants and $u_{nk}$ is a
random variable such that there exists $\pi_0>0$ with
${\mathbb P}_{f_k}(u_{nk} >0) \geq\pi_0$;
\textup{(d)} $\sup_k v_{nk} \leq\lambda_n$.

Then $\sup_{f \in V} {\mathbb P}_{X_n^{(f)}} (d(\tilde{f}, f) \geq
\delta
/2 ) \geq\pi_0/2 $
for any arbitrary estimator $\tilde{f}$.
\end{lemma}

Let now $V=\{ f_{jk}\dvtx0 \leq k \leq2^j-1 \}$ so that $\aleph=
2^j$. Choose $d(f,g) = \| f-g\|$, where $\| \cdot\|$ is the
$L^2$-norm on the unit interval $T$. Then $d(f_{jk}, f_{jk'}) = \gamma_j=
\delta$. Let $v_{nk} = \lambda_n = j \ln2$ and $u_{nk} = \ln
\varLambda_n (f_0, f_{jk}) + j \ln2$. Now, to apply Lemma
\ref{korost}, we need to show that for some $\pi_0 >0$, uniformly
for all $f_{jk}$, we have
$
{\mathbb P}_{f_{jk}} (u_{nk} > 0) = {\mathbb P}_{f_{jk}} ( \ln
\varLambda_n (f_0,
f_{jk}) > -j \ln2 ) \geq\pi_0 >0.
$

Note that in the case of the discrete model,
\begin{eqnarray*}
- \ln\varLambda_n (f_0, f_{jk}) &=& 0.5 \sum_{i=1}^N\sum_{l=1}^M\{ [y (u_l,
t_i) - \gamma_j(\psi_{jk}*g) (u_l, t_i)]^2 - y^2 (u_l, t_i) \}\\
&=& v_{jk} - u_{jk},
\end{eqnarray*}
where
\[
u_{jk} = \gamma_j\sum_{i=1}^N\sum_{l=1}^M(\psi_{jk}*g)(u_l, t_i)
\varepsilon_{li},\qquad
v_{jk} = 0.5 \gamma_j^2 \sum_{i=1}^N\sum_{l=1}^M[(\psi_{jk}*g)
(u_l, t_i)]^2.
\]
Observe that, due to ${\mathbb P}(\varepsilon_{li} >0) = {\mathbb
P}(\varepsilon_{li} < 0)
=0.5$, we have ${\mathbb P}(u_{jk} >0) = 0.5$.
By properties of the discrete Fourier transform and taking into account that
in the case of Meyer wavelets $|\psi_{mjk}| \leq
2^{-j/2}$ [see, e.g., Johnstone et al. (\citeyear{Johnstoneetal04}),
page 565], we derive that
\[
v_{jk} \leq\frac{\gamma_j^2}{4 \pi} \sum_{i=1}^N\sum_{l=1}^M\sum
_{m \in C_j}|\psi_{mjk}|^2 |g_m(u_l)|^2
\leq\frac{ N M \gamma_j^2}{4 \pi2^j} \sum_{m \in C_j}M^{-1} \sum
_{l=1}^M|g_m(u_l)|^2
\equiv B_n,
\]
where $B_n = (4 \pi)^{-1} n 2^{-j} \gamma_j^2 \sum_{m \in C_j}\tau
_1^d(m,\underline{u},M)$.

Let $j=j_n$ be such that $B_n \leq0.5 j \ln2$. Then, by
applying Lemma \ref{korost} and Chebyshev's inequality, we obtain
%
%
\begin{eqnarray} \label{lowbound1}
\inf_{\tilde{f}_n} \sup_{f \in B_{p,q}^s(A)} {\mathbb E}\|\tilde
{f}_n - f\|^2 &\geq&
\inf_{\hat{f}_n} \sup_{f \in V} \frac{1}{4} \gamma_j^{2}
{\mathbb P}(\|\tilde{f}_n - f\| \ge\gamma_j/2 ) \nonumber\\[-8pt]\\[-8pt]
&\geq& 0.25
\gamma_j^{2} \pi_0.\nonumber
\end{eqnarray}
Thus we just need to choose the smallest possible $j=j_n$
satisfying $B_n \leq0.5 j \ln2$, evaluate $\gamma_j= c 2^{-j s'}$ and
to plug
it into (\ref{lowbound1}). By direct calculations, we derive, under
condition (\ref{cond1}), that
%
%
\begin{equation} \label{err1}
\sum_{m \in C_j}\tau_1^d(m,\underline{u},M) \leq
\cases{
K \varepsilon_n 2^{-j(2\nu-1)} j^{-\lambda}, \qquad\mbox{if $\alpha=0$}, \cr
K \varepsilon_n 2^{-j(2\nu+ \beta-1)} j^{-\lambda} \exp\bigl( -\alpha
(2\pi/3)^\beta
2^{j \beta}\bigr),\cr
\mbox{\hspace*{107pt}if $\alpha>0$}.}
\end{equation}
Hence if $\alpha=0$, then $2^{j_n} = C ( {n^*}(\ln{n^*})^{-(\lambda+1)}
)^{1/(2 s'+ 2\nu)}$,
and if $\alpha>0$, then $2^{j_n} = C (\ln{n^*})^{1/\beta}$.
Now, to obtain the lower bound, plug $\gamma_j= c 2^{-j_n s'}$ into
(\ref{lowbound1})
%
%
\begin{equation}\label{lowbousparse}\quad
\inf_{\tilde{f}_n} \sup_{f \in B_{p,q}^s}
{\mathbb E}\|\tilde{f}_n - f\|^2 \geq
\cases{C ({n^*})^{{2 s'}/({2s'+2\nu})} (\ln{n^*})^{{2
s'(\lambda
+1)}/({2s'+2\nu})}, \cr
\hspace*{91pt}\mbox{if $\alpha
=0$},\cr
C (\ln{n^*})^{-{2 s'}/{\beta}}, \qquad\mbox{if $\alpha>0$}.}
\end{equation}

\textit{Dense case}. Let
$\eta$ be the vector with components $\eta_k = \pm1$, $k=0,1,
\ldots,\break 2^j-1$, denote by $\varXi$ the set of all possible vectors
$\eta$ and let $f_{j \eta} = \gamma_j\sum_{k=0}^{2^j-1}\eta_k \psi
_{jk}$. Let also
$\eta^i$ be the vector with components
$\eta^i_k = (-1)^{{\mathbb I}(i=k)}\eta_k $ for $i,k =
0,1,\ldots,$\break
$2^j-1$. Note
that by
(\ref{bpqs}), in order $f_{j \eta} \in B_{p,q}^s (A) $, we need
$\gamma_j\leq A
2^{-j (s+1/2)}$. Set $\gamma_j= c_{\star} 2^{-j (s+1/2)}$, where
$c_{\star}$ is a positive constant such that $c_{\star}<A$, and
apply the following lemma on lower bounds:
\begin{lemma}[{[Willer (\citeyear{W05}), Lemma 2]}] \label{willer}
Let $\varLambda_n (f,g)$ be
defined as in Lemma \ref{korost}, and let $\eta$ and $f_{j \eta}$ be
as described above. Suppose that, for some positive constants $\lambda$
and $\pi_0$, we have
$
{\mathbb P}_{f_{j \eta}} (- \ln\varLambda_n (f_{j \eta^i}, f_{j
\eta}) \leq
\lambda) \geq\pi_0,
$
uniformly for all $f_{j \eta}$ and all $i=0,\ldots, 2^j-1$. Then,
for any arbitrary estimator $\tilde{f}$ and for some constant $L>0$,
one has
$
\max_{\eta\in\varXi} {\mathbb E}_{f_{j \eta}} \|\tilde{f} - f_{j
\eta}\| \geq L
\pi_0 e^{-\lambda} 2^{j/2} \gamma_j.
$
\end{lemma}

Since, by Chebychev's inequality,
\[
{\mathbb P}_{f_{jk}} \bigl( \ln\varLambda_n (f_{j \eta^i}, f_{j \eta}) >
- \lambda
\bigr)
\geq1 - {\mathbb E}_{f_{jk}} |{\ln\Lambda_n} (f_{j \eta^i}, f_{jk}) |
/ \lambda,
\]
we need to show that
${\mathbb E}_{f_{j \eta}} |{\ln\varLambda_n} (f_{j \eta^i}, f_{j \eta
})| \leq
\lambda_1$,
for a sufficiently small constant $\lambda_1 >0$. Observe that
%
\begin{eqnarray*}
\ln\varLambda_n (f_{j \eta^i}, f_{j \eta})
&=& 0.5 \gamma_j^2 \sum_{i=1}^N\sum_{l=1}^M[(g* f_{j \eta^i} - f_{j
\eta} ) (u_l, t_i)]^2\\
&&{}- \gamma_j \sum_{i=1}^N\sum_{l=1}^M\varepsilon_{li} \bigl[ (g* [f_{j
\eta^i} - f_{j \eta}])
(u_l, t_i) \bigr].
\end{eqnarray*}
Then, due to $|f_{j \eta^i}- f_{j \eta})| = 2|\psi_{jk}|$, one has
${\mathbb E}_{f_{j \eta}} |{\ln\varLambda_n} (f_{j \eta^i}, f_{j \eta
})| \leq
A_n + B_n$
where
\begin{eqnarray*}
A_n &=& 2 \gamma_j {\mathbb E}\Biggl| \sum_{i=1}^N\sum_{l=1}^M(\psi
_{jk}*g)(u_l, t_i) \varepsilon_{li} \Biggr|,\\
B_n &=& 2 \gamma_j^2 \sum_{i=1}^N\sum_{l=1}^M(\psi_{jk}*g)^2(u_l, t_i).
\end{eqnarray*}
Since, by Jensen's inequality, $A_n \leq\sqrt{2 B_n}$, we only need
to construct an upper bound for $B_n$. Note that, similarly to the
sparse case,
one has $B_n = O( n 2^{-j} \gamma_j^2 \sum_{m \in C_j}\tau
_1^d(m,\underline{u},M))$.
%
According to Lemma \ref{willer}, we choose $j = j_n$ that satisfies
the condition $B_n + \sqrt{2 B_n} \leq\lambda_1$.
Using (\ref{err1}), we derive that $2^{j_n} = C ( {n^*}(\ln{n^*})^{-
\lambda}
)^{1/(2s+2\nu+1)}$ if
$\alpha=0$ and $2^{j_n} = C (\ln{n^*})^{1/\beta}$ if $\alpha>0$.
Then, Lemma~\ref{willer} and Jensen's inequality yield
%
%
\begin{equation}\label{lowboudense}\qquad
\inf_{\tilde{f}_n} \sup_{f \in B_{p,q}^s} {\mathbb E}\|\tilde{f}_n
- f\|^2
\geq
\cases{C ({n^*})^{-{2s}/({2s+ 2\nu+1})} (\ln{n^*})^{-{2s
\lambda}/({2s+
2\nu+1})}, \cr
\hspace*{88.8pt}\mbox{if $\alpha
=0$},\cr
C (\ln{n^*})^{-{2s}/{\beta}}, \qquad\mbox{if $\alpha>0$}.}
\end{equation}

Now, to complete the proof
one just needs to note that ${s^*}= \min(s,s')$, and that
%
%
\begin{equation} \label{cases}
2s/(2s+2\nu+1) \leq2s^*/(2s^* + 2\nu) \qquad\mbox{if }
\nu(2-p) \leq ps^*,
\end{equation}
with the equalities taken
place simultaneously, and then to choose the highest of the lower
bounds (\ref{lowbousparse}) and (\ref{lowboudense}). This completes the
proof of Theorem \ref{th:lower}.
\end{pf*}
\begin{pf*}{Proof of Lemma \protect\ref{l:coef}}
In what follows, we shall only construct the proof for $b_{jk}$ since the
proof for $a_{j_0k}$ is very similar. Again, we construct the proof only
for discrete model, since in the case of continuous model, one can
always choose \mbox{$\varepsilon_n=1$}, so the only difference with Pensky and
Sapatinas (\citeyear{PenskyS09}) is an extra logarithmic factor.

Note that, by (\ref{coefest}), one has
$
\hat{b}_{jk}- b_{jk}= \sum_{m \in C_j} (\hat{f}_m - f_m)
\overline{\psi_{mjk}},
$
with
%
%
\begin{equation} \label{discerror}
\hat{f}_m - f_m= N^{-1/2} \Biggl( \sum_{l=1}^M
\overline{g_{m}(u_l)} z_{ml} \Biggr) \Bigg/ \Biggl( \sum_{l=1}^M|g_{m}(u_l)|^2
\Biggr),
\end{equation}
where $z_{ml}$ are standard (complex-valued) Gaussian random variables, independent for
different $m$ and $l$. Therefore,
\[
{\mathbb E}|\hat{b}_{jk}- b_{jk}|^2 = N^{-1} \sum_{m \in C_j}
|\psi_{mjk}|^2 \Biggl[ \sum_{l=1}^M
|g_{m}(u_l)|^2 \Biggr]^{-1} = O ( n^{-1} \Delta_1 (j) )
\]
since $|C_j| = 4 \pi2^j $ and $|\psi_{mjk}|^2 \leq2^{-j}$.
If $\kappa=2$, then
%
%
\begin{eqnarray}\label{hbbrisk}
{\mathbb E}|\hat{b}_{jk}- b_{jk}|^4 &=& O \biggl( \sum_{m \in C_j}
{\mathbb E}|\hat{f}_m - f_m
|^4 |\psi_{mjk}|^4 \biggr)\nonumber\\
&&{} + O \biggl( \biggl[ \sum_{m \in C_j} {\mathbb E}|\hat{f}_m - f_m|^2 |\psi_{mjk}|^2
\biggr]^2 \biggr) \nonumber\\[-8pt]\\[-8pt]
& = & O \bigl( 2^{-j} N^{-2} M^{-3} \Delta_2 (j) +
N^{-2} M^{-2} \Delta_1^2 (j) \bigr) \nonumber\\
&=& O \bigl( 2^{-j} n^{-2} M^{-1}
\Delta_2 (j) + n^{-2} \Delta_1^2 (j) \bigr).\nonumber
\end{eqnarray}
Direct calculations show that when $\alpha=0$ one has
$\Delta_2 (j) = O (2^{6j\nu} j^{3 \lambda} \varepsilon_n^{-3})$.
Plugging expressions for $\Delta_1 (j)$ and $\Delta_2 (j)$ into formula
(\ref{hbbrisk})
and taking into account that $2^j \leq2^{J-1} < ({n^*})^{1/(2\nu+1)}$,
one derives
%
%
\begin{eqnarray*}
{\mathbb E}|\hat{b}_{jk}- b_{jk}|^4 &=& O \biggl( \frac{2^{6j\nu} j^{3
\lambda}}{n^2 \varepsilon^3 M_n }
+ \frac{2^{4j\nu} j^{2 \lambda}}{n^2 \varepsilon^2} \biggr)
\\
&=& O \bigl( n ({n^*})^{{6 \nu}/({2\nu+1}) - 3} (\ln n)^{3 \lambda}
+ ({n^*})^{{4 \nu}/({2\nu+1}) - 2} (\ln n)^{2 \lambda} \bigr).
\end{eqnarray*}
To complete the proof, observe that in the last expression, the second
term is
asymptotically smaller than the first.
\end{pf*}
\begin{pf*}{Proof of Lemma \protect\ref{l:deviation}}
Again we carry out the proof only for the discrete case. The proof for
the continuous case can be obtained as
a minor variation of the proof below. Consider the set of vectors
$
\varOmega_{jr} = \{ v_k, k \in U_{jr}\dvtx\sum_{k \in U_{jr}}|v_k|^2
\leq1 \}
$
and the centered Gaussian process defined by
$
Z_{jr} (v) = \sum_{k \in U_{jr}}v_k (\hat{b}_{jk}- b_{jk}).
$
The proof of the lemma is based on the following inequality:
\begin{lemma}[{[Cirelson, Ibragimov and Sudakov (\citeyear{CIS76})]}]
\label{l:cirel}
Let $D$ be a subset
of $\mathbb R=(-\infty,\infty)$, and let $(\xi_t)_{t \in D}$ be a
centered Gaussian process. If ${\mathbb E}(\sup_{t \in D} \xi_t )
\leq B_1$
and $\sup_{t \in D}\operatorname{Var} (\xi_t) \leq B_2$, then, for
all $x>0$,
we have
$
{\mathbb P}( \sup_{t \in D} \xi_t \geq x + B_1 ) \leq
\exp(-x^2/(2 B_2) ).
$
\end{lemma}

To apply Lemma \ref{l:cirel}, we need to find $B_1$ and $B_2$. Note
that, by Jensen's inequality, we obtain
\begin{eqnarray*}
{\mathbb E}\Bigl[ \sup_{v \in\varOmega_{jr}} Z_{jr} (v) \Bigr] &=&
{\mathbb E}\biggl[ \sum_{k \in U_{jr}}|\hat{b}_{jk}- b_{jk}|^2 \biggr]^{1/2}
\leq\biggl[ \sum_{k \in U_{jr}}
{\mathbb E}|\hat{b}_{jk}- b_{jk}|^2 \biggr]^{1/2}
\\
&\leq&\frac{\sqrt{c_1} 2^{\nu j} j^{\lambda/2} \sqrt{\ln n}}{\sqrt{{n^*}}}
= B_1.
\end{eqnarray*}
[Here $c_1$ is the same positive constant as in (\ref{delta1}) with
$\alpha=0$.] Also, by (\ref{finaleqdis}) and (\ref{discerror}),
we have
$
{\mathbb E}[(\hat{b}_{jk}- b_{jk})(\hat{b}_{jk'} - b_{jk'}) ]
= n^{-1}
\sum_{m \in C_j}\psi_{mjk}
\overline{\psi_{mjk'}} [\tau_1 (m)]^{-1}
$
where $\tau_1 (m)$ is defined in (\ref{taum}). Hence
\begin{eqnarray*}
\sup_{v \in\varOmega_{jr}} \operatorname{Var}(Z_{jr} (v)) & = &
n^{-1} \sup_{v \in
\varOmega_{jr}} \sum_{k \in U_{jr}}\sum_{k' \in U_{jr}} v_k v_{k'}
\sum_{m \in C_j}\psi_{mjk}\overline{\psi_{mjk'}} [\tau_1 (m)]^{-1}
\\
& \leq & c_1 ({n^*})^{-1} 2^{2 \nu j} j^\lambda\sum_{k \in U_{jr}}v^2_k
\leq c_1 ({n^*}
)^{-1} 2^{2 \nu j} j^\lambda= B_2,
\end{eqnarray*}
by using $\sum_{m \in C_j}\psi_{mjk}\overline{\psi_{mjk'}} =
{\mathbb I}(k=k')$ and (\ref{delta1}) for $\alpha=0$.
Therefore, by applying Lemma \ref{l:cirel} with $B_1$ and $B_2$
defined above and
$
x = B_1 ( (2 \sqrt{c_1})^{-1} \mu\sqrt{h_2} - 1 ),
$
and noting that under condition (\ref{mu_cond}), $\ln({n^*}) \geq h_2
\ln
n$, we derive
\[
{\mathbb P}\biggl( \sum_{k \in U_{jr}}|\hat{b}_{jk}- b_{jk}|^2  \geq
\frac{\mu^2 2^{2\nu j} j^\lambda \ln({n^*})} {4 {n^*}} \biggr)
\leq\exp\biggl\{ - \biggl( \frac{\mu\sqrt{h_2}}{2 \sqrt{c_1}} - 1 \biggr)^2 \frac
{B_1^2}{2B_2} \biggr\}
\leq n^{-3},
\]
since (\ref{mu_cond}) implies that $0.5 [\mu\sqrt{h_2}/ (2 \sqrt
{c_1}) -
1]^2 \geq3$.
This completes the proof of Lemma \ref{l:deviation}.
\end{pf*}
%
\begin{pf*}{Proof of Theorem \protect\ref{th:upper}}
First, note that in the case of $\alpha>0$, we have ${\mathbb E}\|
\hat{f}_n- f \|^2 =
R_1 + R_2$, where
%
%
\begin{equation}\label{r1r2}
R_1 = \sum_{j=J}^\infty\sum_{k=0}^{2^j-1}b_{jk}^2,\qquad
R_2 = \sum_{k=0}^{2^{j_0}-1} {\mathbb E}(\hat{a}_{j_0k}- a_{j_0k})^2,
\end{equation}
since ${j_0}= J$. It is well known [see, e.g., Johnstone (\citeyear
{J02}), Lemma
19.1] that if $f \in B_{p,q}^s (A) $, then for some positive constant
$c^{\star}$, dependent on $p$, $q$, $s$ and $A$ only, we have
%
%
\begin{equation}\label{besball}
\sum_{k=0}^{2^j-1}b_{jk}^2 \leq c^{\star} 2^{-2 j {s^*}};
\end{equation}
thus
$R_1 = O ( 2^{-2 J {s^*}} ) = O ( (\ln{n^*})^{-2{s^*}/\beta} )$.
Also, using (\ref{delta1}) and (\ref{ha}), we derive
$
R_2 = O ( n^{-1} 2^{j_0}\Delta_1 ({j_0}) ) = O ( ({n^*})^{-1/2}
(\ln{n^*})^{2\nu/\beta} ) = o ( (\ln{n^*})^{-2{s^*}/\beta} ),
$
thus completing the proof for $\alpha>0$.

Now consider the case of $\alpha=0$. Note that by condition (\ref{ns_prop2})
one has $\ln{n^*}\asymp\ln n$.
Due to the orthonormality of
the wavelet basis, we obtain
%
%
\begin{equation}\label{errtotal}
{\mathbb E}\| \hat{f}_n- f \|^2 = R_1 + R_2 + R_3 + R_4,
\end{equation}
where $R_1$ and $R_2$ are defined in (\ref{r1r2}), and
\begin{eqnarray*}
R_3 & = & \sum_{j=j_0}^{J-1}\sum_{r \in A_j}\sum_{k \in
U_{jr}}{\mathbb E}\bigl[ (\hat{b}_{jk}- b_{jk})^2 {\mathbb I}\bigl(\hat
{B}_{jr}\geq\mu
^2 ({n^*})^{-1} 2^{2 \nu j} \ln({n^*}) j^\lambda\bigr) \bigr], \\
R_4 & = & \sum_{j=j_0}^{J-1}\sum_{r \in A_j}\sum_{k \in
U_{jr}}{\mathbb E}\bigl[ b_{jk}^2 {\mathbb I}\bigl(\hat{B}_{jr}< \mu^2 ({n^*}
)^{-1} 2^{2 \nu j} \ln({n^*}) j^\lambda\bigr) \bigr],
\end{eqnarray*}
where $\hat{B}_{jr}$ and $\mu$ are defined by (\ref{bjr}), (\ref
{lamj}) and (\ref{mu_cond}), respectively.

Let us now examine each term in (\ref{errtotal}) separately. Similarly
to the case of $\alpha>0$, we obtain
$R_1 = O ( 2^{-2 J {s^*}} ) = O ( n^{-2{s^*}/(2 \nu+1)}
)$.
By direct calculations, one can check that
$2{s^*}/(2 \nu+1) > 2s/(2s + 2\nu+1)$, if
$\nu(2-p) < p {s^*}$, and $2{s^*}/(2 \nu+1) \ge2{s^*}/(2{s^*}
+ 2\nu)$, if $\nu(2-p) \geq p {s^*}$. Hence
%
%
\begin{equation} \label{r1a}
R_1 =
\cases{O \bigl( ({n^*})^{-{2s}/({2s+2\nu+1})} \bigr), &\quad if $\nu(2-p) < p {s^*}$,
\cr
O \bigl( ({n^*})^{-{2{s^*}}/({2{s^*}+2\nu})} \bigr), &\quad if
$\nu(2-p) \geq p {s^*}$.}
\end{equation}
Also, by (\ref{ha}) and (\ref{delta1}), we obtain
%
%
\begin{eqnarray} \label{r2}
R_2 &=& O \bigl( ({n^*})^{-1} 2^{(2 \nu+1) {j_0}} \bigr) =
o \bigl( ({n^*})^{-{2s}/({2s+2\nu+1})} \bigr) \nonumber\\[-8pt]\\[-8pt]
&=& o \bigl( ({n^*})^{-
{2s^*}/({2s^*+2\nu
})} \bigr).
\nonumber
\end{eqnarray}

Denote
\[
\Theta_{jr}= \biggl\{ \omega\dvtx\sum_{k \in U_{jr}}|\hat{b}_{jk}- b_{jk}|^2
\geq0.25 \mu^2 ({n^*})^{-1} 2^{2 \nu j} \ln({n^*}) j^\lambda \biggr\}.
\]
To construct the upper bounds for $R_3$ and $R_4$, note that simple
algebra yields
$
R_3 \leq(R_{31} + R_{32}), R_4 \leq(R_{41} + R_{42}),
$
where
\begin{eqnarray*} R_{31} & = & \sum_{j=j_0}^{J-1}
\sum_{r \in A_j}\sum_{k \in U_{jr}}{\mathbb E}[ (\hat{b}_{jk}-
b_{jk})^2 {\mathbb I}( \Theta_{jr} ) ],\\
R_{41}  & =  & \sum_{j=j_0}^{J-1}\sum_{r \in A_j}\sum_{k \in
U_{jr}}{\mathbb E}[ b_{jk}^2 {\mathbb I}(\Theta_{jr}) ],
\\
R_{32} & = & \sum_{j=j_0}^{J-1}\sum_{r \in A_j}\sum_{k \in
U_{jr}}{\mathbb E}\bigl[ (\hat{b}_{jk}- b_{jk})^2 {\mathbb I}\bigl(B_{jr}
> 0.25 \mu^2 ({n^*})^{-1} 2^{2 \nu j} \ln({n^*}) j^\lambda\bigr) \bigr],
\\
R_{42} & = & \sum_{j=j_0}^{J-1}\sum_{r \in A_j}\sum_{k \in
U_{jr}}{\mathbb E}\bigl[ b_{jk}^2 {\mathbb I}\bigl(B_{jr}
< 2.5 \mu^2 ({n^*})^{-1} 2^{2 \nu j} \ln({n^*}) j^\lambda\bigr) \bigr].
\end{eqnarray*}
%
Then, by (\ref{besball}), Lemmas \ref{l:coef} and \ref{l:deviation},
and the
Cauchy--Schwarz inequality, we derive
\begin{eqnarray*}
R_{31}+ R_{41}
& = & O \Biggl( \sum_{j=j_0}^{J-1}\sum_{r \in A_j}\sum_{k \in U_{jr}}\bigl[ \sqrt
{{\mathbb E}(\hat{b}_{jk}- b_{jk})^4} + b_{jk}^2 \bigr]
\sqrt{ {\mathbb P}(\Theta_{jr}) } \Biggr) \\
& = & O \Biggl( \sum_{j=j_0}^{J-1}\bigl[\sqrt{n} (\ln n)^{3 \lambda/2}
({n^*})^{-{3}/({2(2
\nu+1)})} + 2^{-2j{s^*}}\bigr]
n^{-{3/2}} \Biggr) \\
&=& O(({n^*})^{-1}),
\end{eqnarray*}
provided $\mu$ satisfies (\ref{mu_cond}). Hence
%
%
\begin{equation}\label{r3141}
\varDelta_1 = R_{31}+ R_{41} = O ( ({n^*})^{-1} ).
\end{equation}

Now, consider
%
%
\begin{equation}\label{del2}
\varDelta_2= R_{32} + R_{42}.
\end{equation}
First, let us study the dense case, that is, when $\nu(2-p) < p{s^*}$.
Let $j_1$ be such that
%
%
\begin{equation}\label{j1}
2^{j_1} = ({n^*})^{{1}/({2s + 2\nu+1})} (\ln n)^{({(2/p-1)_+ -
\lambda})/({2\nu+ 2s +1})}.
\end{equation}
Then, $\varDelta_2$ can be partitioned as $\varDelta_2 =
\varDelta_{21} + \varDelta_{22}$, where the first component is
calculated over the set of indices ${j_0}\leq j \leq j_1$ and the
second component over $j_1+1 \leq j \leq J-1$. Hence, using (\ref{bjr})
and Lemma \ref{l:coef}, and taking into account that the cardinality
of $A_j$ is $|A_j|=2^j/\ln n$, we obtain
%
%
\begin{eqnarray} \label{del21c1}
\varDelta_{21}
&=& O \Biggl( \sum_{j=j_0}^{j_1}\Biggl[ \frac{2^{(2\nu+1)j} j^\lambda}{{n^*}} +
\sum_{r \in A_j}\frac{
2^{2 \nu j} \ln({n^*}) j^\lambda}{{n^*}} \Biggr] \Biggr)
\nonumber\\[-8pt]\\[-8pt]
&=& O \Biggl( \Biggl[ \frac{(\ln n)^\lambda}{{n^*}} \Biggr]^{ {2s}/({2s + 2\nu+1})}
(\ln
n)^{\varrho} \Biggr),\nonumber
\end{eqnarray}
where $\varrho$ is defined in (\ref{rovalue}).
To obtain an expression for $\varDelta_{22}$, note that for $p \geq
2$, by
(\ref{ns_prop2}) and (\ref{besball}), we have
%
%
\begin{eqnarray} \label{del22c0}
\varDelta_{22}
&=& O \Biggl( \sum_{j=j_1 +1}^{J-1}\sum_{r \in A_j}B_{jr}\Biggr) = O \Biggl( \sum_{j=j_1
+1}^{J-1}
2^{-2js} \Biggr) \nonumber\\[-8pt]\\[-8pt]
&=& O \bigl( ({n^*})^{-{2s}/({2s+2\nu+1})} (\ln n)^{
{2s\lambda}/({2s+2\nu+1})} \bigr).\nonumber
\end{eqnarray}
If $1 \leq p <2$, then
$
B_{jr}^{p/2} = ( \sum_{k \in U_{jr}}b_{jk}^2 )^{p/2} \leq\sum_{k \in
U_{jr}}|b_{jk}|^p,
$
so that by Lemma \ref{l:coef}, and since $\nu(2 - p) < p s^*$, we
obtain
%
%
\begin{eqnarray} \label{del22c1}
\varDelta_{22} &=&
O \Biggl( \sum_{j=j_1 +1}^{J-1}\sum_{r \in A_j}[ ( ({n^*})^{-1} 2^{2\nu j}
j^\lambda\ln n )^{1-p/2}
B_{jr}^{p/2} ] \Biggr)
\nonumber\\
%
& = & O \Biggl( \sum_{j=j_1 +1}^{J-1}( ({n^*})^{-1} 2^{2\nu j} j^\lambda\ln n
)^{1-p/2} 2^{-pjs^*} \Biggr)
\\
&=& O \bigl( ({n^*})^{-{2s}/({2s+2\nu+1})} (\ln n)^{ {2s\lambda
}/({2s+2\nu+1})
+ \varrho} \bigr).\nonumber
%
\end{eqnarray}

Let us now study the sparse case when $\nu(2-p) > p{s^*}$.
Let $j_1$ be defined by
$2^{j_1} = ({n^*})^{{1}/({2s + 2\nu+1})} (\ln n)^{{ - \lambda
}/({2\nu
+ 2s
+1})}$.
Hence, if $B_{jr}\geq0.25 \mu^2 ({n^*})^{-1} 2^{2 \nu j} \times\break\ln({n^*})
j^\lambda$,
then $B_{jr}\leq\sum_{k=0}^{2^j-1}b_{jk}^2 \leq c^* 2^{-2j{s^*}}$
[see (\ref{besball})]
implies that
$j \leq j_2$ where $j_2$ is such that $2^{j_2} = C [{n^*}/(\ln n)^{1+
\lambda
}]^{1/(2{s^*}+ 2\nu)}$,
where $C$ depends on $\mu$ and $c^*$ only.
Again, partition $\varDelta_2 = \varDelta_{21} + \varDelta_{22}$,
where the
first component is calculated over ${j_0}\leq j \leq j_2$ and the
second component over $j_2+1 \leq j \leq J-1$. Then, using similar
arguments to that in (\ref{del22c1}), and taking into account that
$\nu(2 - p)>p {s^*}$, we derive
%
%
\begin{eqnarray} \label{del21c2}
\varDelta_{21} & = & O \Biggl( \sum_{j=j_0}^{j_2}[ ({n^*})^{-1} 2^{2\nu j}
j^\lambda\ln n ] ^{1-p/2}
\sum_{r \in A_j}\sum_{k \in U_{jr}}|b_{jk}|^p \Biggr) \nonumber\\
& = & O \Biggl( \sum_{j=j_0}^{j_2}[ ({n^*})^{-1} 2^{2\nu j} j^\lambda\ln n ]^{1-p/2}
2^{-pj{s^*}} \Biggr)
\\
&=& O \bigl( [ ({n^*})^{-1} (\ln n)^{1 + \lambda} ]^{{2{s^*}}/({2s^*+2\nu
})} \bigr).\nonumber
\end{eqnarray}
To obtain an upper bound for $\varDelta_{22}$, recall (\ref{del2}) and
keep in mind that the portion of $R_{32}$ corresponding to $j_2+1
\leq j \leq J-1$ is just zero. Hence, by (\ref{besball}), we obtain
\begin{eqnarray*}
\varDelta_{22} &=& O \Biggl( \sum_{j=j_2 +1}^{J-1}\sum_{k=0}^{2^j-1}
b_{jk}^2 \Biggr) = O \Biggl( \sum_{j=j_2 +1}^{J-1}2^{-2j{s^*}} \Biggr)
\\
&=& O \bigl( [ ({n^*})^{-1} (\ln n)^{1 + \lambda} ]^{{2{s^*}}/({2s^*+2\nu
})} \bigr).
\end{eqnarray*}

Now, in order to complete the proof, we just need to study the case
when $\nu(2-p) = p{s^*}$. In this situation, we have $2s/(2s+ 2\nu
+ 1) = 2s^* /(2s^*+ 2\nu) = 1-p/2$ and $2 \nu j (1 - p/2) = p j
s^*$. Recalling (\ref{bpqs}) and noting that $s^* \leq s'$, we obtain
$
\sum_{j=j_0}^{J-1}( 2^{p j s^*} \sum_{k=0}^{2^j-1}|b_{jk}|^p )^{q/p}
\leq A^q.
$
Then we repeat the calculations in (\ref{del21c2}) for all indices
${j_0}\leq j \leq J-1$. If $1 \leq p < q$, then, by H\"{o}lder's
inequality, we obtain
%
%
\begin{eqnarray}
\label{del2c3}
\varDelta_{2}
& = & O \Biggl( ( ({n^*})^{-1} (\ln n)^{1 + \lambda} )^{1-p/2}
(\ln n)^{1 - p/q}
\nonumber\\
&&\hspace*{14.2pt}{} \times
\Biggl[ \sum_{j=j_0}^{J-1}\Biggl( 2^{p j s^*} \sum_{k=0}^{2^j-1}|b_{jk}|^p \Biggr)^{q/p}
\Biggr]^{p/q} \Biggr) \\
& = & O \bigl( ( ({n^*})^{-1} (\ln n)^{1 + \lambda} )^{
{2{s^*}}/({2s^*+2\nu})} (\ln n)^{1 - p/q} \bigr).\nonumber
\end{eqnarray}
If $1 \leq q \leq p$, then, by the inclusion $B_{p,q}^s(A) \subset
B_{p,p}^s(A)$, we obtain
%
%
\begin{eqnarray} \label{del2c4FF}
\varDelta_{2} &=& O \Biggl( \sum_{j=j_0}^{J-1}\bigl( (\ln n)^{1 + \lambda}/{n^*}\bigr)^{1-p/2}
2^{p j s^*} \sum_{k=0}^{2^j-1}|b_{jk}|^p \Biggr)
\nonumber\\[-8pt]\\[-8pt]
&=& O \bigl( \bigl( (\ln n)^{1 + \lambda}/{n^*}\bigr)^{{2{s^*}}/({2s^*+2\nu
})} \bigr).\nonumber
\end{eqnarray}
By combining (\ref{r1a}), (\ref{r2}), (\ref{r3141}),
(\ref{del21c1})--(\ref{del2c4FF}), we complete the proof of Theorem
\ref
{th:upper}.
\end{pf*}
\begin{pf*}{Proof of Theorem \protect\ref{th:condisrates}}
The first part of the theorem is identical to Proposition 1 of Pensky
and Sapatinas (\citeyear{PenskyS09}).
The second part can be proved by contradiction. Assume that,
assumptions (\ref{disconHigh}) and (\ref{disconLow}) hold but condition
(\ref{th3cond}) does not take place. It follows from (\ref{disconHigh})
and (\ref{disconLow}) that
\begin{eqnarray*}
|g_{m}(u_*)|^2 &\leq& K |m|^{-2\nu_1} \exp(-\alpha_1
|m|^{\beta_1}), \qquad\nu_1>0 \mbox{ if } \alpha_1 =0,
\\
|g_{m}(u^*)|^2 &\geq& K |m|^{-2\nu_2} \exp(-\alpha_2
|m|^{\beta_2}), \qquad\nu_2>0 \mbox{ if } \alpha_2 =0.
\end{eqnarray*}
Observe that condition (\ref{th3cond}) of Theorem \ref
{th:condisrates} can be
violated only in one of the following ways: $\alpha_1 = \alpha_2 =0$ but
$\nu_2 < \nu_1$, or $\alpha_1 >0$ but $\alpha_2 =0$, or $\alpha
_1>0$ and $ \alpha_2 >0$
but $\beta_2 < \beta_1$.

Applying Theorem \ref{th:lower} with $M=1$, $\varepsilon_n=1$ and $u_1=u_*$,
we arrive
at, as $n \rightarrow\infty$,
\[
R_n^d (B_{p,q}^s (A), u_*, 1) \geq
\cases{
Cn^{-{2s}/({2s+2\nu_1+1})}, &\quad if $\alpha_1=0, \nu_1(2-p) <
p{s^*}$, \cr
C \biggl( \dfrac{\ln n}{n} \biggr)^{{2{s^*}}/({2s^*+2\nu_1})}, &\quad
if $\alpha_1=0, \nu_1(2-p) \geq p{s^*}$, \cr
C (\ln n)^{-{2{s^*}}/{\beta_1}}, &\quad if $\alpha_1>0$.}
\]
On the other hand, applying Theorem \ref{th:upper} with $M=1$,
$\varepsilon
_n=1$ and
$u_1=u^*$, we derive that, as $n \rightarrow\infty$,
\[
\sup_{f \in B_{p,q}^s (A)} {\mathbb E}\|\hat{f}_n^d -f\|^2 \leq
\cases{
Cn^{-{2s}/({2s+2\nu_2+1})} ( \ln n)^{\varrho}, \cr
\hspace*{91pt}\mbox{if
$\alpha_2=0, \nu_2(2-p) < p{s^*}$},
\cr
C \biggl( \dfrac{\ln n}{n} \biggr)^{{2{s^*}}/({2s^*+2\nu_2})}
(\ln n )^{\varrho}, \cr
\hspace*{91pt}\mbox{if $\alpha_2=0, \nu_2(2-p) \geq p{s^*}$},
\cr
C (\ln n)^{-{2{s^*}}/{\beta_2}}, \qquad\mbox{if $\alpha_2>0$},}
\]
where $\rho$ is given by formula (\ref{rovalue}) with $\nu=\nu_2$. Now,
to complete the proof just note that if $\alpha_1 = \alpha_2 =0$ but
$\nu_2 < \nu_1$ or $\alpha_1 >0$ but $\alpha_2 =0$ or $\alpha_1
\alpha_2 >0$
but $\beta_2 < \beta_1$, then the asymptotical minimax lower bounds
for the $L^2$-risk at the point $u=u_*$ are higher than the
corresponding upper bounds at the point $u=u^*$. Hence, in this
case, the convergence rates cannot be independent of the choice of $M$
and the
selection of points $\underline{u}$, arriving at the
required contradiction.
\end{pf*}
\begin{pf*}{Proof of Theorem \protect\ref{th:regular_rates}}
Note that the first inequality in formula (\ref{th41}), as well as
relations (\ref{th42}) and (\ref{th43}) between upper bounds in
discrete and
continuous cases, follow directly from Theorems \ref{th:lower} and
\ref{th:upper} and from inequalities $\tau_1^d (m, u^*,1) \geq K
\tau_1^c (m)$ and $\tau_1^d (m, u^*,1) \geq K \tau_1^d (m,
\underline{u}, M)$.
Hence one only needs to prove the second asymptotic relation in formula
(\ref{th41}).

Let $R_n^d (B_{p,q}^s (A),\underline{u},M)$ be the minimax $L^2$-risk for
fixed values of $\underline{u}$ and $M$, defined by formula
(\ref{risk:point_dis}), and let
\[
H (\underline{u}, M, j)= 2^{-j} \gamma_j^2 \sum_{m \in C_j}\tau
_1^d(m,\underline{u},M).
\]
From the proof of Theorem 1 of Pensky and Sapatinas (\citeyear
{PenskyS09}), it
follows that, in the sparse case [when $ \nu(2-p) \geq p{s^*}$], one
has $R_n^d (B_{p,q}^s (A),\underline{u},M) \geq C 2^{-2 {j_n}{s^*}}$,
as $n
\rightarrow\infty$, where ${j_n}\equiv{j_n}(\underline{u},M)$ is
such that $n
H (\underline{u}, M, j_n)/{j_n}= C$. Similarly, in the dense case
[when $ \nu(2-p) < p{s^*}$], one has $R_n^d (B_{p,q}^s (A),\underline{u},M)
\geq C 2^{-2 {j_n}s}$, as $n \rightarrow\infty$, where ${j_n}\equiv
{j_n}(\underline{u},M)$ is such that $n H (\underline{u}, M, j_n) = C$.

Consider now two different values of $M$, say $M_1$ and $M_2$, and
the corresponding sets of $\underline{u}$'s, say $\underline{u}_1$
and $\underline{u}_2$. If
$\tau_1^d(m,\underline{u}_1,M_1) < \tau_1^d(m,\underline{u}_2,M_2)$
for any $m \in C_j$, then $H (\underline{u}_1, M_1, j) < H (\underline
{u}_2, M_2, j)$.
Observe that, for fixed $M$ and $\underline{u}$, both $H (\underline
{u}, M,
j)$ and $ H (\underline{u}, M, j)/j$ are decreasing functions of
$j$. Hence, if $j_{n1} = {j_n}(\underline{u}_1,M_1)$ and $j_{n2} =
{j_n}(\underline{u}_2,M_2)$ are the values of ${j_n}$ corresponding to
$(\underline{u}_1,
M_1)$ and $(\underline{u}_2, M_2)$, respectively, then $j_{n1} \leq j_{n2}$.
To show that this is true in the dense case, observe that the
opposite, $j_{n1} > j_{n2}$, implies
$
C n^{-1} = H (\underline{u}_1, M_1, j_{n1}) < H (\underline{u}_1,
M_1, j_{n2}) < H (\underline{u}
_2, M_2, j_{n2}),
$
so that $j_{n2}$ cannot be the solution of equation $H (\underline
{u}_2, M_2,
j_{n2}) = C n^{-1}$ and
$j_{n1} > j_{n2}$ cannot be true. In the sparse case,
one just needs to replace $H(\underline{u}, M,j)$ by $ H(\underline
{u}, M,j)/j$.

Now, it follows immediately that in both sparse and dense cases,
$R_n^d (B_{p,q}^s (A)$, $\underline{u}_1,M_1) > R_n^d (B_{p,q}^s
(A),\underline{u}_2,M_2)$. Therefore, the infimum of $R_n^d (B_{p,q}^s
(A),\underline{u},M)$ is attained at $\tilde{M}$ and $\tilde
{\underline{u}}$ such that
$\tau_1^d(m,\tilde{\underline{u}},\tilde{M}) =
\sup_{\underline{u},M} \tau_1^d(m,\underline{u},M)$. Since, for any
choice of $M$ and any selection of points $\underline{u}$,
one has $\tau_1^d(m,\underline{u},M) \leq K \tau_1^{d}(m; u^*, 1)$, the
validity of the theorem follows from Theorem 1 in Pensky and
Sapatinas (\citeyear{PenskyS09}).
\end{pf*}
\begin{pf*}{Proof of Lemma \protect\ref{th:answer}}
Recall that
$
\tau_1^c(m) = \int_{a}^{b}|g_m(u)|^2 \,du$ and
$\tau_1^{d}(m, u^*, 1) = |g_m(u^*)|^2$.
Observe that since $\nu(\cdot)$, $\alpha(\cdot)$ and $\beta(\cdot)$
are continuous functions on the interval $U=[a,b]$, then there exist
$\nu_1 \leq\nu_2$, $\alpha_1 \leq\alpha_2$ and $\beta_1 \leq
\beta_2$
such that $\nu_1 \leq\nu(u) \leq\nu_2$, $\alpha_1 \leq\alpha(u)
\leq\alpha_2$
and $\beta_1 \leq\beta(u) \leq\beta_2$, $u \in[a,b]$. Moreover, in
the inequalities
above, either $\alpha_1 = \alpha_2 = 0$ and $\nu_1 >0$, or $\alpha
_1 >0$ and
$\beta_1 >0$.
Consider the cases when (a) $\alpha(u) \equiv0$ and (b) $\alpha(u) > 0$,
$\beta(u) \neq\mbox{const}$.

\textit{Case} 1: $\alpha(u) \equiv0$. Then $|g_m(u)|^2 \asymp
|m|^{-2\nu(u)}$, so that $|g_m(u)|^2 \leq K |g_m(u^*)|^2$ and
$|g_m(u^*)|^2 \asymp|m|^{-2 \nu(u^*)}$. Hence, in the discrete
case, the asymptotical minimax lower bounds in (\ref{low:reg_equi1})
and the asymptotical minimax upper bounds in (\ref{up:reg_equi1}), for
the $L^2$-risk, follow directly from Theorems \ref{th:lower} and
\ref{th:upper}, respectively.

In order to complete the proof, we need to obtain the asymptotical
minimax lower and upper bounds for the $L^2$-risk in the continuous
model. For this purpose, observe that, under conditions
(\ref{equi:reg}) and (\ref{k_val}), one has [see, e.g., Bender and Orzag
(\citeyear{BO78}), pages 266--267]
%
%
\begin{equation} \label{tauc_polyn}
\tau_1^c(m) \asymp\int_a^b \exp(-2 |{\ln m}| \nu(u)) \,du
\asymp|m|^{-2 \nu(u^*)} ({\ln}|m|)^{-1/k},
\end{equation}
so that Theorems \ref{th:lower} and \ref{th:upper} yield,
respectively, the asymptotical minimax lower bounds in
(\ref{low:reg_equi}) and the asymptotical minimax upper bounds in
(\ref{up:reg_equi}), for the $L^2$-risk.

\textit{Case} 2: $\alpha(u) > 0$ and $\beta(u) \neq\mbox{const}$. In
this case, $\beta(u^*) = \beta_1$. Therefore, one derives
$
K |m|^{-2 \nu_2} \exp(- \alpha_2 |m|^{\beta_1}) \leq|g_m (u^*)|^2
\leq K |m|^{-2 \nu_1} \exp(- \alpha_1 |m|^{\beta_1}).
$
\break Hence the asymptotical minimax lower bounds in (\ref{low:reg_equi1})
and the asymptotical minimax upper bounds in (\ref{up:reg_equi1}), for
the $L^2$-risk, follow directly from
Theorems \ref{th:lower} and \ref{th:upper}, respectively.

To obtain the asymptotical minimax lower and upper bounds for the
$L^2$-risk in the continuous model, note that
%
%
\begin{eqnarray} \label{cont1}
\tau_1^c(m) &\leq& K |m|^{-2 \nu_1} \int_a^b \exp\bigl( -\alpha(u)
|m|^{\beta(u)} \bigr) \,du\nonumber\\[-8pt]\\[-8pt]
&\leq& C_3 (b-a) |m|^{-2 \nu_1}
\exp(-\alpha_1 |m|^{\beta_1}).\nonumber
\end{eqnarray}
On the other hand,
%
%
\begin{equation} \label{cont2}
\tau_1^c(m) \geq K |m|^{-2 \nu_2} \int_a^b \exp\bigl( -\alpha(u)
|m|^{\beta(u)} \bigr) \,du.
\end{equation}
Since $\beta(\cdot)$ is a
continuously differentiable function in some neighborhood of $u^*$,
$|u-u^*|< d$, we have $\beta(u) \leq\beta(u^*) + \beta_* |u-u^*|$,
where $\beta^* = {\max_{|u-u^*|< d}} |\beta^\prime(u)|$. Therefore,
using the inequality $e^z < 1 + 3z$ for $0<z<1$, we obtain
$|m|^{\beta(u)} \leq|m|^{\beta_1} \exp( \beta_* |u-u^*| {\ln}|m| )
\leq
|m|^{\beta_1} (1 + 3 \beta_* |u-u^*| {\ln}|m|)$ for
$|u - u^*| < {\ln}|m|/\break(3 \beta_*)$.
%
Denote $\Omega_m(u^*) = \{ u \in U\dvtx|u-u^*|<|m|^{- (\beta_1+1)}\}
$. Then
%
%
\begin{eqnarray}\label{cont3}
&&\int_a^b \exp\bigl( -\alpha(u) |m|^{\beta(u)} \bigr) \,du \nonumber\\
&&\qquad \geq e^{ - \alpha_2
|m|^{\beta_1}}
\int_{\Omega_m(u^*)} \exp( - \alpha_2 |m|^{\beta_1} 3 \beta_* |u-u^*|
{\ln}|m|) \,du \\
&&\qquad \geq e^{-1} |m|^{- (\beta_1+1)} \exp(-\alpha_2 |m|^{\beta_1}),\nonumber
\end{eqnarray}
since $ {3 \beta_* \alpha_2 |m|^{-1} \ln}|m| <1$ for $|m|$ large enough.
Combining (\ref{cont1})--(\ref{cont3}), we derive that
%
%
\begin{eqnarray} \label{tauc_expo}
&&K e^{-1} |m|^{-(2 \nu_2 +\beta_1 +1)} \exp( -\alpha_2 |m|^{\beta
_1} )\nonumber\\[-8pt]\\[-8pt]
&&\qquad\leq
\tau_1^c(m) \leq K (b-a) |m|^{-2 \nu_1} \exp( -\alpha_1
|m|^{\beta_1}),\nonumber
\end{eqnarray}
so that Theorems \ref{th:lower} and \ref{th:upper} yield,
respectively, the asymptotical minimax lower bounds in
(\ref{low:reg_equi}) and the asymptotical minimax upper bounds in
(\ref{up:reg_equi}), for the $L^2$-risk.
\end{pf*}
\begin{pf*}{Proof of Theorem \protect\ref{th:sufficient}}
First consider the case when $\alpha(u) \equiv0$.
From (\ref{equi:reg}), it follows that
$
\tau_1^d(m,\underline{u},M_n) \leq K |m|^{-2 \nu(u^*)},
$
so that
$\varepsilon_n \leq K ({\ln}|m|)^{\lambda_1}$. Since, in this case,
${\ln}|m| \asymp\ln n$ and ${\ln}|m| > 1$ as $n \rightarrow\infty$,
one has
$\varepsilon_n =O((\ln n)^{\lambda_3})$, where $\lambda_3 = \max
(\lambda_1, 0)$.
The latter, in combination with (\ref{eps_log}),
implies that condition (\ref{ns_prop2}) holds and, moreover, that
$C n (\ln n)^{-\lambda_2} \leq{n^*}\leq C n (\ln n)^{\lambda_3}$.
Then Theorems \ref{th:lower} and \ref{th:upper} imply
that under conditions (\ref{equi:reg}), (\ref{suf_pol}) and (\ref
{eps_log}), one has
$R_n^d (B_{p,q}^s (A), \underline{u}, M) \geq R_n^d (B_{p,q}^s (A))$, where
$R_n^d (B_{p,q}^s (A))$ is given by expression (\ref{low:reg_equi1}) and
that, as \mbox{$n\rightarrow\infty$},
\[
\sup_{f \in B_{p,q}^s (A)} {\mathbb E}\|\hat{f}_n^{d} -f\|^2
\leq
\cases{C ( n^{-1} (\ln n)^{\lambda_1 + \lambda_2} )^{-
{2s}/({2s+2\nu
(u^*)+1})} ( \ln n )^{\varrho},
\vspace*{2pt}\cr
\qquad\mbox{if $\nu(u^*)(2-p) < p{s^*}$},
\vspace*{2pt}\cr
C ( n^{-1} (\ln n)^{1+\lambda_1 + \lambda_2} )^{
{2{s^*}}/({2s^*+2\nu(u^*)})}
( \ln n )^{\varrho}, \vspace*{2pt}\cr
\qquad\mbox{if $\alpha(u) = 0, \nu(u^*)(2-p)
\geq p{s^*}$},}
\]
where $\rho$ is defined in (\ref{rovalue}). If, moreover, (\ref{suf_opp})
holds, then
Theorem \ref{th:lower} yields, as $n\rightarrow\infty$,
\[
R_n^d (B_{p,q}^s (A), \underline{u}, M) \geq
\cases{C ( n^{-1} (\ln n)^{\lambda_1 + \lambda_2} )^{-
{2s}/({2s+2\nu(u^*)+1})},\vspace*{2pt}\cr
\qquad\mbox{if $\nu(u^*)(2-p) < p{s^*}$},
\vspace*{2pt}\cr
C ( n^{-1} (\ln n)^{1+\lambda_1 + \lambda_2} )^{
{2{s^*}}/({2s^*+2\nu(u^*)})}, \vspace*{2pt}\cr
\qquad\mbox{if $\nu(u^*)(2-p) \geq p{s^*}$}.}
\]
To complete the proof of this part, compare the above upper and lower
bounds with (\ref{low:reg_equi}) and (\ref{up:reg_equi}).

Now, let $\alpha(u) >0$. Then, due to assumption (\ref{ns_prop2}) one has
$\ln{n^*}\asymp\ln n$. Under condition (\ref{equi:reg}), by Theorem
\ref
{th:lower},
$R_n^d (B_{p,q}^s (A), \underline{u}, M) \geq R_n^d (B_{p,q}^s (A))
\geq
C (\ln n)^{- 2{s^*}/\beta(u^*)}$, as $n\rightarrow\infty$.
Also, by Theorem \ref{th:upper},
$\sup_{f \in B_{p,q}^s (A)} {\mathbb E}\|\hat{f}_n^{d} -f\|^2 \leq C
(\ln
n)^{- 2{s^*}/\beta(u^*)}$,
as $n\rightarrow\infty$.
To complete the proof, compare the above lower and upper bounds for the
$L^2$-risks with the corresponding bounds in
(\ref{low:reg_equi}) and (\ref{up:reg_equi}).
\end{pf*}
\begin{pf*}{Proof of Theorem \protect\ref{th:necessary}}
Note that conditions (\ref{necc_pol}), (\ref{necc_exp}) and Theorem
\ref
{th:lower} imply that, as $n\rightarrow\infty$,
%
%
\begin{equation}\label{nec_rates}
R_n^d (B_{p,q}^s (A), \underline{u}, M) \geq
\cases{
C ({n^*})^{-{2s}/({2s+2\nu+1})} (\ln{n^*})^{{2s \lambda
}/({2s+2\nu+1})},\cr
\hspace*{92.6pt}\mbox{if $\alpha=0, \nu(2-p) < p{s^*}$},
\cr
C \biggl( \dfrac{\ln{n^*}}{{n^*}} \biggr)^{{2{s^*}}/({2s^*+2\nu})}
(\ln{n^*})^{{2{s^* \lambda}}/({2s^*+2\nu})}, \cr
\hspace*{92.6pt}\mbox{if $\alpha=0, \nu(2-p) \geq p{s^*}$},
\cr
C (\ln{n^*})^{-{2{s^*}}/{\beta}}, \qquad\mbox{if $\alpha>0$}.}
\hspace*{-37pt}
\end{equation}
Denote the ratio between the upper bound for the $L^2$-risk (\ref
{up:reg_equi1})
in the continuous model\vspace*{-2pt} and the lower bound (\ref{nec_rates}) by
$
\Delta_n = \sup_{f \in B_{p,q}^s (A)} {\mathbb E}\|\hat{f}_n^{d*}
-f\|^2 /
R_n^d (B_{p,q}^s (A), \underline{u}, M),
$
and observe that the convergence rates in the discrete model are inferior
to the convergence rates in the continuous model if
$
\lim_{n \rightarrow\infty}(\ln n)^{h {\mathbb I}(\alpha(u) \equiv
0)} \Delta_n=0
$
for any $h>0$.

Let $\alpha(u) \equiv0$ and consider the case when $ \nu(2-p) < p{s^*}$.
Then, taking into account that
under condition (\ref{ns_prop2}) one has $\ln{n^*}\asymp\ln n$, we obtain
\begin{eqnarray*}
\lim_{n \rightarrow\infty}\Delta_n (\ln n)^h & = &
O \Bigl( \lim_{n \rightarrow\infty}n^{-{2s}/({2s+2\nu(u^*)+1})} ( \ln
n )^{\varrho+h}\\
&&\hspace*{34.87pt}{}\times
({n^*})^{{2s}/({2s+2\nu+1})} (\ln{n^*})^{-{2s \lambda
}/({2s+2\nu+1})}
\Bigr)\\
& = & O \Bigl( \lim_{n \rightarrow\infty}\bigl[ (\ln n)^{\rho+ h - {2s
\lambda}/({2s+2\nu+1})}\\
&&\hspace*{38pt}{}\times
n^{- ({2s}/({2s+2\nu(u^*)+1}) - {2s}/({2s+2\nu+1})(1+
\varepsilon_0) )} \bigr] \Bigr).
\end{eqnarray*}
Now, if $\nu> \nu(u^*)$, then it is easy to see that under condition
(\ref{nec10}) we have
$\lim_{n \rightarrow\infty}\Delta_n (\ln n)^h = 0$ for any $h$, and
the convergence rates
in the discrete model are inferior in this case.
If $\nu= \nu(u^*)$, then
$
\lim_{n \rightarrow\infty}\Delta_n (\ln n)^h =
O ( \lim_{n \rightarrow\infty}[ (\ln n)^{\rho+ h - {2s
\lambda}/({2s+2\nu+1})}
(\varepsilon_n)^{ {2s}/({2s+2\nu+1})} ] ) =0
$
if condition (\ref{nec11})\break holds.
The sparse case when $ \nu(2-p) < p{s^*}$ can be treated in a similar manner.

Now, consider the case when $\alpha(u) > 0$. One has
\begin{eqnarray*}
\lim_{n \rightarrow\infty}\Delta_n &=& \lim_{n \rightarrow\infty
}(\ln n)^{-{2{s^*}}/{\beta(u^*)}}
(\ln{n^*})^{ {2{s^*}}/{\beta}} \\
&=& (1+ \varepsilon_0)^{
{2{s^*}}/{\beta}}
\lim_{n \rightarrow\infty}(\ln n)^{-2{s^*}( {1}/{\beta(u^*)} -
{1}/{\beta} ) }
\end{eqnarray*}
and it is easy to see that under each set of conditions in (\ref{nec2}),
$\lim_{n \rightarrow\infty}\Delta_n =0$.
\end{pf*}
\begin{pf*}{Proof of Corollary \protect\ref{cor1}}
Note that if $M=M_n$ is finite, then for \mbox{$\alpha(u) \equiv0$} one has
$\tau_1^d(m,\underline{u},M_n) \asymp|m|^{-2 \nu}$ where
$\nu= \min(\nu(u_1), \nu(u_2),\ldots, \nu(u_M))$. If $\alpha(u)>0$,
then denote
$l_0 = \arg\min_l \beta(u_l)$, $\beta= \beta(u_{l_0})$, $\nu_0 =
\nu
(u_{l_0})$ and $\alpha_0= \alpha(u_{l_0})$.
In this case,
$
\tau_1^d(m,\underline{u},M_n) \asymp|m|^{-2 \nu_0 } \exp( - \alpha_0
|m|^\beta)
$
and hence the validity of the corollary follows from Theorems \ref
{th:sufficient} and \ref{th:necessary}.
\end{pf*}
\begin{pf*}{Proof of Corollary \protect\ref{cor2}}
Note that
$
\tau_1^d(m,\underline{u},M_n) \geq K (\ln n)^{-\lambda^*} |m|^{-2
\nu(u^*)},
$
and hence the validity of the corollary follows from Theorem \ref{th:sufficient}.
\end{pf*}
\begin{pf*}{Proof of Corollary \protect\ref{cor3}}
Note that, for $u_l$ such that $\beta(u_l) = \beta(u^*)$ one has
$\tau_1^d(m,\underline{u},M_n) \geq K n^{-\tau} |m|^{-2 \nu(u_l)}
\exp(
- \alpha(u_l) |m|^{\beta(u^*)} )$.
Then the validity of the corollary follows from
Theorem \ref{th:sufficient}.
\end{pf*}
\begin{pf*}{Proof of Theorem \protect\ref{th:uniform_sample}}
%
First, consider the case when $\alpha(u) \equiv0$. Denote $v(x)= \nu
(S(x))$, $x^* = q(u^*)$,
and let $l^*$ be the index of a point closest to $x^*$, that is,
$l^* = \arg\min|x^* - (l-1+d)/M|$. Note that $v(x^*) = \nu(u^*)$ and
the function $v(x)$ is continuously differentiable
with $|v^\prime(x)| \leq v_0$ for some $v_0 >0$. Note that
$
\tau_1^d(m,\underline{u},M_n) \leq K |m|^{-2 \nu(u^*)},
$
so if we show that
under condition (\ref{uniform_cond}) we have
%
%
\begin{equation} \label{taudi}
\tau_1^d(m,\underline{u},M_n) \geq K |m|^{-2 \nu(u^*)} (\ln
n)^{-\lambda}
\end{equation}
for some constant $\lambda\in\mathbb R$, then the
validity of the theorem will follow from Theorem~\ref{th:sufficient}.
In order to prove (\ref{taudi}), note that
\begin{eqnarray*}
\tau_1^d(m,\underline{u},M )
&\asymp& \frac{1}{M} \sum_{l=1}^{M }
|m|^{-2 v( {l-1+d})/{M }}\\
&\asymp& \frac{1}{M} \sum_{l=1}^{M } |m|^{-2 [ v ({l-1+d})/{M
} - v ({l^*-1+d})/{M } ]
-2 [ v ({l^*-1+d})/{M } - v (x^*) ] - 2 \nu(u^*) }\\
&\geq& \frac{K}{M} |m|^{-2 \nu(u^*) } \sum_{l=1}^{M } |m|^{-
{2v_0 |l-l^*|}/{M} - {v_0^*}/{M} }\\
&\geq& \frac{K}{M} |m|^{-2 \nu(u^*) } \sum_{k=0}^{M/2-1} |m|^{-
{2v_0 k}/{M} - {v_0^*}/{M} },
\end{eqnarray*}
where $v_0^* = v_0 {\mathbb I}(x^* \neq(l^*-1+d)/M)$.
Now, recall the following statement from Calculus: if $u(z)$, $z \geq
0$, is a continuous,
positive, monotonically decreasing function, then
%
%
\begin{eqnarray} \label{Calc1}
\sum_{k=0}^{M/2-1} u(k) &\geq& \max\biggl( \int_0^{M/2} u(x) \,dx, u(0) + \int
_1^{M/2} u(x) \,dx \biggr) \nonumber\\[-8pt]\\[-8pt]
&\geq&
\frac{1}{2} \biggl( u(0) + \int_0^{M/2} u(x) \,dx \biggr).\nonumber
\end{eqnarray}
Applying (\ref{Calc1}) with $u(x)= |m|^{-2 \nu(u^*) -
{v_0^*}/{M} } |m|^{- {2v_0 x}/{M}}$,
and taking into account that
$
\int_0^{M/2} u(x) \,dx \asymp|m|^{-2 \nu(u^*) - {v_0^*}/{M} }
M^{-1} {\ln}|m|,
$
we obtain
\[
\tau_1^d(m,\underline{u},M ) \geq K |m|^{-2 \nu(u^*)} ({\ln}
|m|)^{-1} ( {1+ M^{-1} \ln}|m| ) \times
\exp( {- v_0^* M^{-1} \ln}|m| ).
\]
Now, recall that ${\ln}|m| \asymp\ln n$ in this case and note that
under the first assumption in (\ref{uniform_cond}),
$\tau_1^d(m,\underline{u},M )$ satisfies condition
(\ref{suf_pol}) of Theorem \ref{th:sufficient} with $\lambda_1 = -1$
and $\varepsilon
_n=1$. Hence, the convergence rates
in the discrete and the continuous models almost coincide.

Now, consider the case when $\alpha(u)> 0$. Denote $v(x)= \beta(S(x))$
and let, as before, $x^* = q(u^*)$
and $l^* = {\arg\min}|x^* - (l-1+d)/M|$. Note that $v(x^*) = \beta
(u^*)$ and that the function $v(x)$ is continuously differentiable
with $|v^\prime(x)| \leq v_0$ for some constant $v_0 >0$. Denote, as
before, $v_0^* = v_0 {\mathbb I}(x^* \neq
(l^*-1+d)/M)$.
Note that
$
\tau_1^d(m,\underline{u},M ) \geq M^{-1} K |m|^{-2\nu_1} \sum
_{l=1}^M \exp( - \alpha_1 |m|^{\beta(u_l)} ),
$
where $\nu_1 = \max\nu(u)$, $\alpha_1 = \max\alpha(u)$, $u \in U$,
and, in order to prove the statement, we need to construct a lower
bound for
$S(m,M) = \sum_{l=1}^M \exp( - \alpha_1 |m|^{\beta(u_l)} )$. Similarly
to the polynomial case, we obtain that
$
S(m,M) \geq\break K \sum_{k=0}^{M/2-1} \exp( - \alpha_1 |m|^{\beta(u^*)+
{v_0 k}/{M} + {v_0^*}/({2 M}) } ).
$
Denote $\alpha_m = \alpha_1 |m|^{\beta(u^*)+ {v_0^*}/({2 M})}$, and
apply inequality (\ref{Calc1}) with
$u(x) = \exp( - \alpha_m |m|^{{v_0 k}/{M}} )$.
Observe that
\[
\int_0^{M/2} u(x) \,dx = \frac{M}{{v_0^* \ln}|m|} \int_1^{|m|^{v_0/2}}
z^{-1} \exp(-\alpha_m z)\,dz
\geq\frac{M}{{2 v_0^* \ln}|m| \alpha_m} \exp(-\alpha_m)
\]
and recall that ${\ln}|m| \asymp\ln\ln n$. Hence, under the second of
conditions in (\ref{uniform_cond}),
as $M \rightarrow\infty$ and $|m| \rightarrow\infty$, we derive that
\[
S(m,M)
\geq
K M (\ln n)^{-1} |m|^{- (\beta(u^*)+1) } \exp\bigl( - \alpha_1
|m|^{\beta(u^*)} \exp[{0.5 v_0^* M^{-1} \ln}|m|] \bigr).
\]
Note that due to assumption (\ref{uniform_cond}) and due to ${\ln}|m|
\asymp\ln\ln n$, there exists $\tau_3>0$ such that
${M_n^{-1} \ln}|m| \leq\tau_3$ when $n$ is large enough. Therefore,
$
\tau_1^d(m,\underline{u},M ) \geq C (\ln n)^{-1} |m|^{- (2\nu_1+
\beta(u^*)+1) }
\exp( - \alpha_1 \exp(0.5 v_0^* \tau_3 ) |m|^{\beta(u^*) } ).
$
Application of Theorem \ref{th:sufficient} with $\nu= \nu_1+ 0.5
\beta(u^*)+ 0.5$,
$\alpha= \alpha_1 \exp(0.5 v_0^* \tau_3 )$ and $\varepsilon_n =
(\ln
n)^{-1}$ completes the proof of this
part of the statement.

To prove the last statement in Theorem \ref{th:uniform_sample}, recall
a simple fact from Calculus: if function $F(x)$, $x \in[0,1]$, is
continuously differentiable with $F_0 = {\max_x }| F^\prime(x) |$, then
for any $d$
such that $0 \leq d \leq1$ one has
%
%
\begin{equation} \label{Calc}
\Biggl| M^{-1} \sum_{l=1}^M F \bigl( (l-1+d)/M \bigr) - \int_0^1 F(x)\, dx \Biggr|
\leq0.5 M^{-1} F_0.
\end{equation}
Let $\alpha(u) \equiv0$. Note that since $| q'(u)|$ is bounded and
separated from zero, one has
$
{\int_0^1} |g_m(S(x))|^2 \,dx = \int_a^b |g_m(u)|^2 q'(u) \,du \asymp\tau_1^c(m).
$
Therefore, if
%
%
\begin{equation} \label{residual}
R(m, n) = \tau_1^d(m,\underline{u},M_n) - \int_0^1 |g_m(S(x))|^2 \,dx
= o
( \tau_1^c(m) )
\end{equation}
as $n \rightarrow\infty$, then $\tau_1^d(m,\underline{u},M_n)
\asymp\tau_1^c(m)$
and the theorem is proved.
Applying formula (\ref{Calc}) to $F(x) = |g_m(S(x))|^2$ and noting that
$|S' (x)| \leq s_2$, we obtain
\[
R(m, n) \leq0.5 s_2 M_n^{-1} \max_{u \in[a,b]} \biggl| \frac
{d}{du} |g_m(u)|^2 \biggr|
= O \bigl( M_n^{-1} {|m|^{-2\nu(u^*)} \ln}|m| \bigr).
\]
Comparing the last expression with $\tau_1^c(m)$ given by formula
(\ref{tauc_polyn}), we confirm
that condition (\ref{residual}) holds and the theorem is valid in this case.
\end{pf*}
\begin{pf*}{Proof of Lemma \protect\ref{prop:case1}}
Recall that $0 < a < b <\infty$, $\beta_1 \leq\beta(u) \leq
\beta_2$, $u \in[a,b]$, for some $0 < \beta_1 \leq\beta_2 <
\infty$, and $u_l=a+(b-a)l/M$, $l=1,2,\ldots,M$. Consider first the
case when $a \in\mathbb N$. Then,
$
4\pi^2 m^2 M \tau_1^d(m,\underline{u},M) =\break \sum_{l=1}^M \beta
^2(u)\sin
^2(2 \pi\times m u_l) \geq
\beta^2_1 \sum_{l=1}^M \sin^2 (2\pi m(b-a) l/M ).
$
Using the formula 1.351.1 of Gradshtein and Ryzhik (\citeyear{GR80}) with
$x=2\pi m (b-a)/M$, $n=M$ and $k=l$, we obtain
%
%
\begin{eqnarray}
\label{eq:box-car-bound1}
&&\tau_1^d(m,\underline{u},M)
\geq
\frac{\beta^2_2}{4\pi^2 m^2 M}
\biggl( \frac{M}{2}
-\bigl(\cos\bigl(2(M+1)\pi m (b-a)/M\bigr)\nonumber\\
&&\hspace*{191.4pt}{}\times \sin\bigl(2 \pi m (b-a)\bigr)\bigr)\\
&&\hspace*{153.28pt}{}\times\bigl({2\sin\bigl(2 \pi m
(b-a)/M\bigr)}\bigr)^{-1} \biggr).\hspace*{-25pt}\nonumber
\end{eqnarray}
Since $2^j \leq|m|< 2\pi/3 2^j$ and $\ln n \leq2^j \leq n^{1/3}$,
the condition $M_n \geq(32 \pi/3)(b-a)n^{1/3}$ guarantees that $|2
\pi m (b-a)/M| \leq\pi/2$.
In this case, using the inequality $y \leq2 \sin(y)$, $0 \leq y
\leq\pi/2$ [see, e.g., Lang (\citeyear{L66}), page 41], we derive
%
%
\begin{equation} \label{eq:box-car-bound2}\qquad\quad
2\sin\bigl(2 \pi m (b-a)/M\bigr) \geq2 \pi m (b-a)/M \geq\bigl(4 \pi^2 (b-a) \ln n\bigr)/(3M).
\end{equation}
%
Hence, combining (\ref{eq:box-car-bound1}) and (\ref{eq:box-car-bound2}),
for $n$ large enough, we arrive at $\tau_1^d(m,\underline{u}$,\break $M) \geq
K m^{-2}$.

Consider now the case when $a \notin\mathbb N$. A standard
trigonometrical identity yields
%
%
\begin{eqnarray}
\label{eq:box-car-bound3} &&\sum_{l=1}^M \sin^2\bigl(2 \pi m a + 2 \pi m
(b-a) l/M \bigr) \nonumber\\[-8pt]\\[-8pt]
&&\qquad= \frac{1}{2} \Biggl(M -\sum_{l=1}^M \cos\bigl(4 \pi m a + 4
\pi m (b-a) l/M \bigr)\Biggr).\nonumber
\end{eqnarray}
Using formula 1.341.3 of Gradshtein and Ryzhik (\citeyear{GR80})
with $x=4
\pi m a$, $y=4 \pi m(b-a)/M$, $n=M$ and $k=l$, we derive for $M
\geq4 (b-a)|m|$,
%
%
\begin{eqnarray}\qquad
\label{eq:box-car-bound4}
&&\sum_{l=1}^M \cos\bigl(4 \pi m a +
4 \pi m (b-a) l/M \bigr) \nonumber\\
&&\qquad = \frac{\cos(4\pi m a + 2\pi m (b-a) (M-1)/M)
\sin
(2\pi m(b-a))}{\sin(2\pi m (b-a)/M)} \\
&&\qquad \leq\frac{M}{\pi m (b-a)}.\nonumber
\end{eqnarray}
Hence combining (\ref{eq:box-car-bound3}) and (\ref{eq:box-car-bound4})
in a manner similar to the first part of the proof, for $|m|$ large
enough, we arrive at $\tau_1^d(m,\underline{u},M) \geq K m^{-2}$
which completes the proof.
\end{pf*}
\begin{pf*}{Proof of Theorem \protect\ref{th:box-car-minimax}}
The proof follows directly from the discussion of Section
\ref{box-car}, by combining Theorems \ref{th:lower},
\ref{th:upper} and Lemma\vspace*{1pt} \ref{prop:case1}, taking $A_j = C_j =\{
m\dvtx
\psi_{mjk} \neq0\}$, and noting that, for the Meyer wavelets, $C_j
\subseteq2\pi/3 [-2^{j+2}$,\break $-2^j] \cup[2^j, 2^{j+2}]$ with $|C_j| =
4\pi2^j$ [see Johnstone et al. (\citeyear{Johnstoneetal04}),
page 565].
\end{pf*}
\end{appendix}

\section*{Acknowledgments}

The authors would like to thank Athanasia Petsa for her help in
carrying out
the simulation study. Finally, we would like to thank an Associate
Editor and two anonymous referees for their suggestions
which helped to significantly improve the paper.

\printaddresses

\end{document}